\documentclass[reqno]{amsart}
\usepackage{amssymb}
\usepackage{mathrsfs}
\usepackage{amsmath,amstext,amsfonts,verbatim}

\usepackage[latin1]{inputenc}

\usepackage{amsmath,amsthm,amssymb}

\usepackage{amsfonts}

\usepackage{amsmath,amssymb}
\usepackage{graphicx}
\usepackage{color}
\usepackage{indentfirst}

\theoremstyle{definition}

\theoremstyle{remark}

\numberwithin{equation}{section}

\newcommand{\neweq}[1]{\begin{equation}\label{#1}}

\def\phi{\varphi}

\def\incep{\left\{\begin{array}{cl} }
 \def\termin{\end{array}\right. }
\def\2af{2^*_\alpha}

\begin{document}

\title [*******]{\textbf{ A note on weak mean equicontinuity and strong mean sensitivity}}

\author{Zhongxuan Yang$^*$}
\address{College of Mathematics and Statistics, Chongqing University, Chongqing 401331, China}
\thanks{$^*$ Corresponding author.}
\email{20220601008@stu.cqu.edu.cn}
\thanks{The research was supported by NSF of China (No. 11671057) and NSF of Chongqing (Grant No. cstc2020jcyj-msxmX0694).}

\author{Xiaojun Huang}
\address{College of Mathematics and Statistics, Chongqing University, Chongqing 401331, China}
\email{hxj@cqu.edu.cn}


% General info

\date{}

\begin{abstract}
In this paper, we study the weak mean metric and give some properties by replacing the Besicovitch pseudometric with weak mean metric in the definition of mean equicontinuity and mean sensitivity. We study an opposite side of weak mean equicontinuity, strong mean sensitivity and we obtain a version of Auslander-Yorke dichotomies: minimal topological dynamical systems are either weak mean equicontinuous or strong mean sensitive, and transitive topological dynamical systemss are either almost weak mean equicontinuous  or strong mean sensitive. Furthermore, motivated by the localized idea of sensitivity, we  introduce some notions of new version sensitive tuples and  study the properties of these sensitive tuples, we show that  a transitive dynamical system is strong mean sensitive  if and only if it admits a strong mean sensitive tuple. Finally, We introduce the notions of  weakly mean equicontinuity  of a topological dynamical system respect to a given continuous function $f$, and  we show that a topological dynamical system is weakly mean equicontinuity then it is weakly mean equicontinuity with respect to every continuous function. 
 
\end{abstract}

\maketitle

\section{Introduction}

Throughout this paper, a topological dynamical system is a pair  $(X, T)$, where  $X$  is a nonempty compact metric space with a metric  $d $ and $ T$  is a continuous map from  $X $ to itself.

Equicontinuous systems are renowned for their straightforward dynamical behaviors. A topological dynamical system  $(X, T)$  is called equicontinuous if for any  $\varepsilon>0 $ there exists  $\delta>0$  such that  $d\left(T^{n} x, T^{n} y\right)<\varepsilon$  for any  $n \in \mathbb{N} $ whenever  $d(x, y)<\delta $. The examination of topological dynamical systems whose structures exhibit proximity to equicontinuous ones is a natural pursuit.  To address this, Fomin [2] introduced a weaker form of equicontinuity known as mean-$L$-stable. Recall that a topological dynamical system  $(X, T)$  is mean-$L$-stable [3] if for any  $\epsilon>0 $, there is a  $\delta>0$  such that if  $x, y \in X$  with  $d(x, y)<\delta $ then  $d\left(T^{k} x, T^{k} y\right)<\epsilon$  for all $ k \in \mathbb{Z}_{+}$ except a set of upper density no larger than  $\epsilon $, and Fomin established the noteworthy result that if a minimal system is mean-L-stable then it is uniquely ergodic. This distinctive property aligns with the concept of mean equicontinuity introduced by Li et al. [1]. The seminal work by Li et al. demonstrates that mean equicontinuity is stronger than just diserete spectrum. We say that a topological dynamical system  $(X, T)$  is mean equicontinuous if for any  $\epsilon>0$  there is a  $\delta>0$  such that if  $x, y \in X$  with $ d(x, y)<\delta$  then
$$
\limsup _{n \rightarrow \infty} \frac{1}{n} \sum_{i=0}^{n-1} d\left(T^{i} x, T^{i} y\right)<\epsilon .
$$

The phenomenon of mean equicontinuity has garnered heightened scholarly interest, emerging as a pivotal area within recent years. Its profound connections with the ergodic properties intrinsic to measurable dynamical systems underscore its significance. For a thorough exploration of the intricacies surrounding mean equicontinuity and its correlated domains, scholars are encouraged to delve into the detailed survey outlined in [1]. This survey provides a nuanced elucidation of mean equicontinuity, complemented by additional characterizations and insights available in a multitude of references [4-16]. The dynamics of a mean equicontinuous topological dynamical system maintain simplicity, characterized by a discrete spectrum [1] and bounded "max-mean" complexity [13]. Notably, in the study of topological dynamical systems with bounded mean complexity, Huang et al. [13] introduced the notion of equicontinuity in the mean, and established that a minimal topological dynamical system is mean equicontinuous if and only if it is equicontinuous in the mean. This characterization holds true for general topological dynamical systems [5]. 

To provide a natural localization of equicontinuity, Akin et al. [20] introduced the concept of almost equicontinuity. Notably, an almost equicontinuous topological dynamical system exhibits simplicity in the sense of zero entropy [21]. Extending this idea, Li et al. [1] introduced the notion of almost mean equicontinuity, revealing that the dynamics of some almost mean equicontinuous topological dynamical systems can become intricate, showcasing positive entropy. As an illustration, García-Ramos et al. [22] presented an example of an almost mean equicontinuous topological dynamical system that is Devaney chaotic and exhibits positive entropy.

In addition, the concept of sensitive dependence on initial conditions, often referred to as sensitivity, constitutes a foundational and pivotal element within chaos theory and the broader field of dynamical systems theory. It is pertinent to recall that within the framework of a topological dynamical system $(X, T)$, sensitivity is established when there exists $\delta > 0$ such that for any nonempty open subset $U$ of $X$, there exist $x, y \in U$ and $k\in\mathbb{N}$ satisfying $ d(T^{k} x, T^{k} y)>\delta$. This dynamical system manifests dichotomous properties depending on its minimality or transitivity. Specifically, if $(X, T)$ is minimal, it inherently assumes either equicontinuity or sensitivity [24]. In the case of transitivity, $(X, T)$ exhibits either almost equicontinuity or sensitivity [20]. Expanding our considerations to the realm of mean sensitivity, Li et al. [1] introduced the concept and demonstrated its compatibility with the aforementioned dichotomy properties, maintaining consistency with (almost) mean equicontinuity. Furthermore, Li et al. [8] introduced the notion of sensitivity in the mean, providing evidence that the established dichotomy properties persist concerning sensitivity in the mean and (almost) equicontinuity in the mean, and they showed that for a minimal topological dynamical system is mean sensitive if and only if it is sensitive in the mean. In this paper, we will show a  topological dynamical system is mean sensitive if and only if it is sensitive in the mean [see Appendix, Theorem 1].

Blanchard et al. [24-26] pioneered the introduction of the concepts of entropy pairs and entropy pairs concerning a measure $\mu$. Subsequently, these notions have undergone comprehensive studied by numerous researchers. Building upon these foundational ideas, Huang et al. [27] extended the discourse by proposing the concepts of entropy tuples and entropy tuples associated with a measure $\mu$. Motivated by the intricate theory surrounding  entropy tuples, Ye and Zhang [28] introduced the concept of a sensitive tuple. Specifically, they established that a transitive topological dynamical system demonstrates sensitivity if and only if it possesses a sensitive tuple. Furthermore, it was revealed that an entropy tuple of a transitive topological dynamical system is a sensitive tuple. Drawing inspiration from these fundamental considerations, Li et al. [8] introduced the notions of mean-sensitive tuple and sensitive in the mean tuple. They demonstrated that a transitive topological dynamical system attains mean sensitivity (sensitivity in the mean) if and only if it accommodates a mean-sensitive tuple (sensitive in the mean tuple). Yu [29] and Li [9] delved deeper into the study of mean sensitivity and its localized representation, namely, weakly 2-sensitive in the mean (density sensitivity) and Weakly sensitive in the mean tuple (density-sensitive tuple). Their findings established that a transitive dynamical system exhibits weak 2-sensitivity in the mean (density sensitivity) if and only if it possesses a weakly sensitive in the mean tuple (density-sensitive tuple).

In the seminal work by García-Ramos and Kwietniak [18], the fundamental concepts of FK-sensitivity and FK-continuity were introduced to characterize models of zero entropy loosely Bernoulli systems, representing a paradigm shift by replacing the Besicovitch pseudometric with Feldman-Katok pseudometric in the definitions of mean sensitivity and mean equicontinuity. This adaptation not only led to the establishment of Feldman-Katok continuity but also unveiled the innovative notion of FK-sensitivity. The authors rigorously demonstrated that, for a minimal dynamical system, the inherent properties gravitate towards either Feldman-Katok continuity or FK-sensitivity. It is noteworthy that Feldman-Katok continuity, as juxtaposed with mean equicontinuity, manifests itself with a discernible level of weakened constraints. Furthermore, in the meticulous examination of the Feldman-Katok pseudometric, intriguing allowances for 'time delay' and 'space jump' are explicitly accommodated.

In a recent scholarly investigation conducted by Zheng et al. [12], the nuanced aspects pertaining to the temporal sequencing of points within orbits were deliberately excluded. This strategic omission is predicated upon the realization that considerations of such temporal order lack substantive significance in the exploration of statistical properties across protracted temporal durations. Rather than accord priority to the order of points, the authors introduced a novel mathematical construct termed the weak-mean pseudometric $\overline{F}$, which is analogous to the pseudometric $\widetilde{d}_{FK}$ defined in [17]. This introduction brought forth the concept of weak mean equicontinuity.  We say  $(X, T)$ is weak mean equicontinuous if  for any  $\epsilon>0$  there is a  $\delta>0$  such that if  $x, y \in X$  with $ d(x, y)<\delta$  then  $\overline{F}(x, y)  =\limsup _{n \rightarrow+\infty} \inf _{\sigma \in S_{n}} \frac{1}{n} \sum_{k=1}^{n} d\left(T^{k} x, T^{\sigma(k)} y\right)<\varepsilon.$
Consequently, it can be deduced that weak mean equicontinuity represents a weaker concept  compared to Feldman-Katok continuity. In their seminal contribution [12], the authors conclusively established that systems exhibiting weak mean equicontinuity precisely align with those wherein the time-average operator adeptly preserves the continuity of observed.

In the seminal work [11], conducted to delve into the characterization of discrete spectrum, García-Ramos et al. introduced the concept of mean equicontinuity with respect to a function. Let $(X, T)$ be a dynamical system, and consider a function $f: X \rightarrow \mathbb{R}$. The system $(X, T)$ is mean equicontinuous with respect to $f$, if for any $\varepsilon > 0$, there exists $\delta > 0$ such that for any $x, y \in X$ satisfying $d(x, y) \leq \delta$, it holds that $\limsup_{n \to \infty} \frac{1}{n} \sum_{i=1}^{n} \left| f(T^{i} x) - f(T^{i} y) \right| \leq \varepsilon$. Remarkably, García-Ramos and collaborators established a profound result, demonstrating that a dynamical system $(X, T)$ achieves mean equicontinuity if and only if it is mean equicontinuous with respect to every continuous function defined on $X$.

 This paper is organized as follows. In Section 2, we give some basic notions 
 and results of topological dynamical system. In Section 3, we  study some properties of  weakly mean equicontinuity, and we show that for a transitive weakly mean equicontinuous system, it is minimal if and only if it has full support. In Section 4, we study the localization of weak mean equicontinuity, namely, almost weakly  mean equicontinuous (resp.  almost weakly  equicontinuous in the mean ), and we show that for a almost weakly  mean equicontinuous (resp.  almost weakly  equicontinuous in the mean ) system, then every transitive point is weakly mean equicontimuous (resp. weakly  equicontinuous in the mean ). In Section 5, we consider  the opposite side of equicontinuity, namely,  strong  sensitivity in the mean and strong mean equicontinuity, and we show that  a minmimal system is either weakly  mean equicontinuous (resp. weakly  equicontinuous in the mean )  and strong mean sensitive (resp.  strong  sensitive in the mean) and for a dynamical system is strong mean sensitive if and only if it is strong  sensitive in the mean without the need for a   minimality condition. In Section 6, we introduce the notions of strong sensitive tuple, strong sensitive in the mean tuple etc. , and we show that a transitive topological dynamical system attains strong mean sensitivity (strong sensitivity in the mean) if and only if it accommodates a  strong mean sensitive tuple (strong sensitive in the mean tuple). In Section 7, we introduce the notions of weakly density-$t$-equicontinuity and density $\overline{F}$-sensitivity, and we show that a topological dynamical system  is weakly mean  equicontinuous if and only if it is weakly density-$t$-equicontinuous for every  $t \in[0,1)$. Finally, we consider weakly mean equicontinuity and weakly equicontinuity in the mean  with respect to a continuous function, and  we show that a topological dynamical system is weakly mean equicontinuity then the system is weakly mean equicontinuity with respect to every continuous function.

\section{Preliminaries}

In this section we recall some notions and results of topological dynamical system.  which are needed in our paper. Note that  $\mathbb{N}$  denotes the set of all non-negative integers and  $\mathbb{N}^{+}$ denotes the set of all positive integers in this paper.

{\bf 2.1} Let  $F \subset \mathbb{N}$, we define the upper density  $\overline{D}(F)$  of  $F$  by
$$
\overline{D}(F)=\limsup _{n \rightarrow+\infty} \frac{\#(F \cap[0, n-1])}{n},
$$
where  $\#(\cdot)$  is the number of elements of a set. Similarly, the lower density  $\underline{D}(F)$  of $ F$  is defined by
$$
\underline{D}(F)=\liminf _{n \rightarrow+\infty} \frac{\#(F \cap[0, n-1])}{n} .
$$
We call  $ F$  has density $ D(F)$  if  $\overline{D}(F)=\underline{D}(F) $.

{\bf 2.2} Suppose  $(X, T)$  is a topological dynamical system. The  $\sigma$ -algebra of Borel subsets of  $X$  will be denoted According to  $\mathscr{B}(X) $. Let  $M(X) $ be the collection of all regular Borel probability measures defined on the measurable space  $(X, \mathscr{B}(X))$. In the weakly* topology, $ M(X)$  is a nonempty compact set.

We say  $\mu \in M(X)$  is  $T$-invariant if  $\mu\left(T^{-1}(A)\right)=\mu(A)$  holds for any  $A \in \mathscr{B}(X) $. Denote According to $ M(X, T)$  the collection of all $ T $-invarant regular Borel probability measures defined on the measurable space  $(X, \mathscr{B}(X))$. In the weakly  ${ }^{*} $ topology,  $M(X, T)$  is a nonempty compact convex set.

We say  $\mu \in M(X, T)$  is ergodic if for any  $A \in \mathscr{B}(X)$  with  $T^{-1} A=A$, $ \mu(A)=0  $ or $ \mu(A)=1 $ holds. Denote by  $E(X, T)$  the collection of all ergodic measures on $ (X, T)$. As well known, $ E(X, T)$  is the collection of all extreme points of  $M(X, T)$ and  $E(X, T)$  is nonempty. 

We say  $(X, T)$  is uniquely ergodic if  $E(X, T)$  is singleton. Since  $E(X, T)$  is the set of extreme points of  $M(X, T)$, then  $(X, T)$  is unique ergodic if and only if  $M(X, T) $ is singleton. For  $\mu \in M(X, T)$, the support of  $\mu$  is defined by  $\operatorname{supp}(\mu)=\{x \in X: \mu(U)>  0 \text{\ for any neighborhood\ }  U  \text{\ of\ }  x\} $. A measure  $\mu$  on  $X$  has full support if  $\operatorname{supp}(\mu)=X $. It is known that the support of an ergodic measure is a transitive subsystem. The support of a topological dynamical system  $(X, T)$  is defined by  $\operatorname{supp}(X, T)=\overline {\bigcup \{\operatorname{supp}(\mu): \mu \in M(X, T)\}} $.

Given  $x \in X $, we have  $\left\{\frac{1}{n} \sum_{k=0}^{n-1} \delta_{T^{k} x}\right\}_{n\in\mathbb{N}} \subset M(X)$, where $ \delta_{x} $ is the Dirac measure supported on $x$. Denote by $ M_{x}$  the collection of all limit points of  $\left\{\frac{1}{n} \sum_{k=0}^{n-1} \delta_{T^{k} x}\right\}_{n\in\mathbb{N}} .$  Since  $M(X) $ is compact, we have  $M_{x} \neq \emptyset$. Moreover, $ M_{B} \subset M(X, T)$. We call  $M_{x}$  the measure set generated by  $x$.

A point $ x \in X $ is called  generic point (see [30]) if for any  $f \in C(X)$, the time average
$$
f^{*}(x)=\lim _{n \rightarrow+\infty} \frac{1}{n} \sum_{k=1}^{n} f\left(T^{k} x\right)
$$
exist.

The following Lemmas are important throughout the entire article.

{\bf Lemma 2.1}(see [12]) Let  $(X, T)$  be a topological dynamical system. Then

(1) For any sequences  $\left\{x_{k}\right\}_{k=0}^{n-1}$  and  $\left\{y_{k}\right\}_{k=0}^{n-1}$  of $ X$, we have
$$
\inf _{\sigma \in S_{n}} \sum_{k=0}^{n-1} d\left(x_{k}, y_{\sigma(k)}\right)=\inf _{\sigma \in S_{n}} \sum_{k=0}^{n-1} d\left(y_{k}, x_{\sigma(k)}\right)
$$
In particular, for any  $x, y \in X$, we have
$$
\inf _{\sigma \in S_{n}} \frac{1}{n} \sum_{k=0}^{n-1} d\left(T^{k} x, T^{\sigma(k)} y\right)=\inf _{\sigma \in S_{n}} \frac{1}{n} \sum_{k=0}^{n-1} d\left(T^{k} y, T^{\sigma(k)} x\right) .
$$

(2) For any sequences  $\left\{x_{k}\right\}_{k=0}^{n-1},\left\{y_{k}\right\}_{k=0}^{n-1}$  and $\left\{z_{k}\right\}_{k=0}^{n-1}$  of  $X$ , we have
$$
\inf _{\sigma \in S_{n}} \sum_{k=0}^{n-1} d\left(x_{k}, z_{\sigma(k)}\right) \leq \inf _{\sigma \in S_{n}} \sum_{k=0}^{n-1} d\left(x_{k}, y_{\sigma(k)}\right)+\inf _{\sigma \in S_{n}} \sum_{k=0}^{n-1} d\left(y_{k}, z_{\sigma(k)}\right).$$
In particular, for any  $x, y, z \in X$ , we have
$$
\inf _{\sigma \in S_{n}} \frac{1}{n-m} \sum_{k=0}^{n-1} d\left(T^{k} x, T^{\sigma(k)} z\right) \leq \inf _{\sigma \in S_{n}} \frac{1}{n} \sum_{k=0}^{n-1} d\left(T^{k} x, T^{\sigma(k)} y\right)+\inf _{\sigma \in S_{n}} \frac{1}{n} \sum_{k=0}^{n-1} d\left(T^{k} y, T^{\sigma(k)} z\right).
$$

(3) For any  $x, y \in X$, we have
$$
\overline{F}(x, y)=\overline{F}(y, x)
$$
and
$$
\underline{F}(x, y)=\underline{F}(y, x).
$$

(4) For any  $x, y, z \in X$ , we have
$$
\overline{F}(x, z) \leq \overline{F}(x, y)+\overline{F}(y, z)
$$
and
$$
\underline{F}(x, z) \leq \underline{F}(x, y)+\overline{F}(y, z).
$$

{\bf Lemma 2.2}(see [12]) Let  $(X, T)$  be a topological dynamical system. For any  $x, y \in X$  and  $r, s \in \mathbb{N}$, we have
$$
\overline{F}(T^{r} x, T^{s} y)=\overline{F}(x, y)
$$
and
$$
\underline{F}\left(T^{r} x, T^{s} y\right)=\underline{F}(x, y)
$$
If $F(x, y)$  exists, we also have
$$
F\left(T^{r} x, T^{s} y\right)=F(x, y).
$$

\section{weakly mean equicontinuity}
In this section,  we  will  study some properties of  weakly mean equicontinuity. Let  $(X, d)$  be a compact metric space. For  $x \in X$  and  $\varepsilon>0$, denote  $B(x, \varepsilon)=\{y \in X: d(x, y)<\varepsilon\}$. Denote According to the product space  $X \times X=   \{(x, y): x, y \in X\}$  and the diagonal  $\Delta_X=\{(x, x): x \in X\} .$  A subset of  $X$  is called a  $G_{\delta}$  set if it can be expressed as a countable intersection of open sets;  a Borel set if it can be formed from open sets (or, equivalently, from closed sets) through the operations of countable union, countable intersection, and relative complement; a residual set if it contains the intersection of a countable collection of dense open sets. According to the Baire category theorem, a residual set is also dense in $ X$.

Let  $(X, T)$  be a dynamical system. The orbit of a point  $x \in X,\left\{x, T x, T^{2} x, \ldots,\right\} $, is denoted by   $\operatorname{Orb}(x, T)$. The $ \omega$-limit set of $ x$  is the set of limit points of the orbit sequence
$$
\omega(x, T)=\bigcap_{N \geq 0} \overline{\left\{T^{n} x: n \geq N\right\}}.
$$

If  $A$  is a nonempty closed subset of  $X$  and  $T A \subset A $, then  $\left(A,\left.T\right|_{A}\right)$  is called a subsystem of  $(X, T)$, where $ \left.T\right|_{A}$  is the restriction of  $T$  on  $A$. If there is no ambiguity, we will use the notation  $T$  instead of  $\left.T\right|_{A} .$

We say that a point  $x \in X$  is recurrent if  $x \in \omega(x, T)$. The system  $(X, T)$  is called  transitive if  $\overline{\operatorname{Orb}(x, T)}=X$  for some  $x \in X$, and such a point $ x $ is called a transitive point. Denote by $ \operatorname{Trans}(X, T)$  the set of transitive points of  $(X, T)$. With a Baire category argument, one can show that if  $(X, T)$  is transitive then  $\operatorname{Trans}(X, T) $ is a dense  $G_{\delta}$  subset of  $x$. For  $x \in X$  and  $U, V \subset X $, let  $N(x, U)=\left\{n \in \mathbb{Z}_{+}: T^{n} x \in U\right\}$  and  $N(U, V)=\left\{n \in \mathbb{Z}_{+}: U \cap T^{-n} V \neq \emptyset\right\} $. If $ U$  is a neighborhood of  $x$, then the set  $N(x, U)$  is called the set of return times of the point  $x$  to the neighborhood  $U $. Recall that a dynamical system  $(X, T)$  is called topologically transitive  if for every two nonempty open subsets  $U, V$  of  $X$  the set  $N(U, V)$  is infinite. The system $(X, T)$  is weakly mixing if $(X\times X, T\times T) $ is topologically transitive, $(X, T)$  is totally transitive if for any $n\in \mathbb{N}$, $(X, T^n)$ is topologically transitive.

The system  $(X, T)$  is said to be minimal if every point of  $X$  is a transitive point (i.e.,  $\operatorname{Trans}(X, T)=X$). A subset $ Y $ of  $X $ is called minimal if  $(Y, T)$  forms a minimal subsystem of  $(X, T)$. A point  $x \in X$  is called minimal if  the subsystem  ($\overline{\operatorname{Orb}(x, T)}, T|_{\overline{\operatorname{Orb}(x, T)}}$)  is minimal.

Let  $(X, T)$  and  $(Y, S)$   be  two topological dynamical systems and  $\pi: X \rightarrow Y $ is a continuous onto map which intertwines the actions (i.e.,  $\pi \circ T=S \circ \pi$), one says that  $(Y, S)$  is a factor of  $(X, T)$  or  $(X, T)$  is an extension of  $(Y, S)$, and  $\pi$  is a factor map. If  $\pi$  is a homeomorphism, then we say that  $\pi$  is a conjugacy and that the dynamical systems  $(X, T)$ and  $(Y, S)$  are conjugate.

A pair of points  $(x, y) \in X \times X$  is said to be proximal if for any  $\varepsilon>0$, there exists a positive integer  $n $ such that  $d\left(T^{n} x, T^{n} y\right)<\varepsilon$. Let  $P(X, T)$  denote the collection of all proximal pairs in  $(X, T)$. The dynamical system  $(X, T)$  is called proximal if any pair of two points in $X $ is proximal, i.e.,  $P(X, T)=X \times X$.
If  $(x, y) \in X \times X $ is not proximal, then it is said to be distal. A dynamical system  $(X, T)$  is called distal if any pair of distinct points in  $(X, T)$  is distal. A pair of points  $(x, y)$  is called regionally proximal if for every  $\varepsilon>0$, there exist two points  $x^{\prime}, y^{\prime} \in X$  with  $d\left(x, x^{\prime}\right)<\varepsilon$  and  $d\left(y, y^{\prime}\right)<\varepsilon$, and a positive integer  $n$  such that  $d\left(T^{n} x^{\prime}, T^{n} y^{\prime}\right)<\varepsilon$. Let $ RP(X, T)$  be the set of all regionally proximal pairs in  $(X, T)$.

A pair  $(x, y) \in X \times X$  is called  weakly mean asymptotic if
$$
\lim _{n \rightarrow \infty}\inf _{\sigma \in S_{n}} \frac{1}{n} \sum_{i=1}^{n} d\left(T^{i} x, T^{\sigma(i)}  y\right)=0.
$$

A pair  $(x, y) \in X \times X$  is called weakly mean proximal if
$$
\liminf _{n \rightarrow \infty} \inf _{\sigma \in S_{n}}\frac{1}{n} \sum_{i=1}^{n} d\left(T^{i} x, T^{\sigma(i)}  y\right)=0.
$$

A pair  $(x, y) \in X \times X$  is called strong mean proximal if
$$
\liminf _{n \rightarrow \infty} \sup_{\sigma \in S_{n}}\frac{1}{n} \sum_{i=1}^{n} d\left(T^{i} x, T^{\sigma(i)}  y\right)=0.
$$
The weakly mean asymptotic relation, weakly mean proximal  relation and strong mean proximal relation of  $(X, T)$, denoted by   $\operatorname{WMAsym}(T)$, $ \operatorname{WMProx}(T)$ and $ \operatorname{SMProx}(T)$, are the set of all weakly mean asymptotic pairs, weakly mean proximal pairs and strong mean proximal pairs respectively. 

For  $x \in X $, the weakly  mean asymptotic cell,  weakly mean proximal cell and strong mean proximal cell of  $x$  are defined by
$$
\operatorname{WMAsym}(T, x)=\{y \in X:(x, y) \in \operatorname{WMAsym}(T)\},
$$
$$
\operatorname{WMProx}(T, x)=\{y \in X:(x, y) \in \operatorname{WMProx}(T)\}
$$
and
$$
\operatorname{SMProx}(T, x)=\{y \in X:(x, y) \in \operatorname{SMProx}(T)\}
$$
respectively.

{\bf Proposition 3.1 }
Let  $(X, T)$  be a topological dynamical system and $\delta>0$. Then 

(1)  for every $x\in X$, weakly mean proximal cell of $x$ is a Borel subset of $X$ ;

(2) the weakly mean proximal relation of $(X, T)$  is a Borel subset of $X \times X $.

{\bf Proof } (1)   Since
$$
\operatorname{WMAsym}(T, x)=\bigcap_{m=1}^{\infty} \bigcup_{\ell=1}^{\infty}\left(\bigcap_{n \geq \ell}\bigcup_{\sigma \in S_{n}}\left\{y\in X: \frac{1}{n} \sum_{k=1}^{n} d\left(T^{k} x, T^{\sigma(k)} y\right)<\frac{1}{m}\right\}\right) .
$$
Then,  $\operatorname{WMAsym}(T, x)$  is a Borel subset of $X$.

(2) Since 
$$
\operatorname{WMAsym}(T)=\bigcap_{m=1}^{\infty} \bigcup_{\ell=1}^{\infty}\left(\bigcap_{n \geq \ell}\bigcup_{\sigma \in S_{n}}\left\{(x, y) \in X \times X: \frac{1}{n} \sum_{k=1}^{n} d\left(T^{k} x, T^{\sigma(k)} y\right)<\frac{1}{m}\right\}\right) .
$$
Then,  $\operatorname{WMAsym}(T)$  is a Borel subset of $X$.

{\bf Proposition 3.2}
Let  $(X, T)$  be a topological dynamical system and $\delta>0$. Then 

(1)  for every $x\in X$, strong mean proximal cell of $x$ is a Borel subset of $X$;

(2) the strong mean proximal relation of $(X, T)$  is a Borel subset of $X \times X $.

{\bf Proof } (1)   Since
$$
\operatorname{SMAsym}(T, x)=\bigcap_{m=1}^{\infty} \bigcap_{\ell=1}^{\infty}\left(\bigcup_{n \geq \ell}\bigcap_{\sigma \in S_{n}}\left\{y\in X: \frac{1}{n} \sum_{k=1}^{n} d\left(T^{k} x, T^{\sigma(k)} y\right)<\frac{1}{m}\right\}\right) .
$$
Then,  $\operatorname{SMAsym}(T, x)$  is a Borel subset of $X$.

(2) Since 
$$
\operatorname{SMAsym}(T)=\bigcap_{m=1}^{\infty} \bigcap_{\ell=1}^{\infty}\left(\bigcup_{n \geq \ell}\bigcap_{\sigma \in S_{n}}\left\{(x, y) \in X \times X: \frac{1}{n} \sum_{k=1}^{n} d\left(T^{k} x, T^{\sigma(k)} y\right)<\frac{1}{m}\right\}\right) .
$$
Then,  $\operatorname{SMAsym}(T)$  is a Borel subset of $X$.

{\bf Proposition 3.3} Let  $(X, T)$  be a topological dynamical system. If  $(X, T)$  is weakly mean equicontinuous, then every proximal pair is weakly mean asymptotic.

{\bf Proof} Let  $\left(x_{0}, y_{0}\right) \in X \times X$  be a proximal pair. As  $(X, T)$  is weakly mean equicontinuous, for every  $\varepsilon>0$  there exists  $\delta>0$  such that for any  $x, y \in X$  with  $d(x, y)<\delta ,$
$$
\limsup _{n \rightarrow \infty} \inf _{\sigma \in S_{n}}\frac{1}{n} \sum_{i=1}^{n} d\left(T^{i} x, T^{\sigma(i)}  y\right)<\varepsilon.
$$
Since $\left(x_{0}, y_{0}\right)$  is proximal, there exists  $k \in \mathbb{N}$  such that  $d\left(T^{k} x_{0}, T^{k} y_{0}\right)<\delta$. Then
$$
\begin{aligned}
\limsup _{n \rightarrow \infty}  \inf _{\sigma \in S_{n}}\frac{1}{n}\sum_{i=1}^{n} d\left(T^{i} x_{0}, T^{\sigma(i)} y_{0}\right)&=\limsup _{n \rightarrow \infty}\inf _{\sigma \in S_{n}} \frac{1}{n} \sum_{i=1}^{n} d\left(T^{i}\left(T^{k} x_{0}\right), T^{\sigma(i)}\left(T^{k} y_{0}\right)\right)\\
&<\varepsilon .
\end{aligned}
$$
As  $\varepsilon>0$  is arbitrary,  $\left(x_{0}, y_{0}\right)$  is weakly  mean asymptotic.

{\bf Proposition 3.4} Let  $(X, T)$  be a topological dynamical system. If  $(X, T)$  is weakly mean equicontinuous, then every regionally proximal pair is weakly mean asymptotic.

{\bf Proof} Let  $\left(x_{0}, y_{0}\right) \in X \times X$  be a regionally proximal pair. As  (X, T)  is weakly mean equicontinuous, for every  $\varepsilon>0$  there exists  $\delta>0$  such that whenever  $x, y \in X$  with  $d(x, y)<\delta ,$
$$
\limsup _{n \rightarrow \infty} \inf _{\sigma \in S_{n}}\frac{1}{n} \sum_{i=1}^{n} d\left(T^{i} x, T^{\sigma(i)}  y\right)<\frac{\varepsilon}{3}.
$$

For given $\delta$, there exist $u,v\in X$ and $k\in \mathbb{N}$ such that $d(x_0, u)<\delta$, $d(y_0, v)<\delta$ and $d(T^{k} u, T^{k} v)<\delta$. Then
$$ 
\begin{aligned}
	\limsup _{n \rightarrow \infty} \inf _{\sigma \in S_{n}}\frac{1}{n} \sum_{i=1}^{n} d\left(T^{i} x_{0}, T^{\sigma(i)} y_{0}\right)
	& \leq \limsup _{n \rightarrow \infty} \inf _{\sigma \in S_{n}}\frac{1}{n} \sum_{i=1}^{n} d(T^{i} x_{0}, T^{\sigma(i)} u)+\limsup _{n \rightarrow \infty} \inf _{\sigma \in S_{n}}\frac{1}{n} \sum_{i=1}^{n} d(T^{i} u, T^{\sigma(i)} v)) \\
	&+ \limsup _{n \rightarrow \infty} \inf _{\sigma \in S_{n}}\frac{1}{n} \sum_{i=1}^{n} d(T^{i} v, T^{\sigma(i)} y_0)) \\
	& \leq \limsup _{n \rightarrow \infty} \inf _{\sigma \in S_{n}}\frac{1}{n} \sum_{i=1}^{n} d(T^{i} x_{0}, T^{\sigma(i)} u)+\limsup _{n \rightarrow \infty} \inf _{\sigma \in S_{n}}\frac{1}{n} \sum_{i=1}^{n} d(T^{i} y_0, T^{\sigma(i)} v)) \\
	&+ \limsup _{n \rightarrow \infty} \inf _{\sigma \in S_{n}}\frac{1}{n} \sum_{i=1}^{n} d(T^{i} (T^{k}u), T^{\sigma(i)} (T^{k}v)))\\
	& <\frac{\varepsilon}{3}+\frac{\varepsilon}{3}+\frac{\varepsilon}{3}=\epsilon .
\end{aligned}
$$
As  $\varepsilon>0$  is arbitrary,  $\left(x_{0}, y_{0}\right)$  is weakly  mean asymptotic.

{\bf  Theorem 3.5} 
 Let  $\pi:(X, T) \rightarrow(Y, S)$  be a factor map between two topological dynamical systems. If  $(X, T)$  is weakly mean equicontimuous, then  $(Y, S) $ is weakly mean equicontimuous.
  
{\bf Proof} Assume the contrary that  $(Y, S)$  is not weakly mean equicontinuous. Then there is  $\delta>0$  and two sequences  $\left\{y_{k}\right\},\left\{z_{k}\right\}$  of points in  $Y $ with  $d_1\left(y_{k}, z_{k}\right)<\frac{1}{k} $ such that for each  $k \geq 1$ ,
$$
\limsup _{n \rightarrow \infty} \inf _{\sigma \in S_{n}}\frac{1}{n} \sum_{i=1}^{n} d_1\left(S^{i}\left(y_{k}\right), S^{\sigma(i)}\left(z_{k}\right)\right) > 2 \delta.
$$
According to the compactness of  $Y$, we may assume that  $\lim _{k \rightarrow \infty} y_{k}=\lim _{k \rightarrow \infty} z_{k}=y \in Y$ . For each  $k \geq 1 ,$
$$
\begin{aligned}
\underset{n \rightarrow \infty}{\limsup } \inf _{\sigma \in S_{n}}\frac{1}{n} \sum_{i=0}^{n-1}  d_1\left(S^{i}\left(y_{k}\right), S^{\sigma(i)}(y)\right) &+\underset{n \rightarrow \infty}{\limsup } \inf _{\sigma \in S_{n}}\frac{1}{n} \sum_{i=0}^{n-1}  d_1\left(S^{i}\left(z_{k}\right), S^{\sigma(i)}(y)\right)\\
 &\geq \limsup _{n \rightarrow \infty} \inf _{\sigma \in S_{n}}\frac{1}{n} \sum_{i=0}^{n-1}  d_1\left(S^{i}\left(y_{k}\right), S^{\sigma(i)}\left(z_{k}\right)\right) .
\end{aligned}
$$

Then either  
$$\limsup _{n \rightarrow \infty} \inf _{\sigma \in S_{n}}\frac{1}{n} \sum_{i=1}^{n}  d_1\left(S^{i}\left(y_{k}\right), S^{\sigma(i)}(y)\right) > \delta$$  
or  
$$\limsup _{n \rightarrow \infty} \inf _{\sigma \in S_{n}}\frac{1}{n} \sum_{i=1}^{n}  d_1\left(S^{i}\left(z_{k}\right), S^{\sigma(i)}(y)\right) > \delta . $$

Without loss of generality, assume that
$$
\limsup _{n \rightarrow \infty} \inf _{\sigma \in S_{n}}\frac{1}{n} \sum_{i=1}^{n}  d_1\left(S^{i}\left(y_{k}\right), S^{\sigma(i)}(y)\right) > \delta
$$
holds for all $ k \geq 1$. Choose a sequence  $\left\{x_{k}\right\}$  in  $X$  with  $\pi\left(x_{k}\right)=y_{k} $. According to the compactness of  $X$, we can assume that  $\lim _{k \rightarrow \infty} x_{k}=x$. Then $ \pi(x)=y .$

Choose an  $\eta>0$  such that  $\eta<\delta /(1+\operatorname{diam}(Y))$. According to the continuity of  $\pi$, there exists  $\theta \in(0, \eta) $ such that if $ d(u, v)<\theta$  then  $ d_1(\pi(u), \pi(v))<\eta $. Since  $(X, T) $  is weakly mean equicontinuous, there is  $\varepsilon>0 $ such that
$$
\underset{n \rightarrow \infty}{\limsup }  \inf _{\sigma \in S_{n}}\frac{1}{n} \sum_{i=1}^{n} d\left(T^{i} w, T^{\sigma(i)} x\right)<\theta^{2}
$$
for all  $w \in B(x, \varepsilon)$. Choose  $x_{j} \in B(x, \varepsilon) $. Let  $E_{\sigma}=\left\{1\leq i\leq n: d\left(T^{i} x_{j}, T^{\sigma(i)} x\right) \geq \theta\right\}$.  One has
$$
\theta^{2}>\limsup_{n \rightarrow \infty} \inf _{\sigma \in S_{n}}\frac{1}{n} \sum_{i=1}^{n} d\left(T^{i} x_{j}, T^{\sigma(i)} x\right) \geq \limsup _{n \rightarrow \infty} \frac{1}{n}(\theta\inf _{\sigma \in S_{n}} \#( E_{\sigma}))=\theta \limsup _{n \rightarrow \infty} \frac{1}{n}(\inf _{\sigma \in S_{n}} \#( E_{\sigma})),
$$
 and then $ \limsup _{n \rightarrow \infty} \frac{1}{n}(\inf _{\sigma \in S_{n}} \#( E_{\sigma}))<\theta$. Let  $F_{\sigma}=\left\{1\leq i\leq n:  d_1\left(S^{i} y_{j}, S^{\sigma(i)} y\right) \geq \eta\right\} $. According to the choice of  $\theta,$  we have  $F_{\sigma} \subset E_{\sigma} $. Then  $\limsup _{n \rightarrow \infty} \frac{1}{n}(\inf _{\sigma \in S_{n}} \#( F_{\sigma}))<\theta$  and
$$
\begin{aligned}
	\limsup _{n \rightarrow \infty} \inf _{\sigma \in S_{n}}\frac{1}{n} \sum_{i=1}^{n}  d_1\left(S^{i}\left(y_{j}\right), S^{\sigma(i)}(y)\right) & \leq \limsup _{n \rightarrow \infty} \frac{1}{n}(\operatorname{diam}(Y) \inf _{\sigma \in S_{n}}\#( F_{\sigma})+\eta n) \\
	& \leq \operatorname{diam}(Y) \limsup _{n \rightarrow \infty} \frac{1}{n}(\inf _{\sigma \in S_{n}} \#( F_{\sigma}))+\eta\\
	 &\leq \operatorname{diam}(Y) \theta+\eta \\
	& \leq(\operatorname{diam}(Y)+1) \eta<\delta,
\end{aligned}
$$
which is a contradiction. Thus  $(Y, S)$  is  weakly mean equicontinuous.

{\bf Theorem 3.6} Let  $(X, T)$  and  $(Y, S)$  be two topological dynamical systems. If $(X \times Y, T \times S)$  is weakly mean equicontinuous, then  $(X, T)$  and  $(Y, S)$  are weakly mean equicontinuous.

{\bf Proof} Assume that  $(X \times Y, T \times S)$  is weakly mean equicontinuous, then for any  $\varepsilon>0$  and any $x\in X, y\in Y$, as $(X \times Y, T \times S)$  is weakly mean equicontinuous, then there exists a $\delta>0$, such that $$
\limsup _{n \rightarrow \infty} \inf _{\sigma \in S_{n}}\frac{1}{n} \sum_{i=1}^{n} \rho\left((T\times S)^{i} (x,y), (T\times S)^{\sigma(i)} (u,v)\right)<\varepsilon .
$$
for any $(u,v)\in X\times Y$ with $\rho ((x,y),(u,v))<\delta.$

Since
$$
\begin{aligned}
	\rho((T\times S)^{i} (x,y), (T\times S)^{\sigma(i)} (u,v) 
	=d (T^{i}(x), T^{\sigma(i)}(u))+d_{1}(S^{i}(y), S^{\sigma(i)}(v)) 	
\end{aligned}
$$
and
$$\rho ((x,y),(u,v))=d (x, u)+d_{1}(y, v)<\delta. $$
Then,
$$
\limsup _{n \rightarrow \infty} \inf _{\sigma \in S_{n}}\frac{1}{n} \sum_{i=1}^{n}d \left(T^{i}\left(x)\right), T^{\sigma(i)} \left(u\right)\right)<\varepsilon \quad \text { and } \quad \limsup _{n \rightarrow \infty} \inf _{\sigma \in S_{n}}\frac{1}{n} \sum_{i=1}^{n}d_{1}\left(S^{i}\left(y)\right), S^{\sigma(i)} \left(v\right)\right)<\varepsilon
$$
for any $d(u,x)<\delta/2$ and $d(v,y)<\delta/2$. This implies $(X, T)$  and  $(Y, S)$  are weakly mean equicontinuous. The proof is completed.

{\bf Remark 3.7}
In [1], the authors showed that  $(X, T)$ is mean equicontinuous if and only if  $(X \times X, T \times T)$ is  mean equicontinuous; in [4], the authors showed that $(X, T)$ is mean equicontinuous if and only if $(X \times X, T \times T)$ is  weakly mean equicontinuous; as we know weakly mean equicontinuity is weakly than mean equicontinuity, then we get that the fact $(X, T)$ is weakly mean equicontinuous is not equivalent to the fact that $(X \times X, T \times T)$ is weakly mean equicontinuous.

In [2], Fomin introduced a notion called mean-$L$-stable, and showed that if a minimal system is mean-$L$-stable then it is uniquely ergodic. Then Li et al. in [1] introduced the notion of mean equicontinuity, and proved that a dynamical system is mean equicontinuous if and only if it is mean-$L$-stable. Inspired by the aboved ideas, we have the similar results in weakly mean equicontinuous dynamical system as follows.

{\bf Theorem 3.8}
A dynamical system $(X, T)$  is weakly mean equicontinuous if and only if if for every $\varepsilon>0 $, there is a $ \delta>0 $ such that  $d(x, y)<\delta$  implies  $$\limsup_{n \rightarrow \infty} \frac{1}{n} \inf _{\sigma \in S_{n}} \#\{1 \leq i \leq n: d\left(T^{i} x, T^{\sigma(i)}  y\right)>\varepsilon\}<\varepsilon.$$

{\bf Proof}  Assume that  $(X, T) $ is weakly mean equicontinuous. For every  $\varepsilon>0$  there exists a  $\delta>0$  such that 
$$
\limsup _{n \rightarrow \infty}\inf _{\sigma \in S_{n}}  \frac{1}{n} \sum_{i=1}^{n} d\left(T^{i} x, T^{\sigma(i)}  y\right)<\varepsilon^{2}
$$
for all  $x, y \in X $ with  $d(x, y)<\delta $. Let  $E_{n}^{\sigma}=\left\{1\leq i\leq n: d\left(T^{i} x, T^{\sigma(i)}  y\right)> \varepsilon\right\}$. One has
$$
\begin{aligned}
\varepsilon^{2}&>\limsup _{n \rightarrow \infty}\inf _{\sigma \in S_{n}}  \frac{1}{n} \sum_{i=1}^{n} d\left(T^{i} x, T^{\sigma(i)}  y\right) \\&\geq \limsup _{n \rightarrow \infty} \inf _{\sigma \in S_{n}} \frac{1}{n}(\varepsilon \#(E_{n}^{\sigma}))\\&=\varepsilon \limsup_{n \rightarrow \infty} \frac{1}{n} \inf _{\sigma \in S_{n}} \#\{1 \leq i \leq n: d\left(T^{i} x, T^{\sigma(i)}  y\right)>\varepsilon\}
\end{aligned}
$$
which implies  $\limsup_{n \rightarrow \infty} \frac{1}{n} \inf _{\sigma \in S_{n}} \#\{1 \leq i \leq n: d\left(T^{i} x, T^{\sigma(i)}  y\right)>\varepsilon\}<\varepsilon$.

Conversely, fix a positive number  $\varepsilon>0 $ and choose a positive number  $\eta<\frac{\varepsilon}{\operatorname{diam}(X)+1} $. There is a $ \delta>0 $ such that  $d(x, y)<\delta$  implies  $$\limsup_{n \rightarrow \infty} \frac{1}{n} \inf _{\sigma \in S_{n}} \#\{1 \leq i \leq n: d\left(T^{i} x, T^{\sigma(i)}  y\right)>\eta\}<\eta.$$
Let  $F_{n}^{\sigma}=\{1\leq i\leq n: d\left(T^{i} x, T^{\sigma(i)}  y\right)>\eta\} $. Then
$$
\begin{aligned}
	\limsup _{n \rightarrow \infty}\inf _{\sigma \in S_{n}}  \frac{1}{n} \sum_{i=1}^{n} d\left(T^{i} x, T^{\sigma(i)}  y\right) & \leq \limsup _{n \rightarrow \infty} \frac{1}{n}(\operatorname{diam}(X) \inf _{\sigma \in S_{n}}\#(F_{n}^{\sigma})+\eta n) \\
	& \leq \operatorname{diam}(X) \eta+\eta\\
	&<\varepsilon,
\end{aligned}
$$
which implies  $(X, T)$  is weakly mean equicontinuous.

{\bf Lemma  3.9}(see [4])  Let  $(X, T)$  be a topological dynamical system. If  $(X, T)$  is  weakly mean equicontinuous, then for every  $x \in X$, ($\overline{\operatorname{Orb}(x, T)}, T|_{\overline{\operatorname{Orb}(x, T)}}$)  is uniquely ergodic. In particular, if  $(X, T)$  is also transitive, then  $(X, T)$  is uniquely ergodic.

{\bf Lemma  3.10} Let  $(X, T)$  be a weakly mean equicontinuous system and  $\mu$  be an ergodic measure on  $X$, then every point of  $\operatorname{supp}(\mu)$  is minimal.

{\bf Proof}  Since  $\mu$  is an ergodic measure on  $X$, then $(\operatorname{supp}(\mu), T)$ is a transitive system. According to Lemma  3.9,  $(\operatorname{supp}(\mu), T)$  is uniquely ergodic, hence, it is minimal.

{\bf Lemma  3.11}  Let  $(X, T)$  be a transitive and weakly mean equicontinuous topological dynamical system. Then  $(X, T)$  is minimal if and only if  $\operatorname{supp}(X, T)=X $.

{\bf Proof}  Necessity is easily accessible. Now we only need to indicate the sufficiency. Choose any transitive point  $x \in X $, then  ($\overline{\operatorname{Orb}(x, T)}, T|_{\overline{\operatorname{Orb}(x, T)}}$)  is uniquely ergodic (by Lemma 3.9) with measure  $\mu $. Then  $(\operatorname{supp}(\mu), T|_{\operatorname{supp}(\mu)})$  is a minimal system. Since  $\operatorname{supp}(X, T)=X$, then  $\operatorname{supp}(\mu)=X $. Thus  $(X, T)$  is a minimal system.

\section{Almost weakly mean equicontinuity and almost weakly  equicontimuous in the mean }
Inspired by the existing ideas, in this section,  we  will introduce the  notions of almost weakly mean equicontinuity, almost weakly  equicontimuous in the mean and  study some properties of almost weakly mean equicontinuity and almost weakly  equicontimuous in the mean. 

Let  $(X, T)$  be a topological dynamical system. A point  $x \in X$  is called weakly mean equicontimuous if for every  $\varepsilon>0 $, there is  $\delta>0$  such that for every  $y \in B(x, \delta) $,
$$
\limsup _{n \rightarrow \infty} \inf _{\sigma \in S_{n}}\frac{1}{n} \sum_{i=1}^{n} d\left(T^{i} x, T^{\sigma(i)}  y\right)<\varepsilon.$$

A point  $x \in X$  is called weakly  equicontimuous in the mean if for every  $\varepsilon>0 $, there is  $\delta>0$  such that for every  $y \in B(x, \delta) $,
$$
\sup _{n \in \mathbb{N}} \inf _{\sigma \in S_{n}}\frac{1}{n} \sum_{i=1}^{n} d\left(T^{i} x, T^{\sigma(i)}  y\right)<\varepsilon.$$

We say that a topological dynamical system $(X, T)$  is weakly mean equicontinuous (resp. weakly  equicontimuous in the mean) if and only if every point in  $X$  is weakly mean equicontinuous (resp. weakly  equicontimuous in the mean). 

A  topological dynamical   system $ (X, T) $ is called almost weakly mean equicontinuous (resp. almost weakly  equicontimuous in the mean)  if it is transitive and has at least one weakly mean equicontinuous point  (resp. weakly  equicontimuous in the mean point). We show that for a transitive system  $(X, T)$, the set of weakly mean equicontinuous points or weakly  equicontimuous in the mean points is either empty or residual. If in addition  $(X, T)$  is almost  weakly mean equicontinuous (resp. weakly  equicontimuous in the mean), then every transitive point is weakly mean equicontinuous (resp. weakly  equicontimuous in the mean). 

Let  $\mathscr{E}$  (resp. $\mathscr{E}^{*}$) denote the set of all weakly mean equicontinuous points (resp. all weakly equicontinuous in the mean points). For every  $\varepsilon>0 $, let
$$
\mathscr{E}_{\varepsilon}=\left\{x \in X: \exists \delta>0, \forall y, z \in B(x, \delta),\limsup _{n \rightarrow \infty} \inf _{\sigma \in S_{n}}\frac{1}{n} \sum_{i=1}^{n} d\left(T^{i} y, T^{\sigma(i)}  z\right)<\varepsilon\right\}.
$$
$$
\mathscr{E}^{*}_{\varepsilon}=\left\{x \in X: \exists \delta>0, \forall y, z \in B(x, \delta),\sup _{n \in \mathbb{N}} \inf _{\sigma \in S_{n}}\frac{1}{n} \sum_{i=1}^{n} d\left(T^{i} y, T^{\sigma(i)}  z\right)<\varepsilon\right\}.
$$
{\bf Proposition 4.1}
Let  $(X, T)$  be a topological dynamical system and  $\varepsilon>0$. Then  $\mathscr{E}_{\varepsilon}$  is open and inversely invariant, that is  $T^{-1} \mathscr{E}_{\varepsilon} \subset \mathscr{E}_{\varepsilon}$. Moreover,  $\mathscr{E}=\bigcap_{m=1}^{\infty} \mathscr{E}_{\frac{1}{m}}$  is a  $G_{\delta}$  subset of $ X$.

{\bf Proof} Let  $x \in \mathscr{E}_{\varepsilon}$. Choose $ \delta>0$  satisfying the condition from the definition of  $\mathscr{E}_{\varepsilon}$  for $ x$. Fix  $y \in B\left(x, \frac{\delta}{2}\right)$. If  $u, v\in B\left(y, \frac{\delta}{2}\right)$, then  $u, v \in B(x, \delta)$, and
$$
\limsup _{n \rightarrow \infty}  \inf _{\sigma \in S_{n}}\frac{1}{n} \sum_{i=1}^{n} d\left(T^{i} u, T^{\sigma(i)} v\right)<\varepsilon
$$
This shows that  $B\left(x, \frac{\delta}{2}\right) \subset \mathscr{E}_{\varepsilon}$  and thus  $\mathscr{E}_{\varepsilon}$  is open.
Let  $x \in X$  with  $T x \in \mathscr{E}_{\varepsilon}$. Choose  $\delta>0$  satisfying the condition from the definition of  $\mathscr{E}_{\varepsilon}$  for  $T x$. According to the continuity of  $T$, there exists  $\eta>0$  such that $ d(T u, T v)<\delta$  for any  $u, v \in B(x, \eta)$. If  $u, v \in B(x, \eta)$ , then  $T u, T v \in B(T x, \delta)$. Thus
$$
\underset{n \rightarrow \infty}{\limsup } \inf _{\sigma \in S_{n}}\frac{1}{n} \sum_{i=1}^{n} d\left(T^{i} u, T^{\sigma(i)} v\right)=\underset{n \rightarrow \infty}{\limsup } \inf _{\sigma \in S_{n}}\frac{1}{n} \sum_{i=1}^{n} d\left(T^{i}(T u), T^{\sigma(i)}(T v)\right)<\varepsilon,
$$
This implies  $x \in \mathscr{E}_{\varepsilon} .$

If $ x \in X$  belongs to all  $\mathscr{E}_{\frac{1}{m}} $, then clearly  $x \in \mathscr{E}$. Conversely, if $ x \in \mathscr{E}$  and  $m>0$, then there exists  $\delta>0$  such that for all  $y \in B(x, \delta) $
$$
\limsup _{n \rightarrow \infty} \inf _{\sigma \in S_{n}}\frac{1}{n} \sum_{i=1}^{n} d\left(T^{i} x, T^{\sigma(i)} y\right)<\frac{1}{2 m}.$$

If  $y, z \in B(x, \delta)$, one has
$$
\begin{aligned}
	\limsup _{n \rightarrow \infty} \inf _{\sigma \in S_{n}}\frac{1}{n} \sum_{i=1}^{n} d\left(T^{i} y, T^{\sigma(i)} z\right)
	& \leq \limsup _{n \rightarrow \infty} \inf _{\sigma \in S_{n}}\frac{1}{n} \sum_{i=1}^{n} d(T^{\sigma(i)} x, T^{i} y)+\limsup _{n \rightarrow \infty} \inf _{\sigma \in S_{n}}\frac{1}{n} \sum_{i=1}^{n} d(T^{i} x, T^{\sigma(i)} z)) \\
	& \leq \limsup _{n \rightarrow \infty} \inf _{\sigma \in S_{n}}\frac{1}{n} \sum_{i=1}^{n} d(T^{i} x, T^{\sigma(i)} y)+\limsup _{n \rightarrow \infty} \inf _{\sigma \in S_{n}}\frac{1}{n} \sum_{i=1}^{n} d(T^{i} x, T^{\sigma(i)} z)) \\
	& <\frac{1}{2 m}+\frac{1}{2 m}=\frac{1}{m} .
\end{aligned}
$$
Hence,  $x \in \mathscr{E}_{\frac{1}{m}}$. This finishes the proof.

{\bf Proposition 4.2 }
 Let  $(X, T) $ be a transitive topological dynamical system.
 
(1) The set of weakly mean equicontinuous points is either empty or residual. If in addition  $(X, T)$  is almost weakly mean equicontimuous, then every transitive point is weakly mean equicontimuous.

(2) If  $(X, T)$  is minimal and almost weakly  mean equicontimuous, then it is weakly mean equicontinuous.

{\bf Proof} According to the transitivity of  $(X, T)$, every  $\mathcal{E}_{\varepsilon}$  is either empty or dense, since  $\mathcal{E}_{\varepsilon}$  is open and inversely invariant. Then  $\mathcal{E}$  is either empty or residual according to the Baire Category theorem.

If  $\mathcal{E}$  is residual, then every  $\mathcal{E}_{\varepsilon}$  is open and dense. Let  $x \in X$  be a transitive point and  $\varepsilon>0$. Then there exists some $ k \in \mathbb{Z}_{+}$ such that  $T^{k} x \in \mathcal{E}_{\varepsilon}$, and as  $\mathcal{E}_{\varepsilon}$  is inversely invariant,  $x \in \mathcal{E}_{\varepsilon}$. Thus  $x \in \mathcal{E} .$

Furthermore, if we assume $(X, T)$  is minimal then we deduce that $\mathscr{E}_{\varepsilon}=X$ for every $\varepsilon>0$, and therefore $(X, T)$  is weakly mean equicontinuous.

{\bf Proposition  4.3} Let  $(X, T)$  be a topological dynamical system and  $x \in X $. If  $y \in \mathscr{E} \cap \overline{\operatorname{Orb}(x, T)}$  then  $x \in \mathscr{E}$  and  $\overline{F}(x, y) =0 $.

{\bf Proof} Let  $\varepsilon>0$  and  $y \in \mathscr{E} \subset \mathscr{E}_{\varepsilon} $. There is  $\delta>0$  such that $ B(y,\delta) \subset \mathscr{E}_{\varepsilon}$  and  $\overline{F}(y, z) < \varepsilon$  for every  $z \in B(y,\delta)$. Since  $y\in \overline{\operatorname{Orb}(x, T)} $, there is  $n \geq 0$  such that  $T^{n} x \in B(y,\delta)$. This implies that  $x \in \mathscr{E}_{\varepsilon}$  by the inversely invariance of  $\mathscr{E}_{\varepsilon}$ and  $\overline{F}(x, y)< \varepsilon $. Since  $\varepsilon>0$  was arbitrary, we conclude that  $x \in \mathscr{E}$  and  $\overline{F}(x, y)=0 $.

By similar proofs of Proposition 4.1 and Proposition 4.2, we can deduce the following results.

{\bf Proposition 4.4}
Let  $(X, T)$  be a topological dynamical system and  $\varepsilon>0$. Then we have

(1) $\mathscr{E}^{*}_{\varepsilon}$  is open and inversely invariant, that is  $T^{-1} \mathscr{E}^{*}_{\varepsilon} \subset \mathscr{E}^{*}_{\varepsilon}$. Moreover,  $\mathscr{E}^{*}=\bigcap_{m=1}^{\infty} \mathscr{E}^{*}_{\frac{1}{m}}$  is a  $G_{\delta}$  subset of $ X$.

(2) The set of weakly equicontinuous in the mean points is either empty or residual. If in addition  $(X, T)$  is almost weakly equicontimuous in the mean, then every transitive point is weakly equicontimuous in the mean.

(3) If  $(X, T)$  is minimal and almost weakly equicontimuous in the mean, then it is weakly equicontinuous in the mean. 

Note that a topological dynamical system is mean equicontinuous if and only if it is equicontinuous in the mean[5,13]. The following results show the relationship between weakly equicontinuous in the mean and weakly mean equicontinuous.

{\bf Theorem  4.5} Let  $(X, T)$  be a minimal topological dynamical system. Then  $(X, T)$  is  weakly mean equicontinuous if and only if weakly equicontinuous in the mean.

{\bf Proof}  It is easy to show that weakly equicontinuity in the mean implies weakly mean equicontinuity. Assume that  $(X, T)$  is weakly mean equicontinuous. For any  $\varepsilon>0$  there is  $\delta_{1}>0$  such that if  $d(x, y)<\delta_{1}$  then
$$
\limsup _{n \rightarrow \infty} \inf _{\sigma \in S_{n}}\frac{1}{n} \sum_{i=1}^{n} d\left(T^{i} x, T^{\sigma(i)}y\right)<\frac{\varepsilon}{8}.
$$
Fix  $z \in X$. For each  $m\in \mathbb{N}$ , let
$$
A_{m}=\left\{x \in \overline{B\left(z, \delta_{1} / 2\right)}: \inf _{\sigma \in S_{n}}\frac{1}{n} \sum_{i=1}^{n} d\left(T^{i} x, T^{\sigma(i)} z\right) \leq \frac{\varepsilon}{4}, n=m, m+1, \ldots\right\} .
$$
Then $A_{m}$  is closed and  $\overline{B\left(z, \delta_{1} / 2\right)}=\bigcup_{m=1}^{\infty} A_{m}$. By the Baire Category Theorem, there is  $m_{1} \in \mathbb{N}$  such that  $A_{m_{1}} $contains an open subset $ U $ of $ X$. By the minimality,  we have that there is  $m_{2} \in \mathbb{N}$  with  $\bigcup_{i=0}^{m_{2}-1} T^{-i} U=X $. Let  $\delta_{2}$  be the Lebesgue number of the open cover  $\left\{T^{-i} U: 0 \leq i \leq m_{2}-1\right\}$  of  $X $. Let  $m=\max \left\{m_{1}, 2m_{2}\right\} $. With the continuity of  $T$, there exists  $\delta_{3}>0$  such that if  $d(x, y)<\delta_{3}$  implies that  $d\left(T^{i} x, T^{i} y\right)<\frac{\varepsilon}{4}$  for any  $0 \leq i \leq m $. Put  $\delta=\min \left\{\delta_{2}, \delta_{3}\right\} $. Let  $x, y \in X$  with  $d(x, y)<\delta$  and  $n \in \mathbb{N} $. 

If  $n \leq m$, then
$$
\inf _{\sigma \in S_{n}}\frac{1}{n} \sum_{i=1}^{n} d\left(T^{i} x, T^{\sigma(i)} y\right) \leq \frac{1}{n} \cdot n \cdot \frac{\varepsilon}{4}<\varepsilon.
$$

If  $n>m$, there exists  $0 \leq k \leq m_{2}-1 $ such that  $x, y \in T^{-k} U$, then  $T^{k} x, T^{k} y \in U$, one has
$$
\begin{aligned}
	\inf _{\sigma \in S_{n}}\frac{1}{n} \sum_{i=1}^{n} d\left(T^{i} x, T^{\sigma(i)} y\right) &\leq \inf _{\sigma \in S_{k}}\frac{1}{n} \sum_{i=1}^{k} d\left(T^{i} x, T^{\sigma(i)} y\right)+\inf _{\sigma \in S_{n}}\frac{1}{n} \sum_{i=1}^{n} d\left(T^{i} T^{k} x, T^{\sigma(i)} T^{k} y\right) \\
	& \leq \frac{\varepsilon}{4}+\inf _{\sigma \in S_{n}}\frac{1}{n} \sum_{i=1}^{n} d\left(T^{i} T^{k} x, T^{\sigma(i)} z\right)+\inf _{\sigma \in S_{n}}\frac{1}{n} \sum_{i=1}^{n} d\left(T^{\sigma(i)} T^{k} y, T^{i} z\right) \\
    & = \frac{\varepsilon}{4}+\inf _{\sigma \in S_{n}}\frac{1}{n} \sum_{i=1}^{n} d\left(T^{i} T^{k} x, T^{\sigma(i)} z\right)+\inf _{\sigma \in S_{n}}\frac{1}{n} \sum_{i=1}^{n} d\left(T^{i} T^{k} y, T^{\sigma(i)} z\right) \\
	& \leq \frac{\varepsilon}{4}+\frac{\varepsilon}{4}+\frac{\varepsilon}{4}<\varepsilon .
\end{aligned}
$$
Therefore,  $\sup _{n \in \mathbb{N}} \inf _{\sigma \in S_{n}}\frac{1}{n} \sum_{i=1}^{n} d\left(T^{i} x, T^{\sigma(i)}  y\right)<\varepsilon$. This implies that  $(X, T)$  is weakly equicontinuous in the mean.

{\bf Remark  4.6} Notice that in [4], the authors showed that  for a general topological dynamical system  that dose not require minimality,  weakly mean equicontinuity is equivalent to  weakly equicontinuity in the mean.

Let  $\pi: X \rightarrow $ Y  be a map. The map  $\pi$  is called open if for every nonempty open subset $ U$  of  $X$, $\pi(U)$  is open in  $Y$, and semi-open if  $\pi(U)$  has nonempty interior in  $Y$. We say that  $\pi$  is open at a point  $x \in X$  for for every neighborhood  $U $ of $ x$, $ \pi(U)$  is a neighborhood of  $\pi(x)$. 

{\bf Theorem  4.7}
Let  $\pi:(X, T) \rightarrow(Y, S)$  be a factor map between two topological dynamical systems. Let  $x \in X$  be a weakly mean equicontinuous point and suppose that  $\pi$  is open at  $x.$ Then  $y=\pi(x)$  is also a  weakly mean equicontinuous point. 

{\bf Proof} If  $y$  is not a weakly mean equicontinuous point, then there exists  $\delta>0$  and a sequence  $y_{j} \rightarrow y$  such that for any  $j$, one has 
$$
\limsup _{n \rightarrow \infty} \inf _{\sigma \in S_{n}}\frac{1}{n} \sum_{i=1}^{n} d_1\left(S^{i} y_{j}, S^{\sigma(i)} y\right)>\delta .
$$

According to $\pi$  is open at  $x$, there exists a sequence  $x_{k} \rightarrow x$  with $ \pi\left(x_{k}\right)=y_{k} $. 
Choose an  $\eta>0$  such that  $\eta<\delta /(1+\operatorname{diam}(Y))$. According to the continuity of  $\pi$, there exists  $\theta \in(0, \eta) $ such that if $ d(u, v)<\theta$  then  $ d_1(\pi(u), \pi(v))<\eta $. 

Since $ x$  is a weakly mean equicontinuous point, there exists $\epsilon>0$, for large enough  $k$ with $x_k \in B(x, \varepsilon)$, we have
$$
\underset{n \rightarrow \infty}{\limsup }  \inf _{\sigma \in S_{n}}\frac{1}{n} \sum_{i=1}^{n} d\left(T^{i} x_k, T^{\sigma(i)} x\right)<\theta^{2}.
$$
 Choose  $x_{k} \in B(x, \varepsilon) $. Let  $E_{\sigma}=\left\{1\leq i\leq n: d\left(T^{i} x_{k}, T^{\sigma(i)} x\right) \geq \theta\right\}$.  we get that
$$
\theta^{2}>\limsup_{n \rightarrow \infty} \inf _{\sigma \in S_{n}}\frac{1}{n} \sum_{i=1}^{n} d\left(T^{i} x_{k}, T^{\sigma(i)} x\right) \geq \limsup _{n \rightarrow \infty} \frac{1}{n}(\theta\inf _{\sigma \in S_{n}} \#( E_{\sigma}))=\theta \limsup _{n \rightarrow \infty} \frac{1}{n}(\inf _{\sigma \in S_{n}} \#( E_{\sigma})),
$$
and then $ \limsup _{n \rightarrow \infty} \frac{1}{n}(\inf _{\sigma \in S_{n}} \#( E_{\sigma}))<\theta$. 

Set  $F_{\sigma}=\left\{1\leq i\leq n:  d_1\left(S^{i} y_{k}, S^{\sigma(i)} y\right) \geq \eta\right\} $. According to the choice of  $\theta$,  we deduce that  $F_{\sigma} \subset E_{\sigma} $. Then  $\limsup _{n \rightarrow \infty} \frac{1}{n}(\inf _{\sigma \in S_{n}} \#( F_{\sigma}))<\theta$  and
$$
\begin{aligned}
	\limsup _{n \rightarrow \infty} \inf _{\sigma \in S_{n}}\frac{1}{n} \sum_{i=1}^{n}  d_1\left(S^{i}\left(y_{k}\right), S^{\sigma(i)}(y)\right) & \leq \limsup _{n \rightarrow \infty} \frac{1}{n}(\operatorname{diam}(Y) \inf _{\sigma \in S_{n}}\#( F_{\sigma})+\eta n) \\
	& \leq \operatorname{diam}(Y) \limsup _{n \rightarrow \infty} \frac{1}{n}(\inf _{\sigma \in S_{n}} \#( F_{\sigma}))+\eta\\
	&\leq \operatorname{diam}(Y) \theta+\eta \\
	& \leq(\operatorname{diam}(Y)+1) \eta<\delta.
\end{aligned}
$$
 This is a contradiction. Thus,  $y$  is a weakly mean equicontinuous point. The proof is finished.

{\bf Theorem 4.8}
 Let  $\pi:(X, T) \rightarrow(Y, S)$  be a factor map between transitive systems. Suppose that  $(X, T)$  is almost weakly mean equicontimuous and  $\pi$  is semi-open. Then  $(Y, S)$  is also almost weakly mean equicontimuous.

{\bf Proof}  Since  $\pi$  is semi-open, then the set  $\{x \in X: \pi \  \text{is open at} \  x\}$  is residual in  $X.$  Pick a transitive point $ x \in X$  such that  $\pi$  is open at  $x$, then  $\pi(x)$  is also a transitive point in  $Y$. Since  $(X, T) $ is almost weakly mean equicontinuous,  and $x$  is a weakly mean equicontinuous point. According to Theorem 4.8, $\pi(x)$  is a weakly mean equicontinuous point, which implies that  $(Y, S)$  is almost weakly mean equicontinuous.

\section{ Strong Mean sensitivity and Strong sensitivity in the mean}
In this section, we will study the opposite side of weakly mean equicontinuity(resp. weakly  equicontinuity in the mean) as strong mean sensitivity (resp. strong  sensitivity in the mean). It follows that if a topological dynamical system  $(X, T)$  is minimal then  $(X, T)$  is either weakly mean equicontinuous (resp. weakly equicontinuity in the mean) or strong mean sensitive (resp. strong  sensitivity in the mean), and furthermore if  $(X, T) $ is a transitive system then  $(X, T)$  is either almost weakly mean equicontinuous (resp. almost weakly  equicontinuity in the mean) or strong mean sensitive (resp. strong  sensitivity in the mean). Moreover, we demonstrate that  a dynamical system is strong  sensitivity in the mean if and only if it is strong mean sensitivity, and it follows that for a transitive dynamical system,  almost weakly mean equicontinuity is equivalent to  weakly equicontinuity in the mean.

A topological dynamical system  $(X, T)$  is strong mean sensitive if there exists  $\delta>0$  such that for every  $x \in X$  and every  $\varepsilon>0$  there is  $y \in B(x, \varepsilon)$  satisfying
	$$
	\limsup _{n \rightarrow+\infty} \inf _{\sigma \in S_{n}} \frac{1}{n} \sum_{k=1}^{n} d\left(T^{k} x, T^{\sigma(k)} y\right)>\delta.$$
	
A topological dynamical system  $(X, T)$  is strong sensitive in the mean if there exists  $\delta>0$  such that for every  $x \in X$  and every  $\varepsilon>0$  there exsit $n\in \mathbb{N}$ and $y \in B(x, \varepsilon)$  satisfying
	$$
	\inf _{\sigma \in S_{n}} \frac{1}{n} \sum_{k=1}^{n} d\left(T^{k} x, T^{\sigma(k)} y\right)>\delta.$$
	
{\bf Proposition 5.1 }
Let  $(X, T)$  be a  topologicai dynamical system. Then the following conditions are equivalent:

(1)  $(X, T)$  is strong mean sensitive;

(2) there exists  $\delta>0$  such that for every nonempty open subset  $U $ of  $X$  there are  $x, y \in U$  satisfying
$$
	\limsup _{n \rightarrow+\infty} \inf _{\sigma \in S_{n}} \frac{1}{n} \sum_{k=1}^{n} d\left(T^{k} x, T^{\sigma(k)} y\right)>\delta.
$$

(3) there exists  $\eta>0$  such that for every  $x \in X $, the set
$$
D_{\eta}(x)=\left\{y \in X: \limsup _{n \rightarrow+\infty} \inf _{\sigma \in S_{n}} \frac{1}{n} \sum_{k=1}^{n} d\left(T^{k} x, T^{\sigma(k)} y\right) \geq \eta\right\}
$$
is a dense  Borel subset of  $X$;

(4) there exists  $\eta>0$  such that
$$
D_{\eta}=\left\{(x, y) \in X \times X: \limsup _{n \rightarrow+\infty} \inf _{\sigma \in S_{n}} \frac{1}{n} \sum_{k=1}^{n} d\left(T^{k} x, T^{\sigma(k)} y\right) \geq \eta\right\}
$$
is a dense Borel subset of  $X \times X .$

{\bf Proof } (1)  $\Rightarrow$  (4) Since  $(X, T)$  is strong mean sensitive, there exists  $\delta>0$  such that for every  $x \in X $ and every  $\varepsilon>0$  there is  $y \in B(x, \varepsilon)$  satisfying
$$
\limsup _{n \rightarrow+\infty} \inf _{\sigma \in S_{n}} \frac{1}{n} \sum_{k=1}^{n} d\left(T^{k} x, T^{\sigma(k)} y\right)>\delta
$$
Let  $\eta=\frac{\delta}{2}$. Since
$$
D_{\eta}=\bigcap_{m=1}^{\infty} \bigcap_{\ell=1}^{\infty}\left(\bigcup_{n \geq \ell}\bigcap_{\sigma \in S_{n}}\left\{(x, y) \in X \times X: \frac{1}{n} \sum_{k=1}^{n} d\left(T^{k} x, T^{\sigma(k)} y\right)>\eta-\frac{1}{m}\right\}\right) .
$$
Then  $D_{\eta}$  is a Borel  subset of  $X \times X$. If $ D_{\eta} $ is not dense in  $X \times X$, then there exist two nonempty open subsets  $U$  and  $V$  of  $X$  such that
$$
U \times V \subset\left\{(x, y) \in X \times X: \limsup _{n \rightarrow+\infty} \inf _{\sigma \in S_{n}} \frac{1}{n} \sum_{k=1}^{n} d\left(T^{k} x, T^{\sigma(k)} y\right)<\eta\right\}.$$
Pick  $x \in U$, $z \in V$  and  $\varepsilon>0 $ such that  $B(x, \varepsilon) \subset U$. Then for every  $y \in B(x, \varepsilon)$, one has
$$
\begin{aligned}
	\limsup _{n \rightarrow+\infty} \inf _{\sigma \in S_{n}} \frac{1}{n} \sum_{k=1}^{n} d\left(T^{k} x, T^{\sigma(k)} y\right) 
	& \leq \limsup _{n \rightarrow+\infty} \inf _{\sigma \in S_{n}} \frac{1}{n} \sum_{k=1}^{n} d\left(T^{k} x, T^{\sigma(k)} z\right)+\limsup _{n \rightarrow+\infty} \inf _{\sigma \in S_{n}} \frac{1}{n} \sum_{k=1}^{n} d\left(T^{k} z, T^{\sigma(k)}  y\right) \\
	& = \limsup _{n \rightarrow+\infty} \inf _{\sigma \in S_{n}} \frac{1}{n} \sum_{k=1}^{n} d\left(T^{k} x, T^{\sigma(k)} z\right)+\limsup _{n \rightarrow+\infty} \inf _{\sigma \in S_{n}} \frac{1}{n} \sum_{k=1}^{n} d\left(T^{k} y, T^{\sigma(k)}  z\right) \\
	& <\eta+\eta=\delta,
\end{aligned}
$$
which is a contradiction. Thus  $D_{\eta} $ is a dense  Borel  subset of  $X \times X$.

(4)  $\Rightarrow$  (3) Assume that there exists  $\eta>0$  such that  $D_{2 \eta}$  is a dense  a Borel  subset of  $X \times X$. Notice that  $D_{\eta}(x)$  is a Borel subset of $ X$. Thus it just indicate that  $D_{\eta}(x)$  is dense. For every a nonempty open subset $ U$  of  $X$, there exist two points  $y, z \in U$  with  $(y, z) \in D_{2 \eta}$, that is
$$
\limsup _{n \rightarrow+\infty} \inf _{\sigma \in S_{n}} \frac{1}{n} \sum_{k=1}^{n} d(T^{k} y, T^{\sigma(k)}z) >2 \eta
$$
Then
$$
\begin{aligned}
	2 \eta<\limsup _{n \rightarrow+\infty} \inf _{\sigma \in S_{n}} \frac{1}{n} \sum_{k=1}^{n} d(T^{k} y, T^{\sigma(k)}z) 
	& \leq \limsup _{n \rightarrow+\infty} \inf _{\sigma \in S_{n}} \frac{1}{n} \sum_{k=1}^{n} d\left(T^{k} y, T^{\sigma(k)} x\right)+\limsup _{n \rightarrow+\infty} \inf _{\sigma \in S_{n}} \frac{1}{n} \sum_{k=1}^{n} d\left(T^{k} x, T^{\sigma(k)} z\right)\\
	& = \limsup _{n \rightarrow+\infty} \inf _{\sigma \in S_{n}} \frac{1}{n} \sum_{k=1}^{n} d\left(T^{k} x, T^{\sigma(k)} y\right)+\limsup _{n \rightarrow+\infty} \inf _{\sigma \in S_{n}} \frac{1}{n} \sum_{k=1}^{n} d\left(T^{k} x, T^{\sigma(k)} z\right)  .
\end{aligned}
$$
Thus, either
$$
\limsup _{n \rightarrow+\infty} \inf _{\sigma \in S_{n}} \frac{1}{n} \sum_{k=1}^{n} d\left(T^{k} x, T^{\sigma(k)} y\right)>\eta
$$
or
$$
\limsup _{n \rightarrow+\infty} \inf _{\sigma \in S_{n}} \frac{1}{n} \sum_{k=1}^{n} d\left(T^{k} x, T^{\sigma(k)} z\right)>\eta,
$$
this implies either  $y \in D_{\eta}(x) $ or  $z \in D_{\eta}(x)$.

(3)  $\Rightarrow$  (2) is obvious.

(2)  $\Rightarrow$  (1) For any  $x \in X $ and $ \varepsilon>0 $, there are  $y, z \in B(x, \varepsilon)$  satisfying
$$
\limsup _{n \rightarrow+\infty} \inf _{\sigma \in S_{n}} \frac{1}{n} \sum_{k=1}^{n} d\left(T^{k} y, T^{\sigma(k)} z\right)>\delta .
$$
Then either
$$
\limsup _{n \rightarrow+\infty} \inf _{\sigma \in S_{n}} \frac{1}{n} \sum_{k=1}^{n} d\left(T^{k} x, T^{\sigma(k)} y\right)>\frac{\delta}{2}
$$
or
$$
\limsup _{n \rightarrow+\infty} \inf _{\sigma \in S_{n}} \frac{1}{n} \sum_{k=1}^{n} d\left(T^{k} x, T^{\sigma(k)}z\right)>\frac{\delta}{2},
$$
which implies that  $(X, T)$  is strong mean sensitive. 

A point  $x \in X$  is strong mean sensitive if there exists  $\delta>0$  such that for every  $\varepsilon>0$  there is  $y \in B(x, \varepsilon)$  satisfying
$$
\limsup _{n \rightarrow+\infty} \inf _{\sigma \in S_{n}} \frac{1}{n} \sum_{k=1}^{n} d\left(T^{k} x, T^{\sigma(k)} y\right)>\delta.$$

A point  $x \in X$  is strong sensitive in the mean if there exists  $\delta>0$  such that for  every  $\varepsilon>0$  there exsit $n\in \mathbb{N}$ and $y \in B(x, \varepsilon)$  satisfying
$$
\inf _{\sigma \in S_{n}} \frac{1}{n} \sum_{k=1}^{n} d\left(T^{k} x, T^{\sigma(k)} y\right)>\delta.$$
It is clear that a point $x\in X$ is either weakly equicontinuous in the mean or strong sensitive in the mean and a point $x\in X$ is either weakly mean equicontinuous  or strong mean sensitive. 

{\bf Remark 5.2 } Notice that every strong mean sensitive topological dynamical system is mean sensitive, every strong mean sensitive topological dynamical system is strong sensitive in the mean. In additional, a topological dynamical system is strong mean sensitive if and only if there exists $\varepsilon>0$ satisfying $\mathscr{E}_{\varepsilon}$ is empty; a topological dynamical system is strong sensitive in the mean if and only if there exists $\varepsilon>0$ satisfying $\mathscr{E}^{*}_{\varepsilon}$ is empty.

{\bf Proposition 5.3}
Let  $(X, T)$  be a topological  dynamical system. Then  $(X, T)$  is strong sensitive in the mean if and only if there is  $\delta>0$  such that for any  nonempty open subset  $U $ of  $X$ and any  $N \in \mathbb{N}$  there are  $m \geq N$  and   $x,y \in U$  such that
$$ \inf _{\sigma \in S_{m}}\frac{1}{m} \sum_{k=1}^{m} d\left(T^{k} x, T^{\sigma(k)} y\right) > \delta.
$$

{\bf Proof}
Sufficiency is obvious. Now assume that  $(X, T)$  is  strong sensitive in the mean. If necessity is not true, then for any  $\delta>0$  there are a nonempty open subset  $U $ of  $X$ and   $N \in \mathbb{N}$  such that for any  $m \geq N$  and any  $x,y\in U$  we have
$$ \inf _{\sigma \in S_{m}}\frac{1}{m} \sum_{k=1}^{m} d\left(T^{k} x, T^{\sigma(k)} y\right) \leq \delta.
$$

Since  $U\subset \overline{U}$,and $\overline{U}$  is compact and  $T$  is continuous, then there are finite open sets  $U_{1}, \ldots, U_{k} \subset X$ such that $\overline{U}\subset \bigcup_{j=1}^{k}U_j$ and  for any $1\leq j\leq k$, $\max_{1\leq i\leq N} d(T^{i}y,T^{i}z)\leq \delta$ holds for all $y, z \in U_{j} $. Since  $U\subset \overline{U}\subset\bigcup_{j=1}^{k}U_j$, there is  $j_{0} \in\{1, \ldots, k\} $ such that $U \cap U_{j_{0}}\neq \emptyset$ and $U \cap U_{j_{0}}$ is an open set. 
Set $V=U \cap U_{j_{0}}$, for any $x,y\in V$, any $m \geq 1$, we have 
$$\inf _{\sigma \in S_{m}}\frac{1}{m} \sum_{k=1}^{m} d\left(T^{k} x, T^{\sigma(k)} y\right) \leq \delta,$$
this leads to a contradiction with the assumption that  $(X, T)$  is strong sensitive in the mean. The proof is completed.

{\bf Proposition 5.4 } Let  $(X, T)$  be a transitive dynamical system. If there exists a transitive point which is  strong mean sensitive, then $ (X, T)$  is strong mean sensitive.

{\bf Proof} Let  $x \in X$  be a strong mean sensitive transitive point, then there exists  $\delta>0$  such that for every  $\varepsilon>0$  there is  $y \in B(x, \varepsilon)$  satisfying
$$
\limsup _{n \rightarrow+\infty} \inf _{\sigma \in S_{n}} \frac{1}{n} \sum_{k=1}^{n} d\left(T^{k} x, T^{\sigma(k)} y\right)>\delta .
$$
Fix a nonempty open subset $ U$  of  $X$. With  $x$  is a transitive point, there exist  $\varepsilon>0$  and  $i \in \mathbb{Z}_{+}$ such that  $T^{i}(B(x, \varepsilon)) \subset U$. There exists  $y \in B(x, \varepsilon)$  satisfying
$$
\limsup _{n \rightarrow+\infty} \inf _{\sigma \in S_{n}} \frac{1}{n} \sum_{k=1}^{n} d\left(T^{k} x, T^{\sigma(k)} y\right)>\delta .
$$

Let $ u=T^{i} x$  and  $v=T^{i} y $. Then  $u, v \in U$  and
$$
\begin{aligned}
	\limsup _{n \rightarrow+\infty} \inf _{\sigma \in S_{n}} \frac{1}{n} \sum_{k=1}^{n} d\left(T^{k} u, T^{\sigma(k)} v\right) &= \limsup _{n \rightarrow+\infty} \inf _{\sigma \in S_{n}} \frac{1}{n} \sum_{k=1}^{n}d\left(T^{k}\left(T^{i} x\right), T^{\sigma(k)}\left(T^{i} y\right)\right) \\
	& =\limsup _{n \rightarrow+\infty} \inf _{\sigma \in S_{n}} \frac{1}{n} \sum_{k=1}^{n} d\left(T^{k} x, T^{\sigma(k)} y\right)>\delta .
\end{aligned}
$$
Therefore,  $(X, T)$  is strong mean sensitive. The proof is completed.

{\bf Theorem 5.5} Let  $(X, T)$  be a topological dynamical system. If  $(X, T)$  is transitive, then  $(X, T)$  is either almost weakly mean equicontimuous or strong mean sensitive.

{\bf Proof} Let  $x \in X$  be a transitive point. If  $x$  is  strong mean sensitive, then  $(X, T)$  is strong mean sensitive. If $ x$  is not  strong mean sensitive, then it is weakly mean equicontinuous, hence, $(X, T)$  is almost  weakly mean equicontinuous.

{\bf Corollary 5.6} Let  $(X, T)$  be a minimal dynamical system. Then  $(X, T)$  is either weakly mean equicontimuous or strong mean sensitive.

{\bf Theorem 5.7} Let  $(X, T)$  be a topological dynamical system. If  $(X, T)$  is transitive, then  $(X, T)$  is either almost weakly  equicontimuous in  the mean or strong sensitive in the mean.

{\bf Proof}  Assume $(  X, T)$  is not almost  weakly equicontinuous in the mean. Then there is a transitive point  $x \in X$  which is not weakly equicontinuous in the mean. Let $ U $ be a nonempty open subset of  $X$. Then there are  $k \in \mathbb{N}$  and  $\epsilon>0 $ such that  $T^{k} B\left(x, \epsilon\right) \subset U $. Since  $x$  is not a weakly equicontinuous in the mean point, then there is a  $\delta>0$  such that for any  $r \in \mathbb{N}$  with  $0<1 / r<\epsilon
 $  there are  $n_{r} \in \mathbb{N}$  and $ y_{r} \in X$  with  $d\left(x, y_{r}\right)<1 /r$  such that
$$
\inf _{\sigma \in S_{n_r}}\frac{1}{n_{r}} \sum_{i=1}^{n_{r}} d\left(T^{i} x, T^{\sigma(i)} y_{r}\right)>\delta.$$

We can require $ n_{r} \rightarrow \infty$  when  $r \rightarrow \infty$. Otherwise, there is $ n \in \mathbb{N}$  such that  $n=n_{r}$  for infinitely many  $r$. Since  $X$  is compact, then  $\max _{0 \leq i \leq n} d\left(T^{i} x, T^{i} y_{r}\right)<\delta $ whenever $ d\left(x, y_{r}\right)<1 / r $ for some large enough  $r>0$.
For above fixed  $k$, pick large enough  $n_{r} \in \mathbb{N}$  such that 
$$
\sup_{\sigma \in S_{k}}\frac{1}{n_{r}} \sum_{i=1}^{k} d\left(T^{i} x, T^{\sigma(i)} y_{r}\right) \leq \frac{k \operatorname{diam}(X)}{n_{r}}<\frac{\delta}{3} .
$$
Denote $ u=T^{k} x $ and  $v_{r}=T^{k} y_{r} $. Then we have  $u, v_{r} \in U$  and 
$$
\begin{aligned}
	\inf _{\sigma \in S_{n_r}}\frac{1}{n_{r}} \sum_{i=1}^{n_{r}} d\left(T^{i} u, T^{\sigma(i)} v_{r}\right) & =\inf _{\sigma \in S_{n_r}}\frac{1}{n_{r}} \sum_{i=1}^{n_{r}} d\left(T^{i} T^{k}x, T^{\sigma(i)}T^{k}y_r\right) \\
	& \geq \inf _{\sigma \in S_{n_r}}\frac{1}{n_{r}} \sum_{i=1}^{n_{r}} d\left(T^{i} x, T^{\sigma(i)} y_{r}\right)-\sup_{\sigma \in S_{k}}\frac{1}{n_{r}} \sum_{i=1}^{k} d\left(T^{i} x, T^{\sigma(i)} y_{r}\right) \\
	& >\delta-\frac{\delta}{3}=\frac{2\delta}{3} .
\end{aligned}
$$
This implies that  $(X, T)$  is strong sensitive in the mean.

{\bf Corollary 5.8} Let  $(X, T)$  be a minimal dynamical system. Then  $(X, T)$  is either weakly equicontimuous  in the mean or strong sensitive in the mean.

{\bf Proof} It follows from Proposition 4.4 and Theorem 5.7.

{\bf Theorem 5.9} Let  $(X, T)$  be a minimal dynamical system. Then  $(X, T)$ is strong sensitive in the mean if and only if $(X, T)$ is strong mean sensitive.

{\bf Proof} It follows from Remark 4.6 and Corollary 5.8.

By Theorem 5.9, we have that for a minimal system, the equivalence between  strong mean sensitivity and  strong sensitivity in the mean holds,  while we have found that above conclusion remains true even for general dynamical system (see Appendix Theorem 1 and Theorem 2). And then we get  an equivalence between  almost weakly mean equicontinuity and weakly equicontinuity in the mean for a transitive dynamical system.

{\bf Theorem 5.10} Let  $(X, T)$  be a topological dynamical system. If there exist at least two $T$-invariant ergodic Borel probability measures with full support, then  $(X, T)$  is strong mean sensitive.

{\bf Proof} Assume that  $(X, T)$  is not strong mean  sensitive. Let  $\mu$  and  $\nu$  be two ergodic Borel probability measures with full support. Assume that  $x\in X$  is a generic point for  $\mu$, and  $y$  is a generic point for  $\nu$. Since  $x$  and  $y$  are generic points of fully supported measures, and their orbits are dense in  $X$. In particular,  $(X, T)$  is transitive, then according to Theorem 5.7, we have that  $(X, T)$  is almost weakly mean equicontinuous. Furthermore,  we deduce that  $x$  and  $y$  are both weakly mean equicontinuous points and  $\overline{F}\left(x, y\right)=0$. Hence, we have  $\mu=\nu$ .

{\bf Corollary 5.11} Let  $(X, T)$  be a topological dynamical system. If there exist at least two $T$-invariant ergodic Borel probability measures with full support, then  $(X, T)$  is strong sensitive in the mean and mean sensitive.

{\bf Proposition 5.12} Let  $\pi:(X, T) \rightarrow(Y, S)$  be a factor map between two dynamical systems. If  $y \in Y $ is a  strong mean sensitive point, then there exists a  strong mean sensitive point  $x \in \pi^{-1}(y)$.

{\bf Proof} By Theorem 3.5, we can  obtain it.

{\bf Theorem 5.13} Let  $\pi:(X, T) \rightarrow(Y, S)$  be a factor map between two dynamical systems. Assume that  $y \in Y$  is a  strong mean sensitive point and  $x \in \pi^{-1}(y)$. If  $\pi$  is open at  $x$, then $ x$  is also a  strong mean sensitive point.

{\bf Proof} By Theorem 4.7, we can obtain it.

{\bf Theorem 5.14} Let  $\pi:(X, T) \rightarrow(Y, S)$  be a factor map between two transitive dynamical systems. Assume that  $(Y, S)$  is strong mean sensitive. If  $\pi$  is semi-open, then  $(X, T)$  is strong mean sensitive.

{\bf Proof} By Theorem 4.8, we can  obtain it.

{\bf Proposition 5.15} Let  $\pi:(X, T) \rightarrow(Y, S)$  be a factor map between two transitive dynamical systems. If  $(Y, S)$  is strong mean sensitive, then there is a strong mean sensitive subsystem  $(Z, T)$  of  $(X, T) .$

{\bf Proof}  By Zorn's Lemma and the compactness of  $X$, we can find a subsystem  $(Z, T)$  of  $(X, T)$  such that  $\pi(Z)=Y$  and  $Z$  is minimal with respect to this property. Let  $y \in Y$  be a transitive point which is also a strong mean sensitive point. According to Proposition 5.12 there is a strong mean sensitive point  $x \in(Z, T)$  with  $\pi(x)=y$. Let  $Z^{\prime}=\omega(x, T) $. Then  $\left(Z^{\prime}, T\right)$  is a subsystem of  $(Z, T)$  with  $\pi\left(Z^{\prime}\right)=Y$. The minimality of $ Z  $ implies that  $Z^{\prime}=Z$. Then  $x$  is a transitive point of  $(Z, T)$  and  $(Z, T) $ is strong mean sensitive.

{\bf Theorem 5.16} Let  $(X, T)$  be an invertible minimal topological dynamical system. If there exist  $x, y \in X$  and  $\delta>0$  such that
$$
\liminf _{k \rightarrow+\infty}  d\left(T^{-k} x, T^{-k} y\right)=0 
$$
and
$$
\limsup _{n \rightarrow \infty} \inf _{\sigma \in S_{n}}\frac{1}{n} \sum_{i=1}^{n} \rho\left(T^{i} x, T^{\sigma(i)} y\right)>\delta,
$$
then  $(X, T)$  is strong mean sensitive or strong sensitive in the mean.

{\bf Proof} Assume that  $x, y \in X$  and  $\delta>0$  satisfy the conditions as stated. Let  $U$  be a nonempty open subset of  $X$. Since  $X$  is minimal, there is  $l \in \mathbb{N}$  such that  $T^{l} x \in U $. Set  $x_1=T^{l} x$, $y_1=T^{l} y$. It is easy to see that  $x_{1}$,  $y_{1}$  and  $\delta>0$  also satisfy the stated two conditions. As  $U$  is a neighborhood of  $x_{1}$, we can find  $\epsilon>0$  such that  $B\left(x_{1},2 \epsilon\right) \subset U$  with  $\epsilon<\delta / 10$. Set  $V=B\left(x_{1},\epsilon\right) $. On the one hand,  since $\liminf _{k \rightarrow+\infty}  d\left(T^{-k} x_1, T^{-k} y_1\right)=0 $, then we have  that  the set
$$
\left\{k<0:d\left(T^{k} x_{1}, T^{k} y_{1}\right)<\epsilon\right\}
$$
is a thick set  in  $\mathbb{Z}_{-}$. On the other hand,  $(X, T)$  is minimal implies  $\left(X, T^{-1}\right)$  is minimal, so the set
$$
\left\{k<0: T^{k} x_{1} \in V\right\}
$$
is syndetic  in  $\mathbb{Z}_{-}$. Hence, there exists $s<0$  such that
$$
T^{s} x_{1} \in V \text { and }  d\left(T^{s} x_{1}, T^{s} y_{1}\right)<\epsilon.
$$
Therefore,  $T^{s} x_{1},T^{s} y_{1} \in B(V, \epsilon) \subset U$. Set  $x_{2}=T^{s} x_{1}$, $y_{2}=T^{s} y_{1}$. It is easy to see that
$$
\begin{aligned}
\limsup _{n \rightarrow \infty} \inf _{\sigma \in S_{n}}\frac{1}{n} \sum_{i=1}^{n} \rho\left(T^{i} x_2, T^{\sigma(i)} y_2\right)&=\limsup _{n \rightarrow \infty} \inf _{\sigma \in S_{n}}\frac{1}{n} \sum_{i=1}^{n} \rho\left(T^{i} (T^{s+l}x), T^{\sigma(i)}(T^{s+l}y) \right)\\
&=\limsup _{n \rightarrow \infty} \inf _{\sigma \in S_{n}}\frac{1}{n} \sum_{i=1}^{n} \rho\left(T^{i} x, T^{\sigma(i)}y \right)\\
&>\delta.
\end{aligned}
$$ 
This means that  $(X, T)$  is strong mean sensitive, and by Theorem 5.9, then $(X, T)$  is also strong sensitive in the mean. The proof is completed.

{\bf Theorem 5.17} Let  $(X, T)$  and  $(Y, S)$  be two topological dynamical systems. If  $(X, T)$  or  $(Y, S)$  is strong mean sensitive, then $(X \times Y, T \times S)$  is strong mean sensitive.

{\bf Proof}  Without loss of generality, we assume that  $(X, T)$   is strong mean sensitive. For any $\epsilon>0$ and $(x,y)\in X\times Y$, we can choose $0<\epsilon_1<\epsilon$ satisfying $B(x,\epsilon_1)\times B(y,\epsilon_1)\subset B((x,y),\epsilon)$. Since $(X, T)$  is strong mean sensitive, there exists  $\delta>0$ and $u\in B(x,\epsilon_1)$ such that 
$$\limsup _{n \rightarrow \infty} \inf _{\sigma \in S_{n}}\frac{1}{n} \sum_{i=1}^{n}d \left(T^{i}\left(x)\right), T^{\sigma(i)} \left(u\right)\right)>\delta.$$
Let $u\in B(x,\epsilon_1) $ and $v\in B(y,\epsilon_1)$, then $(u,v)\in B((x,y),\epsilon)$. 
$$
\limsup _{n \rightarrow \infty} \inf _{\sigma \in S_{n}}\frac{1}{n} \sum_{i=1}^{n} \rho\left((T\times S)^{i} (x,y), (T\times S)^{\sigma(i)} (u,v)\right)$$
$$
\geq \limsup _{n \rightarrow \infty} \inf _{\sigma \in S_{n}}\frac{1}{n} \sum_{i=1}^{n} d\left(T^{i} x, T^{\sigma(i)} u\right)
>\delta.
$$
This implies  $(X \times Y, T \times S)$  is strong mean sensitive. The proof is completed.

\section{ Mean tuples}
Motivated by the localized idea of sensitivity, Ye and Zhang [28] introduced the notion of sensitive tuples and showed that a transitive topological  dynamical system is sensitive if and only if it has a sensitive tuple. Following the idea of [28], many researchers [8,9,29] further introduced some notions of sensitive tuples like mean sensitive  tuples and sensitive in the mean tuples,  they showed that a transitive dynamical system is mean sensitive if and only if it admits a mean sensitive tuple; a transitive dynamical system is sensitive in the mean if and only if it admits a sensitive in the mean tuple. In this section, we will introduce concepts of  several new version sensitive tuples and furthermore study the properties of these sensitive tuples.

{\bf lemma 6.1}\ Let  $(X, T)$  be a  topological  dynamical system defined on a compact metric space  $(X, d)$. Then  $(X, T)$ is strong mean sensitive if and only if  there exists $\eta>0$  such that for
 any nonempty open set $U\in X$,  there exist $x_1, x_2 \in U$ satisfying 
$$\limsup _{n \rightarrow+\infty} \frac{1}{n} \inf _{\sigma \in S_{n}}\#\{1\leq i \leq n:d\left(T^{i} x_1, T^{\sigma(i)} x_2\right)>\eta
\}>\eta.$$

{\bf Proof\  }   The sufficiency is obvious and we only need to state the necessity. Assume that  $(X, T)$ is strong mean sensitive  with a  sensitive constant  $\delta>0 $. Without loss of generality, we suppose that $\operatorname{diam}(X)=1$. Take  $\eta={\delta}/{2}$. Then, for any nonempty open subset $U\subset X$, there exist  pairwise distinct points $x_1,x_2\in U$ such that
$$
\limsup _{n \rightarrow+\infty} \inf _{\sigma \in S_{n}} \frac{1}{n} \sum_{k=1}^{n} d\left(T^{k} x_1, T^{\sigma(k)} x_2\right)>\delta.
$$
Hence,  there exists an increasing sequence  $\left\{n_{k}\right\}_{k=1}^{\infty} \subset \mathbb{N}$  such that for every  $k \in \mathbb{N}$,
$$
\inf _{\sigma \in S_{n_{k}}}\frac{1}{n_{k}} \sum_{i=1}^{n_{k}}d\left(T^{k} x_1, T^{\sigma(k)} x_2\right) > \delta .
$$
This implies that for every  $k \in \mathbb{N}$,one has
$$
\begin{aligned}
	\delta &<\inf _{\sigma \in S_{n_{k}}}\frac{1}{n_{k}} \sum_{i=1}^{n_{k}}d\left(T^{k} x_1, T^{\sigma(k)} x_2\right)\\
	&\leq  \inf _{\sigma \in S_{n_{k}}}\frac{\#\left\{1 \leq i\leq n_{k}: d\left(T^{k} x_1, T^{\sigma(k)} x_2\right)\leq\eta\right\}}{n_{k}} \eta \\
	&+\inf _{\sigma \in S_{n_{k}}}\frac{\#\left\{1 \leq i\leq n_{k}: d\left(T^{k} x_1, T^{\sigma(k)} x_2\right)>\eta\right\}}{n_{k}}\operatorname{diam}(X) \\
	&\leq \eta+\inf _{\sigma \in S_{n_{k}}}\frac{\#\left\{1 \leq i\leq n_{k}: d\left(T^{k} x_1, T^{\sigma(k)} x_2\right)>\eta\right\}}{n_{k}}.
\end{aligned}
$$
Therefore,
$$
\begin{aligned} 
&\limsup _{n \rightarrow+\infty} \frac{1}{n} \inf _{\sigma \in S_{n}}\#\{1\leq i \leq n:d\left(T^{i} x_1, T^{\sigma(i)} x_2\right)\\
	&\geq \limsup _{k \rightarrow \infty} \inf _{\sigma \in S_{n_{k}}}\frac{\#\left\{1 \leq i\leq n_{k}: d\left(T^{k} x_1, T^{\sigma(k)} x_2\right)>\eta\right\}}{n_{k}}\\
	&> \frac{\delta}{2}=\eta .
\end{aligned}
$$
The proof is completed.

{\bf Definitions 6.2}\ \ Let  $(X, T)$  be a  topological  dynamical system defined on a compact metric space  $(X, d)$. We say that $(x_1, x_2) \in X\times X$ is a $\overline{F}$-mean sensitive tuple if $x_1 \neq x_2$ and for any $\epsilon>0$, there is $c = c(\epsilon ) > 0$ such
that for any nonempty open set $U \subset X $ there exist $y_1, y_2\in U$
such that $\limsup _{n \rightarrow+\infty} \frac{1}{n} \inf _{\sigma \in S_{n}}\#  \{1\leq i \leq n:T^{i}(y_1)\in B(x_1,\epsilon), T^{\sigma(i)}(y_2)\in B(x_2,\epsilon)
\}>c $. We denote the set of all $\overline{F}$-mean sensitive tuples by $MS_{\overline{F}}(X,T)$.

{\bf Theorem  6.3}\ \   Let  $(X, T)$  be a  topological  dynamical system defined on a compact metric space  $(X, d)$. If $(X, T)$  is transitive, then  $(X, T)$ is strong mean sensitive if and only if $MS_{\overline{F}}(X,T)\neq \emptyset$.

{\bf Proof\  }
Assume that $(x_1, x_2) \in MS_{\overline{F}}(X,T)$.
Let $\delta= \frac{1}{2} d(x_1,x_2).$
Then we select some small $\epsilon$ such that $d(B(x_1,\epsilon),B(x_2,\epsilon))>\delta$. By the definition of $\overline{F}$-mean sensitive tuple, for any $\epsilon>0$, there is $c>0$ such that for any nonempty open set $U\in X$ there exist $y_1, y_2 \in U$ such that $$\limsup _{n \rightarrow+\infty} \frac{1}{n} \inf _{\sigma \in S_{n}}\#\{1\leq i \leq n:T^{i}(y_1)\in B(x_1,\epsilon), T^{\sigma(i)}(y_2)\in B(x_2,\epsilon)
\}>c .$$
This deduce that 
$$\limsup _{n \rightarrow+\infty} \frac{1}{n} \inf _{\sigma \in S_{n}} \sum_{k=1}^{n} d\left(T^{k} y_1, T^{\sigma(k)} y_2\right)>c\delta.$$
which yields that  $(X, T)$ is strong mean sensitive.

Conversely, assume that $(X, T)$ is strong mean sensitive. Then according to Lemma 11.1, there exists $\delta >0$,
for any any nonempty open set $U\in X$, there exist $x_1, x_2 \in U$ such that  
$$\limsup _{n \rightarrow+\infty} \frac{1}{n} \inf _{\sigma \in S_{n}}\#\{1\leq i \leq n:d\left(T^{i} x_1, T^{\sigma(i)} x_2\right)>\delta
\}>\delta.$$

Let
$
X_{\delta}=\{(x_1, x_2)\in X\times X: d(x_1,x_2)\geq \delta  \}
$
and it is clear that $X_{\delta}$ is compact in $X\times X$.
Take $x\in X$ as a transitive point. For each $m \in \mathbb{N}$,  we set $W_m=B(x,1/m)$, then there are 
$y_m^1, y_m^2\in W_m$ satisfying
$$\limsup _{n \rightarrow+\infty} \frac{1}{n} \inf _{\sigma \in S_{n}}\#\{1\leq i \leq n:\left(T^{i} y_m^1, T^{\sigma(i)} y_m^2\right)\in X_{\delta}
\}>\delta.$$

Since $X_{\delta}$ is compact, then we can cover $X_{\delta}$ with finite nonempty open sets of 
diameter less than 1, i.e., $X_{\delta}\subset \cup_{i=1}^{N_1}A_1^i$ and $\operatorname{diam}(A_1^i)<1$.
Then for each $m\in \mathbb{N}$ there exits $1\leq N_1^m \leq N_1 $ such that 
$$\limsup _{n \rightarrow+\infty} \frac{1}{n} \inf _{\sigma \in S_{n}}\#\{1\leq i \leq n:\left(T^{i} y_m^1, T^{\sigma(i)} y_m^2\right) \in X_{\delta}\cap \overline { A_1^{N_1^m }} \}>\delta/N_1.$$
Without loss of generality we assume $N_1^m=1$ for all $m\in \mathbb{N}$ and then 
$$\limsup _{n \rightarrow+\infty} \frac{1}{n} \inf _{\sigma \in S_{n}}\#\{1\leq i \leq n:\left(T^{i} y_m^1, T^{\sigma(i)} y_m^2\right) \in X_{\delta}\cap \overline { A_1^1} \}>\delta/N_1.$$

Repeating the procedure above, for $l\geq 1$ we can cover
$X_{\delta} \cap \overline { A_l^{1}} $ with finite nonempty sets of diameter less than $1/(l+1)$, i.e., 
$X_{\delta} \cap \overline { A_l^{1}}\subset \cup_{i=1}^{N_l+1}A_{l+1}^i$ and $\operatorname{diam}(A_{l+1}^i)<1/(l+1)$.
Then for any $m\in \mathbb{N}$ there exists 
$1\leq N_{l+1}^m \leq N_{l+1}$ such that 
$$\limsup _{n \rightarrow+\infty} \frac{1}{n} \inf _{\sigma \in S_{n}}\#\{1\leq i \leq n:\left(T^{i} y_m^1, T^{\sigma(i)} y_m^2\right) \in X_{\delta}\cap \overline { A_{l+1}^{N_{l+1}^m}} \}>\frac{\delta}{N_1\cdots N_{l+1}}.$$

Without loss of generality we assume $N_{l+1}^m=1$ for all $m\in \mathbb{N}$ and then 
$$\limsup _{n \rightarrow+\infty} \frac{1}{n} \inf _{\sigma \in S_{n}}\#\{1\leq i \leq n:\left(T^{i} y_m^1, T^{\sigma(i)} y_m^2\right) \in X_{\delta}\cap \overline { A_{l+1}^{1}} \}>\frac{\delta}{N_1\cdots N_{l+1}}.$$

It is clear than there is a unique point 
$(z_1,z_2)\in  \bigcap_{l=1}^{\infty}\overline{A_l^1} \cap X_{\delta}$.
We will indicate that $(z_1,z_2)\in MS_{\overline{F}}(X,T)$.
In fact, for any $\epsilon>0$ there is $l\in \mathbb{N}$
such that $\overline{A_l^1} \cap X_{\delta} \subset V_1\times V_2$, where $V_i=B(z_i,\epsilon)$ for $i=1,2$.
By the construction, for any $W_m$ there are $y_m^1, y_m^2\in W_m$ such that 
$$\limsup _{n \rightarrow+\infty} \frac{1}{n} \inf _{\sigma \in S_{n}}\#\{1\leq i \leq n:\left(T^{i} y_m^1, T^{\sigma(i)} y_m^2\right) \in X_{\delta}\cap \overline { A_{l}^{1}} \}>\frac{\delta}{N_1\cdots N_{l}}.$$
and furthermore 
$$\limsup _{n \rightarrow+\infty} \frac{1}{n} \inf _{\sigma \in S_{n}}\#\{1\leq i \leq n:\left(T^{i} y_m^1, T^{\sigma(i)} y_m^2\right) \in V_1\times V_2 \}>\frac{\delta}{N_1\cdots N_{l}}.$$
for all $m\in \mathbb{N}$. For any nonempty open set $U\in X$,
as $x$ is a transitive point, there exist $r\in \mathbb{N}$ such that $T^r(x)\in U$. We can pick $m_1\in \mathbb{N}$ such that 
$T^r(W_{m_1})\subset U $, which implies that 
$T^ry_{m_1}^1,T^ry_{m_1}^2\in U$ and 
$$\limsup _{n \rightarrow+\infty} \frac{1}{n} \inf _{\sigma \in S_{n}}\#\{1\leq i \leq n:\left(T^{i} (T^ry_{m_1}^1), T^{\sigma(i)}(T^ry_{m_1}^2) \right) \in V_1\times V_2 \}>\frac{\delta}{N_1\cdots N_{l}}.$$
Hence, we derive that $(z_1,z_2)\in MS_{\overline{F}}(X,T)$. The  proof is completed.

{\bf Definitions 6.4}\ \ Let  $(X, T)$  be a  topological  dynamical system defined on a compact metric space  $(X, d)$. We say that $(x_1, x_2) \in X\times X$ is a $\overline{F}$-sensitive in the mean tuple if $x_1 \neq x_2$ and for any $\epsilon>0$, there is $c = c(\epsilon ) > 0$ such
that for any nonempty open set $U \subset X $ there exist $m\in \mathbb{N}$ and $y_1^m, y_2^m\in U$
such that $$ \frac{1}{m} \inf _{\sigma \in S_{m}}\#\{1\leq i \leq m:T^{i}(y_1^m)\in B(x_1,\epsilon), T^{\sigma(i)}(y_2^m)\in B(x_2,\epsilon)\}>c .$$ We denote the set of all $\overline{F}$-sensitive in the mean tuples by $SM_{\overline{F}}(X,T)$.

{\bf Theorem  6.5}\ \   Let  $(X, T)$  be a  topological  dynamical system defined on a compact metric space $(X, d)$. If $(X, T)$  is transitive, then  $(X, T)$ is strong sensitive in the mean if and only if $SM_{\overline{F}}(X,T)\neq \emptyset$.

Proof it is similar to Theorem 6.3, so we omit it.

{\bf Definitions 6.6}\ \ Let  $(X, T)$  be a topological dynamical system. and  $2 \leq n \in \mathbb{N}$ . We say that  $(X, T)$  is weakly  $\overline{F}$-sensitive in the mean if there is $ \delta>0$  such that for any $\varepsilon>0$  there are  $m \in \mathbb{N}$  and pairwise distinct points  $x_{1}^{m}, x_{2}^{m}$with  $d(x_{1}^{m}, x_{2}^{m})<\varepsilon$  such that
$$
\inf _{\sigma \in S_{m}} \frac{1}{m} \sum_{k=1}^{m} d\left(T^{k} x_{1}^{m}, T^{\sigma(k)} x_{2}^{m}\right)>\delta.$$

{\bf Lemma 6.7}\ \ Let  $(X, T)$  be a topological dynamical system. The following are equivalent:

(1)  $(X, T)$  is weakly $\overline{F}$-sensitive in the mean.

(2) There is  $\delta>0 $ such that for any  $\varepsilon>0 $ and any $ N \in \mathbb{N}$  there are  $m \geq N$  and pairwise distinct points  $x_{1}^{m}, x_{2}^{m}$  with  $d(x_{1}^{m}, x_{2}^{m})<\varepsilon$  such that
$$
\inf _{\sigma \in S_{m}} \frac{1}{m} \sum_{k=1}^{m} d\left(T^{k} x_{1}^{m}, T^{\sigma(k)}x_{2}^{m}\right)>\delta.$$

(3) There is  $\delta>0$  such that for any  $\ell \in \mathbb{N} $ there are  $m_{\ell} \in \mathbb{N}$  and pairwise distinct points $ x_{1}^{m_{\ell}}, x_{2}^{m_{\ell}} $  with  $d(x_{1}^{m_{\ell}}, x_{2}^{m_{\ell}})<1 / \ell$  such that
$$
\lim _{\ell \rightarrow \infty}\inf _{\sigma \in S_{m_{\ell}}} \frac{1}{m_{\ell}} \sum_{k=1}^{m_{\ell}} d\left(T^{k} x_{1}^{m_{\ell}}, T^{\sigma(k)} x_{2}^{m_{\ell}}\right)>\delta .
$$

{\bf Proof\  }  It is easy to check that (3)  $\Leftrightarrow$(2) $\Rightarrow$  (1). If (2) does not hold, then for any  $\delta>0$  there are  $\varepsilon>0 $ and  $N \in \mathbb{N}$  such that for any  $m \geq N $ and any  $x_{1}, x_{2} $ with  $d(x_{1}, x_{2})<\varepsilon$  such that
$$
\inf _{\sigma \in S_{m}} \frac{1}{m} \sum_{k=1}^{m} d\left(T^{k} x_{1}, T^{\sigma(k)}x_{2}\right) \leq \delta .\ \ \ \ \text{(I)}$$

Since $ X$  is compact, we can ask  $\varepsilon$  small enough such that (I) holds for all  $m \geq 1$. So (1) can not be true at the same time. That is,  (1) $\Rightarrow$(2). The  proof is completed.

From the proof of Lemma 6.7, we can easily have the following proposition:

{\bf Proposition 6.8}\ \  $(X, T)$  is weakly mean equicontinuous if and only if it is not weakly $\overline{F}$-sensitive in the mean.

{\bf Definitions 6.9}\ \ Let  $(X, T)$  be a  topological  dynamical system defined on a compact metric space  $(X, d)$. We say that $(x_1, x_2) \in X\times X$ is a weakly $\overline{F}$-sensitive in the mean tuple if $x_1 \neq x_2$ and for any $\epsilon>0$, there is $c = c(\epsilon ) > 0$ such
that for any $\tau>0$ there exist $m\in \mathbb{N}$ and $y_1, y_2\in X$ with $d(y_1, y_2)<\tau$
such that $$ \frac{1}{m} \inf _{\sigma \in S_{m}}\#\{1\leq i \leq m:T^{i}(y_1)\in B(x_1,\epsilon), T^{\sigma(i)}(y_2)\in B(x_2,\epsilon)\}>c .$$ We denote the set of all weakly $\overline{F}$-sensitive in the mean tuples by $WSM_{\overline{F}}(X,T)$.

{\bf Theorem  6.10}\ \   Let  $(X, T)$  be a  topological  dynamical system defined on a compact metric space  $(X, d)$. Then  $(X, T)$ is weakly $\overline{F}$-sensitive in the mean if and only if $WSM_{\overline{F}}(X,T)\neq \emptyset$.

{\bf Proof\  }  Assume that  $\left(x_{1}, x_{2}\right) \in WSM_{\overline{F}}(X,T)$. Put  $\tau= d\left(x_{1}, x_{2}\right)/ 4 $. Then there is  $c=c(\tau)>0$  such that for any  $\varepsilon>0 $ there are  $m \in \mathbb{N}$  and distinct points  $y_{1}, y_{2}$  in  $X$  with  $d\left(y_{1}, y_{2}\right)<\varepsilon$  such that
$$
 \frac{1}{m} \inf _{\sigma \in S_{m}}\#\{1\leq i \leq m:T^{i}(y_1)\in B(x_1,\tau), T^{\sigma(i)}(y_2)\in B(x_2,\tau)\}>c.
$$ 
By the choice of  $\tau $, we have  $d\left(B\left(x_{1}, \tau\right), B\left(x_{2}, \tau\right)\right)>\tau $. This implies that
$$
\inf _{\sigma \in S_{m}} \frac{1}{m} \sum_{k=1}^{m} d\left(T^{k} x_{1}, T^{\sigma(k)} x_{2}\right)>\tau c,
$$
hence,  $(X, T)$  is weakly $\overline{F}$-sensitive in the mean with constant  $\tau c $.

Conversely, we assume  $(X, T) $ is weakly $\overline{F}$-sensitive in the mean. Without loss of generality, we assume that  $\operatorname{diam}(X) \leq 1$. Then there is $ \delta>0 $ such that for any  $\ell \in \mathbb{N} $ there are  $m_{\ell} \in \mathbb{N}$  and  $n$  distinct points  $x_{1}^{m_{\ell}}, x_{2}^{m_{\ell}}$  in  $X$  with  $d\left(x_{1}^{m_{\ell}}, x_{2}^{m_{\ell}}\right)<\varepsilon$  such that
$$
\inf _{\sigma \in S_{m_{\ell}}}\frac{1}{m_{\ell}} \sum_{k=1}^{m_{\ell}}d\left(T^{k} x_{1}^{m_{\ell}}, T^{\sigma(k)} x_{2}^{m_{\ell}}\right)>2 \delta .
$$

Let  $X_{\delta}=\left\{(y_{1}, y_{2}) \in X\times X: d\left(y_{1}, y_{2}\right) \geq \delta\right\}$. For each $ \ell \in \mathbb{N} $, we have
$$\inf _{\sigma \in S_{m_{\ell}}}\frac{1}{m_{\ell}} \sum_{k=1}^{m_{\ell}}d\left(T^{k} x_{1}^{m_{\ell}}, T^{\sigma(k)} x_{2}^{m_{\ell}}\right)\leq \inf _{\sigma \in S_{m_{\ell}}}\frac{\#\left\{1 \leq k \leq m_{\ell}:\left(T^{k} x_{1}^{m_{\ell}},  T^{\sigma(k)}x_{2}^{m_{\ell}}\right) \in X_{\delta}^{n}\right\}}{m_{\ell}}+\delta,$$

and so
$$
 \inf _{\sigma \in S_{m_{\ell}}}\frac{\#\left\{1 \leq k \leq m_{\ell}:\left(T^{k} x_{1}^{m_{\ell}},  T^{\sigma(k)}x_{2}^{m_{\ell}}\right) \in X_{\delta}^{n}\right\}}{m_{\ell}}>\delta.
$$

Since  $X_{\delta}$  is compact, we can pick a finite open cover with diameter less than 1, that is, there are  $N_{1} \in \mathbb{N}$  and  $A_{1}^{1}, \ldots, A_{1}^{N_{1}}$  are open in  $X\times X$  such that  $X_{\delta} \subset \cup_{i=1}^{N_{1}} A_{1}^{i} $ and  $\operatorname{diam}\left(A_{1}^{i}\right)<1 $. Then there is  $1 \leq N_{1}^{\ell} \leq N_{1} $ such that
$$
\inf _{\sigma \in S_{m_{\ell}}}\frac{\#\left\{1 \leq k \leq m_{\ell}:\left(T^{k} x_{1}^{m_{\ell}},T^{\sigma(k)}  x_{2}^{m_{\ell}}\right) \in A_{1}^{N_{1}^{\ell}} \cap X_{\delta}\right\}}{m_{\ell}}>\frac{\delta}{N_{1}} .
$$

Note that  $N_{1}$  is independent of $ \ell$  and $ N_{1}^{\ell} $ has only  $N_{1}$  choices. Since the above procedure holds for all  $\ell \in \mathbb{N} $, without loss of generality we can assume  $N_{1}^{\ell}=1$  for all  $\ell \in \mathbb{N} $.

Repeating the above process, for  $l \geq 1$, we can cover  $\overline{A_{l}^{1}} \cap X_{\delta}$ by finite nonempty open sets with diameter less than  $1 /(l+1) $, that is,  $\overline{A_{l}^{1}} \cap X_{\delta} \subset \cup_{i=1}^{N_{l+1}} A_{l+1}^{i}$  and $ \operatorname{diam}\left(A_{l+1}^{i}\right)<1 /(l+1)$  for $1 \leq i \leq N_{l+1} $. Then there is  $1 \leq N_{l+1}^{\ell} \leq N_{l+1} $ such that
$$
\inf _{\sigma \in S_{m_{\ell}}}\frac{\#\left\{1 \leq k \leq m_{\ell}:\left(T^{k} x_{1}^{m_{\ell}},T^{\sigma(k)}  x_{2}^{m_{\ell}}\right) \in A_{l+1}^{N_{l+1}^{\ell}} \cap X_{\delta}\right\}}{m_{\ell}}>\frac{\delta}{N_{1} \cdots N_{l+1}} .
$$
By the same reason we can assume  $N_{l+1}^{\ell}=1 $ for all  $\ell \in \mathbb{N}$ .

It is clear that there is a unique tuple  $\left(z_{1}, z_{2}\right) \in \cap_{l=1}^{\infty} \overline{A_{l}^{1}} \cap X_{\delta}$. Next we show that  $\left(z_{1}, z_{2}\right) \in WSM_{\overline{F}}(X,T)$. Obviously  $z_{1} \neq z_{2}$. For any  $\tau>0 $, by construction there is  $p=p(\tau) \in \mathbb{N}$  such that  $A_{p}^{1} \cap X_{\delta}^{n} \subset \prod_{i=1}^{2} B\left(z_{i}, \tau\right) $. Let  $\varepsilon>0 $ and then there is $ q \in \mathbb{N} $ such that $ 1 / q<\varepsilon $. For the fixed $ q $, from the above discussion there are  $m_{q} \in \mathbb{N}$  and  $n$  distinct points  $x_{1}^{m_{q}}, x_{2}^{m_{q}} $ in  $X $ with  $d\left(x_{1}^{m_{q}}, x_{2}^{m_{q}}\right)<\varepsilon$  such that
$$
\inf _{\sigma \in S_{m_{q}}}\frac{\#\left\{1 \leq k \leq m_{q}:\left(T^{k} x_{1}^{m_{q}}, T^{\sigma(k)} x_{2}^{m_{q}}\right) \in A_{p}^{1} \cap X_{\delta}^{n}\right\}}{m_{q}}>\frac{\delta}{N_{1} \cdots N_{p}} .
$$
This implies that
$$
\inf _{\sigma \in S_{m_{q}}}\frac{\# \{1 \leq k \leq m_{q}:T^{k} x_{1}^{m_{q}}\in B(z_1,\tau), T^{\sigma(k)} x_{2}^{m_{q}}\in B(z_2,\tau)) \}}{m_{q}}>\frac{\delta}{N_{1} \cdots N_{p}}
$$
and so  $\left(z_{1}, z_{2}\right) \in WSM_{\overline{F}}(X,T)$. The  proof is completed.

Combining Proposition 6.8 with Theorem  6.10 we have the result as follows:

{ \bf Corollary 6.11 }$(X, T)$  is weakly mean equicontinuous if and only if  $WSM_{\overline{F}}(X,T)=\emptyset $.
 
{\bf Proposition 6.12}\ \  Let  $(X, T)$  be a minimal topological dynamical system. Then  $WSM_{\overline{F}}(X,T)= SM_{\overline{F}}(X, T)$.

{\bf Proof\  }  Obviously  $SM_{\overline{F}}(X, T) \subset WSM_{\overline{F}}(X, T)$. It suffices to show $ WSM_{\overline{F}}(X, T) \subset SM_{\overline{F}}(X, T)$. 
Take  $\left(x_{1}, x_{2}\right) \in WSM_{\overline{F}}(X, T) $. Then for any  $\tau>0$  there is  $c=c(\tau)>0 $ such
 that for any  $\ell \in \mathbb{N} $ there are  $m_{\ell} \in \mathbb{N}(m_{\ell} \rightarrow \infty$ 
when  $\ell \rightarrow   \infty)$  and pairwise distinct points  $y_{1}^{m_{\ell}}, y_{2}^{m_{\ell}} \in X$  with  $d\left(y_{1}^{m_{\ell}}, y_{2}^{m_{\ell}}\right)<1 / \ell $ such that
$$
\inf _{\sigma \in S_{m_{\ell}}}\frac{\# \{1 \leq k \leq m_{\ell}:T^{k} y_{1}^{m_{\ell}}\in B(x_1,\tau), T^{\sigma(k)} y_{2}^{m_{\ell}}\in B(x_2,\tau)) \}}{m_{\ell}}>c.
$$

We can choose $ \left\{m_{\ell}: \ell \in \mathbb{N}\right\}$  to be an increasing sequence. Assume  $\lim _{\ell \rightarrow \infty} y_{1}^{m_{\ell}}=y $ for some $ y \in X$. For any nonempty open subset  $U$  of  $X$, we choose a nonempty open subset  $V \subset X $ and some $ \ell_{1} \in \mathbb{N}$  such that  $V \subset \overline{V} \subset B\left(\overline{V}, 1 / \ell_{1}\right) \subset U $. Since  $X $ is minimal, there is  $k \in \mathbb{N}$  such that  $T^{k} y \in V $. By uniform continuity of  $T^{k} $, there is  $\delta>0 $ such that if  $d\left(z_{1}, z_{2}\right)<\delta $ then  $d\left(T^{k} z_{1}, T^{k} z_{2}\right)<1 /\left(2 \ell_{1}\right) $. Choose $ \ell_{2} \in \mathbb{N} $ such that  $d\left(y_{1}^{m_{\ell}}, y\right)<\delta$  for all $ \ell \geq \ell_{2} $ and  $1 / \ell_{2}<\delta $. Then
$$
d\left(T^{k} y, T^{k} y_{i}^{m_{\ell}}\right) \leq d\left(T^{k} y, T^{k} y_{1}^{m_{\ell}}\right)+d\left(T^{k} y_{1}^{m_{\ell}}, T^{k} y_{i}^{m_{\ell}}\right)<1 / \ell_{1}
$$
for all  $\ell \geq \ell_{2}$  and all  $1 \leq i \leq 2 $. This implies that
$$
\left\{T^{k} y_{i}^{m_{\ell}}: 1 \leq i \leq 2\right\} \subset B\left(T^{k} y, 1 / \ell_{1}\right) \subset B\left(\overline{V}, 1 / \ell_{1}\right) \subset U
$$
and so  $\left(T^{k}\times T^{k}\right) \left(y_{1}^{m_{\ell}},y_{2}^{m_{\ell}}\right) \in U\times U $ for all $ \ell \geq \ell_{2} $.

For above  $k$ , let  $\ell_{3} \in \mathbb{N} $ be such that
$$
\sup _{\sigma \in S_{m_{\ell}}}\frac{\#\left\{1 \leq j \leq k: T^{j} y_{1}^{m_{\ell}} \in B\left(x_{1}, \tau\right),T^{\sigma(j)} y_{2}^{m_{\ell}} \in B\left(x_{2}, \tau\right)\right\}}{m_{\ell}} \leq \frac{k}{m_{\ell}}<\frac{c}{3}
$$
for all  $\ell \geq \ell_{3}$. Now for each  $\ell \geq \max \left\{\ell_{2}, \ell_{3}\right\}$  we have  $T^{k} y_{1}^{m_{\ell}}, T^{k} y_{2}^{m_{\ell}} \in U $ and
$$
\begin{aligned}
	&\inf_{\sigma \in S_{m_{\ell}}}\frac{1}{m_{\ell}} \#\{1 \leq j \leq m_{\ell}: T^{j}(T^{k} y_{1}^{m_{\ell}}) \in B(x_{1}, \tau), T^{\sigma(j)}(T^{k} y_{2}^{m_{\ell}}) \in B(x_{2}, \tau)\} \\
	&\geq \inf_{\sigma \in S_{m_{\ell}}}\frac{\#\{1 \leq j \leq m_{\ell}: T^{j} y_{1}^{m_{\ell}} \in B(x_{1}, \tau),T^{\sigma(j)} y_{2}^{m_{\ell}} \in B(x_{2}, \tau)\}}{m_{\ell}}\\
	&-\sup_{\sigma \in S_{m_{\ell}}}\frac{\#\{1 \leq j \leq k: T^{j} y_{1}^{m_{\ell}} \in B(x_{1}, \tau),T^{\sigma(j)} y_{2}^{m_{\ell}} \in B(x_{2}, \tau)\}}{m_{\ell}}\\
	&>c-\frac{c}{3}=\frac{2c}{3} .
\end{aligned}
$$

This implies that  $\left(x_{1}, x_{2}\right) \in SM_{\overline{F}}(X, T) $. The  proof is completed.

{\bf Corollary 6.13}\ \  Let  $(X, T)$  be a minimal topological dynamical system. Then  $SM_{\overline{F}}(X,T)\neq\emptyset$ if and only if $MS_{\overline{F}}(X, T)\neq\emptyset$ if and only if $WSM_{\overline{F}}(X,T)\neq\emptyset$.

{\bf Theorem 6.14}\ \  Let  $(X, T)$  be a minimal topological dynamical system. Then  $(X, T)$ is strong sensitive in the mean if and only if $(X, T)$ is weakly  $\overline{F}$-sensitive in the mean.

{\bf Proof\  }  Obviously  necessity is clear. It suffices to show the sufficiency. Assume $(X, T)$ is weakly  $\overline{F}$-sensitive in the mean. Then  there is  $\delta>0 $ such
that for any  $\ell \in \mathbb{N} $ there are  $m_{\ell} \in \mathbb{N}(m_{\ell} \rightarrow \infty$ 
when  $\ell \rightarrow   \infty)$  and pairwise distinct points  $y_{1}^{m_{\ell}}, y_{2}^{m_{\ell}} \in X$  with  $d\left(y_{1}^{m_{\ell}}, y_{2}^{m_{\ell}}\right)<1 / \ell $ such that
$$
\inf _{\sigma \in S_{m_{\ell}}} \frac{1}{m_{\ell}} \sum_{k=1}^{m_{\ell}} d\left(T^{k} y_{1}^{m_{\ell}}, T^{\sigma(k)}y_{2}^{m_{\ell}}\right)>\delta.$$

We can choose $ \left\{m_{\ell}: \ell \in \mathbb{N}\right\}$  to be an increasing sequence. Assume  $\lim _{\ell \rightarrow \infty} y_{1}^{m_{\ell}}=y $ for some $ y \in X$. For any nonempty open subset  $U$  of  $X$, we choose a nonempty open subset  $V \subset X $ and some $ \ell_{1} \in \mathbb{N}$  such that  $V \subset \overline{V} \subset B\left(\overline{V}, 1 / \ell_{1}\right) \subset U $. Since  $X $ is minimal, there is  $k \in \mathbb{N}$  such that  $T^{k} y \in V $. By uniform continuity of  $T^{k} $, there is  $\delta>0 $ such that if  $d\left(z_{1}, z_{2}\right)<\delta $ then  $d\left(T^{k} z_{1}, T^{k} z_{2}\right)<1 /\left(2 \ell_{1}\right) $. Choose $ \ell_{2} \in \mathbb{N} $ such that  $d\left(y_{1}^{m_{\ell}}, y\right)<\delta$  for all $ \ell \geq \ell_{2} $ and  $1 / \ell_{2}<\delta $. Then
$$
d\left(T^{k} y, T^{k} y_{i}^{m_{\ell}}\right) \leq d\left(T^{k} y, T^{k} y_{1}^{m_{\ell}}\right)+d\left(T^{k} y_{1}^{m_{\ell}}, T^{k} y_{i}^{m_{\ell}}\right)<1 / \ell_{1}
$$
for all  $\ell \geq \ell_{2}$  and all  $1 \leq i \leq 2 $. This implies that
$$
\left\{T^{k} y_{i}^{m_{\ell}}: 1 \leq i \leq 2\right\} \subset B\left(T^{k} y, 1 / \ell_{1}\right) \subset B\left(\overline{V}, 1 / \ell_{1}\right) \subset U
$$
and so  $\left(T^{k}\times T^{k}\right) \left(y_{1}^{m_{\ell}},y_{2}^{m_{\ell}}\right) \in U\times U $ for all $ \ell \geq \ell_{2} $.

For above  $k$ , let  $\ell_{3} \in \mathbb{N} $ be such that
$$
\sup _{\sigma \in S_{m_{\ell}}}\frac{1}{m_{\ell}} \sum_{r=1}^{k} d\left(T^{r} y_{1}^{m_{\ell}}, T^{\sigma(r)}y_{2}^{m_{\ell}}\right) \leq \frac{k\operatorname{diam}(X)}{m_{\ell}}<\frac{\delta}{3}
$$
for all  $\ell \geq \ell_{3}$. Now for each  $\ell \geq \max \left\{\ell_{2}, \ell_{3}\right\}$  we have  $T^{k} y_{1}^{m_{\ell}}, T^{k} y_{2}^{m_{\ell}} \in U $ and
$$
\begin{aligned}
	&\inf _{\sigma \in S_{m_{\ell}}} \frac{1}{m_{\ell}} \sum_{r=1}^{m_{\ell}} d\left(T^{r} T^{k} y_{1}^{m_{\ell}}, T^{\sigma(r)} T^{k}y_{2}^{m_{\ell}}\right) \\
	&\geq \inf _{\sigma \in S_{m_{\ell}}} \frac{1}{m_{\ell}} \sum_{r=1}^{m_{\ell}} d\left(T^{r} y_{1}^{m_{\ell}}, T^{\sigma(r)}y_{2}^{m_{\ell}}\right)
	-\sup _{\sigma \in S_{m_{\ell}}}\frac{1}{m_{\ell}} \sum_{r=1}^{k} d\left(T^{r} y_{1}^{m_{\ell}}, T^{\sigma(r)}y_{2}^{m_{\ell}}\right)\\
	&>\delta-\frac{\delta}{3}=\frac{2\delta}{3} .
\end{aligned}
$$

This implies that  $(X, T)$ is strong sensitive in the mean. The  proof is completed.

{\bf Theorem 6.15}\ \  Let  $\pi:(X, T) \rightarrow(Y, S)$  be a factor map between two topological dynamical systems. Then,

(1)  $\pi\times \pi \left(M S_{\overline{F}}(X, T)\right) \subset M S_{\overline{F}}(Y, S) \cup \Delta(Y)$;

(2)  $\pi\times \pi \left(M S_{\overline{F}}(X, T) \cup \Delta(X)\right)=M S_{\overline{F}}(Y, S) \cup \Delta(Y)$, provided that  $(X, T)$  is minimal.

{\bf Proof\  } Item (1) is easy to be proved by the definition. We only prove item (2).
Supposing that  $\left(y_{1}, y_{2}\right) \in M S_{\overline{F}}(Y, S) $, we will show that there exists  $\left(z_{1}, z_{2}\right) \in M S_{\overline{F}}(X, T)$  such that  $\pi\left(z_{i}\right)=y_{i}$  for $i=1,2 $. Fix  $x \in X$  and let  $U_{m}=B(x, 1 / m) $. Since  $(X, T)$  is minimal,  $\operatorname{int}\left(\pi\left(U_{m}\right)\right) \neq \emptyset$, where  $\operatorname{int}\left(\pi\left(U_{m}\right)\right)$  is the interior of  $\pi\left(U_{m}\right) $. Since  $\left(y_{1}, y_{2}\right) \in M S_{\overline{F}}(Y, S) $, there exists  $\delta>0$  and  $y_{m}^{1}, y_{m}^{2} \in \operatorname{int}\left(\pi\left(U_{m}\right)\right)$  such that
$$\limsup _{n \rightarrow+\infty} \frac{1}{n} \inf _{\sigma \in S_{n}}\#\{1\leq i \leq n:S^{i}(y_{m}^{1})\in \overline{B\left(y_{1}, 1\right)}, S^{\sigma(i)}(y_{m}^{2})\in \overline{B\left(y_{2}, 1\right)}
\}>\delta.$$
Then there exist  $x_{m}^{1}, x_{m}^{2} \in U_{m}$  with  $\pi\left(x_{m}^{i}\right)=y_{m}^{i}$ for $i=1,2$  such that for any  $m \in \mathbb{N} $,
$$\limsup _{n \rightarrow+\infty} \frac{1}{n} \inf _{\sigma \in S_{n}}\#\{1\leq i \leq n:T^{i}(x_{m}^{1})\in \pi^{-1}(\overline{B\left(y_{1}, 1\right)}), T^{\sigma(i)}(x_{m}^{2})\in \pi^{-1}(\overline{B\left(y_{2}, 1\right)})
\}> \delta.$$
Put
$$
A=\pi^{-1}\left(\overline{B\left(y_{1}, 1\right)}\right)\times \pi^{-1}\left(\overline{B\left(y_{2}, 1\right)}\right) ,
$$
and it is clear that  $A$  is a compact subset of $X\times X $. 

We can cover  $A$  with finite nonempty open sets of diameter less than $1 $, that is,  $A \subset \bigcup_{i=1}^{N_{1}} A_{1}^{i} $ and  $\operatorname{diam}\left(A_{1}^{i}\right)<1 $. Then for each  $m \in \mathbb{N} $, there is  $1 \leq N_{1}^{m} \leq N_{1}$  such that
$$\limsup _{n \rightarrow+\infty} \frac{1}{n} \inf _{\sigma \in S_{n}}\#\{1\leq i \leq n:\left(T^{i} x_m^1, T^{\sigma(i)} x_m^2\right) \in  \overline{A_{1}^{N_{1}^{m}}} \cap A \}> \delta/N_1.$$

Without loss of generality, we assume $ N_{1}^{m}=1 $ for all  $m \in \mathbb{N} $. Namely,
$$\limsup _{n \rightarrow+\infty} \frac{1}{n} \inf _{\sigma \in S_{n}}\#\{1\leq i \leq n:\left(T^{i} x_m^1, T^{\sigma(i)} x_m^2\right) \in  \overline{A_{1}^{1}} \cap A \}> \delta/N_1 \quad \text { for all } m \in \mathbb{N} .$$

Repeating the above procedure, for  $l \geq 1 $, we can cover  $\overline{A_{l}^{1}} \cap A $ with finite nonempty open sets of diameter less than  $1 /(l+1) $, that is,  $\overline{A_{l}^{1}} \cap A \subset \bigcup_{i=1}^{N_{l+1}} A_{l+1}^{i} $ and  $\operatorname{diam}\left(A_{l+1}^{i}\right)<1 /(l+1) $. Then for each  $m \in \mathbb{N} $, there is  $1 \leq N_{l+1}^{m} \leq N_{l+1}$  such that
$$\limsup _{n \rightarrow+\infty} \frac{1}{n} \inf _{\sigma \in S_{n}}\#\{1\leq i \leq n:\left(T^{i} x_m^1, T^{\sigma(i)} x_m^2\right) \in  \overline{A_{l+1}^{N_{l+1}^{m}}} \cap A \}>\frac{\delta}{N_{1} N_{2} \cdots N_{l+1}} .$$

Without loss of generality, we assume  $N_{l+1}^{m}=1$  for all  $m \in \mathbb{N} $. Namely,
$$\limsup _{n \rightarrow+\infty} \frac{1}{n} \inf _{\sigma \in S_{n}}\#\{1\leq i \leq n:\left(T^{i} x_m^1, T^{\sigma(i)} x_m^2\right) \in  \overline{A_{l+1}^{1}} \cap A \}> \frac{\delta}{N_{1} N_{2} \cdots N_{l+1}} \quad \text { for all } m \in \mathbb{N} .$$

It is clear that there is a unique point  $\left(z_{1}^{1}, z_{2}^{1}\right) \in \bigcap_{l=1}^{\infty} \overline{A_{l}^{1}} \cap A $. We claim that  $\left(z_{1}^{1},z_{2}^{1}\right) \in M S_{\overline{F}}(X, T) $. In fact, for any  $\tau>0 $, there is  $l \in \mathbb{N} $ such that  $\overline{A_{l}^{1}} \cap A \subset  V_{1} \times V_{2} $, where  $V_{i}=B\left(z_{i}^{1}, \tau\right)$  for  $i=1,2 $. By the construction, for any $ m \in \mathbb{N}$ , there are  $x_{m}^{1}, x_{m}^{2} \in U_{m}$  such that
$$\limsup _{n \rightarrow+\infty} \frac{1}{n} \inf _{\sigma \in S_{n}}\#\{1\leq i \leq n:\left(T^{i} x_m^1, T^{\sigma(i)} x_m^2\right) \in  \overline{A_{l}^{1}} \cap A \}> \frac{\delta}{N_{1} N_{2} \cdots N_{l}}$$
and so
$$\limsup _{n \rightarrow+\infty} \frac{1}{n} \inf _{\sigma \in S_{n}}\#\{1\leq i \leq n:\left(T^{i} x_m^1, T^{\sigma(i)} x_m^2\right) \in  V_{1} \times V_{2} \}> \frac{\delta}{N_{1} N_{2} \cdots N_{l}}$$
for all  $m \in \mathbb{N} $. For any nonempty open set  $U \subset X $, since $x$  is a transitive point, there is  $s \in \mathbb{Z}$  such that  $T^{s} x \in U$. We can choose  $m \in \mathbb{Z}$  such that  $T^{s} U_{m} \subset U $. This implies that  $T^{s} x_{m}^{1}, T^{s} x_{m}^{2} \in U$  and
$$\limsup _{n \rightarrow+\infty} \frac{1}{n} \inf _{\sigma \in S_{n}}\#\{1\leq i \leq n:\left(T^{i}\left(T^{s} x_{m}^{1}\right), T^{\sigma(i)} \left(T^{s} x_{m}^{2}\right)\right) \in  V_{1} \times V_{2} \}> \frac{\delta}{N_{1} N_{2} \cdots N_{l}}.$$
Hence, we have  $\left(z_{1}^{1}, z_{2}^{1}\right) \in M S_{\overline{F}}(X, T) $.
Similarly, for each  $p \in \mathbb{N} $, there exists  $\left(z_{1}^{p},z_{2}^{p}\right) \in M S_{\overline{F}}(X, T) \cap \prod_{i=1}^{2}   \pi^{-1}\left(\overline{B\left(y_{i}, 1 / p\right)}\right) $. Set  $z_{i}^{p} \rightarrow z_{i} $ as  $p \rightarrow \infty $. Then  $\left(z_{1}, z_{2}\right) \in M S_{\overline{F}}(X, T) \cup \Delta(X) $ and  $\pi\left(z_{i}\right)=y_{i}$ for $i=1,2 $.

{\bf Theorem 6.16}\ \ 
 Let $ \pi:(X, T) \rightarrow(Y, S)$  be a factor map. Then:
 
(1)  $\pi\times \pi \left(WSM_{\overline{F}}(X, T)\right) \subset WSM_{\overline{F}}(Y, S) \cup \Delta(Y) $. That is, if  $\left(x_{1},  x_{2}\right) \in WSM_{\overline{F}}(X, T)$  and $ \left(\pi\left(x_{1}\right), \pi\left(x_{2}\right)\right) \notin \Delta(Y)$, then  $\left(\pi\left(x_{1}\right), \pi\left(x_{2}\right)\right) \in WSM_{\overline{F}}(Y, S)$.

(2) If  $WSM_{\overline{F}}(Y, S) \neq \emptyset $, then  $WSM_{\overline{F}}(X, T) \neq \emptyset$  and there exists  $\left(x_{1}, x_{2}\right) \in WSM_{\overline{F}}(X, T)$ such that  $\left(\pi\left(x_{1}\right), \pi\left(x_{2}\right)\right) \in WSM_{\overline{F}}(Y, S) $.

{\bf Proof\  }  Denote  $d^{\prime} $ and $ \rho$  the metrics on  $X$  and  $Y$, respectively. Define a function  $d$  on  $X \times X $ by  $d(x, y)=d^{\prime}(x, y)+\rho(\pi(x), \pi(y))$  for any $ x, y \in X $. Since  $\pi $ is a factor map, it is clear that $ d$  is a compatible metric on  $X$  and  $\pi $ is Lipschitz $1$ with respect to the metrics  $d$  and  $\rho $, i.e., $ d(x, y) \geq \rho(\pi(x), \pi(y))$  for any  $x, y \in X .$

(1) It follows directly from the definitions and the Lipschitz continuity of  $\pi $.

(2) If $ WSM_{\overline{F}}(Y, S) \neq \emptyset $, then by Corollary 6.11 and the fact that each factor of a weakly mean equicontinuous topological dynamical system is weakly mean equicontinuous we have  $WSM_{\overline{F}}(X, T) \neq \emptyset $. Furthermore, if  $\left(y_{1}, y_{2}\right) \in WSM_{\overline{F}}(Y, S)$ , then for  $a:=\rho\left(y_{1}, y_{2}\right) / 4>0 $, there is $ c=c(a)>0 $ such that for any $ \ell \in \mathbb{N}$  there are  $m_{\ell} \in \mathbb{N} $ and distinct points  $y_{1}^{m_{\ell}}, y_{2}^{m_{\ell}} \in Y $ with  $\rho\left(y_{1}^{m_{\ell}}, y_{2}^{m_{\ell}}\right)<1 / \ell$  such that
$$
\inf_{\sigma \in S_{m_{\ell}}}\frac{\#\left\{1 \leq k \leq m_{\ell}: S^{k} y_{1}^{m_{\ell}} \in B\left(y_{1}, a\right),S^{\sigma(k)} y_{2}^{m_{\ell}} \in B\left(y_{2}, a\right)\right\}}{m_{\ell}}>c .
$$

Put $ \delta=\min \{2 a, c\} / 2 $. Then
$$
\inf_{\sigma \in S_{m_{\ell}}}\frac{\#\left\{0 \leq k \leq m_{\ell}-1: \rho\left(S^{k} y_{1}^{m_{\ell}}, S^{\sigma(k)} y_{2}^{m_{\ell}}\right)>2 \delta\right\}}{m_{\ell}}>2 \delta .
$$

Assume that  $\lim _{\ell \rightarrow \infty} y_{1}^{m_{\ell}}=\lim _{\ell \rightarrow \infty} y_{2}^{m_{\ell}}=\overline{y} $. Without loss of generality, by simple triangle computation we can further assume
$$\inf_{\sigma \in S_{m_{\ell}}}\frac{\#\left\{1 \leq k \leq m_{\ell}: \rho\left(S^{k} y_{1}^{m_{\ell}}, S^{\sigma(k)} \overline{y}\right)>\delta\right\}}{m_{\ell}}>\delta$$
for all  $\ell \in \mathbb{N} $. Let  $x_{1}^{m_{\ell}} \in X$  be such that  $\pi\left(x_{1}^{m_{\ell}}\right)=y_{1}^{m_{\ell}} $. Assume  $\lim _{\ell \rightarrow \infty} x_{1}^{m_{\ell}}=\overline{x} $ for some  $\overline{x} \in X $. It is easy to see that $ \pi(\overline{x})=\overline{y} $. Denote  $E_{m_{\ell}}^{\sigma}=\left\{1 \leq k \leq m_{\ell}: \rho\left(S^{k} y_{1}^{m_{\ell}}, S^{\sigma(k)}\overline{y}\right)>\delta\right\} $. Since  $\pi $ is Lipschitz 1 we have
$$
\begin{aligned}
\inf_{\sigma \in S_{m_{\ell}}}\frac{\#\left\{1 \leq k \leq m_{\ell}: d\left(T^{k} x_{1}^{m_{\ell}}, T^{\sigma(k)} \overline{x}\right)>\delta\right\}}{m_{\ell}} &\geq \inf_{\sigma \in S_{m_{\ell}}}\frac{\#\left\{k \in E_{m_{\ell}}^{\sigma}: d\left(T^{k} x_{1}^{m_{\ell}}, T^{\sigma(k)} \overline{x}\right)>\delta\right\}}{m_{\ell}}\\&\geq \inf_{\sigma \in S_{m_{\ell}}}\frac{\#\left\{1 \leq k \leq m_{\ell}: \rho\left(S^{k} y_{1}^{m_{\ell}}, S^{\sigma(k)} \overline{y}\right)>\delta\right\}}{m_{\ell}}\\&>\delta.
\end{aligned}
$$

Taking the same arguments and following the same notations as Theorem 6.10 we can get a sequence of nonempty open subsets  $\left\{A_{l}^{1}\right\}_{l=1}^{\infty} $ in  X  with  $\operatorname{diam}\left(A_{l}^{1}\right)<\frac{1}{l}$  such that there is a unique pair  $\left(z_{1}, z_{2}\right) \in \cap_{l=1}^{\infty} \overline{A_{l}^{1}} \cap X_{\delta} \cap  WSM_{\overline{F}}(X, T) $, where $ X_{\delta}=\left\{\left(x, y\right) \in X\times X: d\left(x, y\right) \geq\right.   \delta  \}$. Moreover, for each  $\tau>0 $ there is  $l=l(\tau) \in \mathbb{N} $ such that for any  $\ell \in \mathbb{N} $ we have
$$
\inf_{\sigma \in S_{m_{\ell}}}\frac{\#\left\{k \in E_{m_{\ell}}^{\sigma}: d\left(T^{k} x_{1}^{m_{\ell}}, z_{1}\right) \leq 1 /(l+1), d\left(T^{\sigma(k)} \overline{x}, z_{2}\right) \leq 1 /(l+1)\right\}}{m_{\ell}}>\frac{\delta}{N_{1} \cdots N_{l+1}} .
$$

Note that  $l=l(\tau) \rightarrow \infty$  when  $\tau \rightarrow 0 $. Put
$$
E_{m_{\ell}}^{\sigma,l}=\left\{k \in E_{m_{\ell}}^{\sigma}: d\left(T^{k} x_{1}^{m_{\ell}}, z_{1}\right) \leq 1 /(l+1), d\left(T^{\sigma(k)} \overline{x}, z_{2}\right) \leq 1 /(l+1)\right\} .
$$

If  $k \in E_{m_{\ell}}^{\sigma,l} $, by Lipschitz continuity of  $\pi $ we have
$$
\rho\left(S^{k} y_{1}^{m_{\ell}}, \pi\left(z_{1}\right)\right) \leq 1 /(l+1) \text { and } \rho\left(S^{\sigma(k)} \overline{y}, \pi\left(z_{2}\right)\right) \leq 1 /(l+1) .
$$

 Note that  $\rho\left(S^{k} y_{1}^{m_{\ell}}, S^{\sigma(k)}\overline{y}\right) \geq \delta$  since  $k \in E_{m_{\ell}}^{\sigma} $. This implies that  $\rho\left(\pi\left(z_{1}\right), \pi\left(z_{2}\right)\right) \geq \delta-2 /(l+1) $. Let  $l \rightarrow \infty $, we get
$$
\rho\left(\pi\left(z_{1}\right), \pi\left(z_{2}\right)\right)=\lim _{l \rightarrow \infty} [\delta-\frac{2}{l+1}]=\delta>0 .
$$

Then by (1) we have that $ \left(\pi\left(z_{1}\right), \pi\left(z_{2}\right)\right) \in WSM_{\overline{F}}(Y, S)$. The  proof is completed.

\section{ weakly density-equicontinuity and density $\overline{F}$-sensitivity}

In [9], Li and Tu introduced the notions of density-$t$-equicontinuity and density sensitivity, and showed that a topological dynamical system is mean  equicontinuous iff it is density-$t$-equicontinuous for every  $t \in[0,1)$; a minimal topological dynamical system is either density-equicontinuous  or density-sensitive. Motivated by there ideas,  in this section, we will introduce concepts of  weakly density-$t$-equicontinuity and density $\overline{F}$-sensitivity, then we get 
dichotomies between weakly density-$t$-equicontinuity and density $\overline{F}$-sensitivity.

Let  $(X, T)$  be a  topological dynamical system. Fix  $t \in[0,1] $. We say that  $(X, T)$  is weakly density-$t$-equicontinuous if for any $ \epsilon>0 $, there is a  $\delta>0 $ such that if for any  $x, y \in X$  with  $d(x, y)<\delta $, we have $$\limsup_{n \rightarrow \infty} \frac{1}{n} \inf _{\sigma \in S_{n}} \#\{1 \leq i \leq n: d\left(T^{i} x, T^{\sigma(i)}  y\right)>\varepsilon\}\leq 1-t.$$ When  $t=1$ we call it weakly density-equicontinuity.

{\bf Proposition  7.1}\ \
$(X, T) $ is weakly mean  equicontinuous if and only if it is weakly density-$t$-equicontinuous for every  $t \in[0,1) $.

{\bf Proof}\ \ Assume that  $(X, T)$  is weakly mean  equicontinuous. Take  $t \in(0,1] $. By Theorem 3.8, for any  $\epsilon>0$  there exists  $\delta(\epsilon)>0 $ such that if  $d(x, y)<\delta(\epsilon)$  then
$$
\limsup_{n \rightarrow \infty} \frac{1}{n} \inf _{\sigma \in S_{n}} \#\{1 \leq i \leq n: d\left(T^{i} x, T^{\sigma(i)}  y\right)>\epsilon\} \leq \epsilon .
$$

When $ \epsilon \in(0, t]$ we choose  $\delta^{\prime}\left(:=\delta^{\prime}(\epsilon, t)\right)=\delta(\epsilon) $. Then  $d(x, y)<\delta^{\prime} $ implies
$$
\limsup_{n \rightarrow \infty} \frac{1}{n} \inf _{\sigma \in S_{n}} \#\{1 \leq i \leq n: d\left(T^{i} x, T^{\sigma(i)}  y\right)>\epsilon\} \leq \epsilon \leq t .
$$

When  $\epsilon \in(t, 1] $ we choose $ \delta^{\prime}\left(:=\delta^{\prime}(\epsilon, t)\right)=\delta(t) $. Then  $d(x, y)<\delta^{\prime} $ implies
$$
\limsup_{n \rightarrow \infty} \frac{1}{n} \inf _{\sigma \in S_{n}} \#\{1 \leq i \leq n: d\left(T^{i} x, T^{\sigma(i)}  y\right)>\epsilon\} \leq \limsup_{n \rightarrow \infty} \frac{1}{n} \inf _{\sigma \in S_{n}} \#\{1 \leq i \leq n: d\left(T^{i} x, T^{\sigma(i)}  y\right)>t\} \leq t .$$
This implies that, for any $ \epsilon>0 $, if we choose  $\delta^{\prime}\left(:=\delta^{\prime}(\epsilon, t)\right)=\min \{\delta(\epsilon), \delta(t)\}$  then
$$
\limsup_{n \rightarrow \infty} \frac{1}{n} \inf _{\sigma \in S_{n}} \#\{1 \leq i \leq n: d\left(T^{i} x, T^{\sigma(i)}  y\right)>\epsilon\} \leq t,
$$
whenever $ d(x, y)<\delta^{\prime}$, which shows that  $(X, T)$  is weakly density-($1-t$)-equicontinuous. Since $ t $ is arbitrary, we deduce that  $(X, T)$  is weakly density-$t$-equicontinuous for every  $t \in[0,1) $.

Conversely, we assume that  $(X, T)$  is weakly density-$t$-equicontinuous for every  $t \in[0,1)$. Let  $\epsilon>0 $. We choose a  $t \in(1-\epsilon, 1)$. Then for the given  $\epsilon$  and $ t$, there is  $\delta(\epsilon, t)>0 $ such that if  $x, y \in X $ with $ d(x, y)<\delta(\epsilon, t)$,  one has  $$\limsup_{n \rightarrow \infty} \frac{1}{n} \inf _{\sigma \in S_{n}} \#\{1 \leq i \leq n: d\left(T^{i} x, T^{\sigma(i)}  y\right)>\varepsilon\}  \leq 1-t<\epsilon .$$ This implies that  $(X, T) $ is weakly mean equicontinuous,  the proof is completed.

A point  $x \in X $ is a weakly density-equicontinuous point if for any  $\epsilon>0 $, there is a  $\delta>0$  such that if  $y \in X $ with  $d(x, y)<\delta$  then  $\limsup_{n \rightarrow \infty} \frac{1}{n} \inf _{\sigma \in S_{n}} \#\{1 \leq i \leq n: d\left(T^{i} x, T^{\sigma(i)}  y\right)>\varepsilon\}=0$; and  $(X, T)$  is almost weakly density-equicontinuous if there exists a transitive and weakly density-equicontinuous point.
Denote by $ D_{\overline{F}}(X, T) $ the set of all weakly density-equicontinuous points in $ X $.

 For every  $\epsilon>0 $, let
$$
D_{\epsilon}(X, T)=\left\{x \in X: \exists \delta>0, \forall y, z \in B(x, \delta), \limsup_{n \rightarrow \infty} \frac{1}{n} \inf _{\sigma \in S_{n}} \#\{1 \leq i \leq n: d\left(T^{i} y, T^{\sigma(i)}  z\right)>\varepsilon\}=0\right\} .
$$

In the same way as Propositions 4.1, Propositions 4.2, it is easy to get the following result. The proof is left to readers.

{\bf Proposition 7.2}\ \
 Let  $(X, T)$  be a topological dynamical system We have the following facts:

(1) For each  $\epsilon>0, D_{\epsilon}(X, T)$  is open,  $T^{-1} D_{\epsilon}(X, T) \subset D_{\epsilon}(X, T) $ and
$$
D_{\overline{F}}(X, T)=\bigcap_{m=1}^{\infty} D_{\frac{1}{m}}(X, T)
$$
is a  $G_{\delta}$  subset of  X .

(2) If  $(X, T) $ is transitive, then  $D_{\overline{F}}(X, T)$  is either empty or residual. If additionally  $(X, T)$  is almost weakly density-equicontinuous, then every transitive point belongs to $ D_{\overline{F}}(X, T)$.

(3) If  $(X, T)$ is minimal and almost weakly density-equicontinuous, then it is weakly density-equicontinuous.

{\bf Definition 7.3 }\ \ Let  $(X, T) $ be a topological dynamical system. We say that $ (X, T) $ is density $\overline{F}$-sensitive if there is $ \delta>0$  such that for any nonempty open subset  $U $ of  $X$  there are  $x,y \in   U  $ such that the set
$$
\limsup_{n \rightarrow \infty} \frac{1}{n} \inf _{\sigma \in S_{n}} \#\{1 \leq i \leq n: d\left(T^{i} x, T^{\sigma(i)}  y\right)>\delta\}>0.
$$ Such a $ \delta$  is referred to be a sensitive constant.

{\bf Definition 7.4}\ \ We say that the tuple  $\left(x_{1}, x_{2}\right) \in X\times X \backslash \Delta(X) $ is a density $\overline{F}$-sensitive tuple if for any $\epsilon>0$  and any nonempty open set $ U \subset X$  there exist  $y_{1}, y_{2}\in U $ such that the set  $$\limsup _{n \rightarrow+\infty} \frac{1}{n} \inf _{\sigma \in S_{n}}\#\{1\leq i \leq n:T^{i}(y_1)\in B(x_1,\epsilon), T^{\sigma(i)}(y_2)\in B(x_2,\epsilon)
 \}>0.$$
We denote the set of all density $\overline{F}$-sensitive tuples by $DS_{\overline{F}}(X,T)$.

Next, we will discuss the relationship between density $\overline{F}$-sensitivity and weakly density-equicontinuity. We have the following dichotomies as follows.

{\bf Proposition 7.5}\ \
 Let  $(X, T)$  be a topological dynamical system.  Then

(1) if  $(X, T)$ is transitive, then it is either almost weakly density-equicontinuous or density $\overline{F}$-sensitive;

(2) if  $(X, T)$  is minimal, then it is either weakly density-equicontinuous or density $\overline{F}$-sensitive.

{\bf Proof}\ \ (1) Assume that $ (X, T)$  is not almost weakly density-equicontinuous. By Proposition 7.2 there is a transitive point  $x \in X$  which is not weakly density-equicontinuous. Let  $U$  be a nonempty open subset of  $X $. Then there are  $k \in \mathbb{N}$  and  $\epsilon_{0}>0$  such that  $T^{k} B\left(x, \epsilon_{0}\right) \subset U $. Since $ x$  is not weakly density-equicontinuous, then there is a $ \delta>0$  such that for any  $\ell \in \mathbb{N} $ with  $0<1 / \ell<\epsilon_{0}  $ there are  $n_{\ell} \in \mathbb{N}$  and $ y_{\ell} \in X$  with  $d\left(x, y_{\ell}\right)<1 / \ell$  such that
$$
\limsup_{n \rightarrow \infty} \frac{1}{n} \inf _{\sigma \in S_{n}} \#\{1 \leq i \leq n: d\left(T^{i} x, T^{\sigma(i)}  y_{\ell}\right)>\delta\}>0 .
$$

Denote  $u=T^{k} x $ and  $v_{\ell}=T^{k} y_{\ell}$ . Then we have  $u, v_{\ell} \in U$  and
$$
\limsup_{n \rightarrow \infty} \frac{1}{n} \inf _{\sigma \in S_{n}} \#\{1 \leq i \leq n: d\left(T^{i} u, T^{\sigma(i)}  v_{\ell}\right)>\delta\}=\limsup_{n \rightarrow \infty} \frac{1}{n} \inf _{\sigma \in S_{n}} \#\{1 \leq i \leq n: d\left(T^{i} x, T^{\sigma(i)}  y_{\ell}\right)>\delta\}>0
$$
This implies that  $(X, T) $ is density $\overline{F}$-sensitive.

(2) It follows from (1) and Proposition 7.2.

The following proposition reveals that density $\overline{F}$-sensitive can be characterized by density $\overline{F}$-sensitive-tuples.

{\bf Proposition 7.6}\ \ Let  $(X, T)$  be a transitive topological dynamical system. Then  $(X, T)$  is density $\overline{F}$-sensitive if and only if  $DS_{\overline{F}}(X, T) \neq \emptyset $.

{\bf Proof}\ \ Assume that  $\left(x_{1}, x_{2}\right) \in DS_{\overline{F}}(X, T) $. Set  $\delta=\frac{1}{2} d(x_{1}, x_{2})$  and it is clear that  $\delta>0 $. Let  $U_{1}, U_{2}  $ be open neighbourhoods of  $x_{1}, x_{2}$  respectively such that  $d\left(U_{1}, U_{2}\right)>\delta$. Since  $\left(x_{1},x_{2}\right) \in DS_{\overline{F}}(X, T) $, then for any nonempty open subset  $U $ of $ X$  there are $ y_{1}, y_{2}\in U $, such that 
$$\limsup _{n \rightarrow+\infty} \frac{1}{n} \inf _{\sigma \in S_{n}}\#\{1\leq i \leq n:T^{i}(y_1)\in U_{1}, T^{\sigma(i)}(y_2)\in U_{2})
\}>0.$$
This implies that
$$
\limsup_{n \rightarrow \infty} \frac{1}{n} \inf _{\sigma \in S_{n}} \#\{1 \leq i \leq n: d\left(T^{i} y_1, T^{\sigma(i)}  y_{2}\right)>\delta\}>0 .
$$
and so  $(X, T) $ is density $\overline{F}$-sensitive.

Conversely, assume that  $(X, T)$  is density $\overline{F}$-sensitive. Set
$$
X_{\delta}=\{\left(x_{1}, x_{2}\right) \in X\times X: d\left(x_{1}, x_{2}\right) \geq \delta\},
$$
and it is clear that $ X_{\delta}$  is a closed subset of  $X\times X $. Let  $x \in X $ be a transitive point. Since  $(X, T) $ is density $\overline{F}$-sensitive, then for every  $m \in \mathbb{N}$  there are  $x_{m}^{1},x_{m}^{2} \in B(x, 1 / m)$  such that
$$\limsup _{n \rightarrow+\infty} \frac{1}{n} \inf _{\sigma \in S_{n}}\#\{1\leq i \leq n:\left(T^{i} x_m^1, T^{\sigma(i)} x_m^2\right)\in X_{\delta}
\}>0.$$

We can choose an open cover $ \left\{A_{1}^{1}, \ldots, A_{1}^{N_{1}}\right\}$  of  $X_{\delta}$ for some  $N_{1} \in \mathbb{N} $ such that  $\max \left\{\operatorname{diam}\left(A_{1}^{i}\right)\right.  :  i=1, \ldots, n\}<1 $. Then for every $ m \in \mathbb{N} $ there is  $1 \leq N_{1}^{m} \leq N_{1} $ such that
$$\limsup _{n \rightarrow+\infty} \frac{1}{n} \inf _{\sigma \in S_{n}}\#\{1\leq i \leq n:\left(T^{i} x_m^1, T^{\sigma(i)} x_m^2\right) \in X_{\delta}\cap \overline { A_1^{N_1^m }} \}>0.$$
Without loss of generality we assume that  $N_{1}^{m}=1$  for all $ m \in \mathbb{N}$  and
$$\limsup _{n \rightarrow+\infty} \frac{1}{n} \inf _{\sigma \in S_{n}}\#\{1\leq i \leq n:\left(T^{i} x_m^1, T^{\sigma(i)} x_m^2\right) \in X_{\delta}\cap \overline{A_{1}^{1}} \}>0.$$

Repeating the above procedure, for every  $l \geq 1$  we can cover  $\overline{A_{l}^{1}} \cap X_{\delta}^{n} $ by finite nonempty open subsets with diameters less than  $\frac{1}{l+1}$, i.e., $ \overline{A_{l}^{1}} \cap X_{\delta}^{n} \subset \bigcup_{i=1}^{N_{l+1}} A_{l+1}^{i}$  and  $\operatorname{diam}\left(A_{l+1}^{i}\right)< 1 /(l+1) $. Then for every  $m \in \mathbb{N}$  there is  $1 \leq N_{l+1}^{m} \leq N_{l+1}$  such that
$$\limsup _{n \rightarrow+\infty} \frac{1}{n} \inf _{\sigma \in S_{n}}\#\{1\leq i \leq n:\left(T^{i} x_m^1, T^{\sigma(i)} x_m^2\right) \in X_{\delta}\cap \overline { A_{l+1}^{N_{l+1}^m}} \}>0.$$

Without loss of generality we assume that  $N_{l+1}^{m}=1$  for all  $m \in \mathbb{N}$  and
$$\limsup _{n \rightarrow+\infty} \frac{1}{n} \inf _{\sigma \in S_{n}}\#\{1\leq i \leq n:\left(T^{i} x_m^1, T^{\sigma(i)} x_m^2\right) \in X_{\delta}\cap \overline { A_{l+1}^{1}} \}>0.$$
It is clear that there is a unique point  $\left(z_{1}, z_{2}\right) \in \bigcap_{l=1}^{\infty} \overline{A_{l}^{1}} \cap X_{\delta}$. Now we shall show that  $\left(z_{1},z_{2}\right) \in DS_{\overline{F}}(X, T) $. For any  $\epsilon>0 $, there is  $l \in \mathbb{N}$  such that  $\overline{A_{l}^{1}} \cap X_{\delta}\subset V_{1} \times V_{2} $, where $ V_{i}=B\left(z_{i}, \epsilon\right)$  for  $i=1, 2$. By the construction, for every  $B(x, 1 / m)$  there are $ x_{m}^{1}, x_{m}^{2} \in B(x, 1 / m) $ such that
$$\limsup _{n \rightarrow+\infty} \frac{1}{n} \inf _{\sigma \in S_{n}}\#\{1\leq i \leq n:\left(T^{i} x_m^1, T^{\sigma(i)} x_m^2\right) \in X_{\delta}\cap \overline { A_{l}^{1}} \}>0.$$
and so
$$\limsup _{n \rightarrow+\infty} \frac{1}{n} \inf _{\sigma \in S_{n}}\#\{1\leq i \leq n:\left(T^{i} x_m^1, T^{\sigma(i)} x_m^2\right) \in V_{1} \times V_{2} \}>0.$$
for all $ m \in \mathbb{N}$. For every nonempty open set  $U \subset X$, since $ x$  is a transitive point, there is  $s \in \mathbb{Z}_{+} $ such that $ T^{s} x \in U$  and so we can choose  $m_{0} \in \mathbb{N} $ such that  $T^{s} B\left(x, 1 / m_{0}\right) \subset U $. This implies that  $T^{s} x_{m_{0}}^{1}, T^{s} x_{m_{0}}^{2} \in U $ and
$$\limsup _{n \rightarrow+\infty} \frac{1}{n} \inf _{\sigma \in S_{n}}\#\{1\leq i \leq n:\left(T^{i} (T^sx_{m_1}^1), T^{\sigma(i)}(T^sx_{m_1}^2) \right) \in V_1\times V_2 \}>0.$$
Hence $\left(z_{1},z_{2}\right) \in DS_{\overline{F}}(X, T) $, completing the proof.

{\bf Theorem 7.7}\ \  Let  $\pi:(X, T) \rightarrow(Y, S)$  be a factor map between two topological dynamical systems. Then,

(1)  $\pi\times \pi \left(DS_{\overline{F}}(X, T)\right) \subset DS_{\overline{F}}(Y, S) \cup \Delta(Y)$;

(2)  $\pi\times \pi \left(DS_{\overline{F}}(X, T) \cup \Delta(X)\right)=DS_{\overline{F}}(Y, S) \cup \Delta(Y)$, provided that  $(X, T)$  is minimal.

{\bf Proof\  } Item (1) is easy to be proved by the definition. We only prove item (2).
Supposing that  $\left(y_{1}, y_{2}\right) \in DS_{\overline{F}}(Y, S) $, we will show that there exists  $\left(z_{1}, z_{2}\right) \in DS_{\overline{F}}(X, T)$  such that  $\pi\left(z_{i}\right)=y_{i}$  for $i=1,2 $. Fix  $x \in X$  and let  $U_{m}=B(x, 1 / m) $. Since  $(X, T)$  is minimal,  $\operatorname{int}\left(\pi\left(U_{m}\right)\right) \neq \emptyset$, where  $\operatorname{int}\left(\pi\left(U_{m}\right)\right)$  is the interior of  $\pi\left(U_{m}\right) $. Since  $\left(y_{1}, y_{2}\right) \in M S_{\overline{F}}(Y, S) $, there exists  $\delta>0$  and  $y_{m}^{1}, y_{m}^{2} \in \operatorname{int}\left(\pi\left(U_{m}\right)\right)$  such that
$$\limsup _{n \rightarrow+\infty} \frac{1}{n} \inf _{\sigma \in S_{n}}\#\{1\leq i \leq n:S^{i}(y_{m}^{1})\in \overline{B\left(y_{1}, 1\right)}, S^{\sigma(i)}(y_{m}^{2})\in \overline{B\left(y_{2}, 1\right)}
\}> 0.$$
Then there exist  $x_{m}^{1}, x_{m}^{2} \in U_{m}$  with  $\pi\left(x_{m}^{i}\right)=y_{m}^{i}$ for $i=1,2$  such that for any  $m \in \mathbb{N} $,
$$\limsup _{n \rightarrow+\infty} \frac{1}{n} \inf _{\sigma \in S_{n}}\#\{1\leq i \leq n:T^{i}(x_{m}^{1})\in \pi^{-1}(\overline{B\left(y_{1}, 1\right)}), T^{\sigma(i)}(x_{m}^{2})\in \pi^{-1}(\overline{B\left(y_{2}, 1\right)})
\}> 0.$$
Put
$$
A=\pi^{-1}\left(\overline{B\left(y_{1}, 1\right)}\right)\times \pi^{-1}\left(\overline{B\left(y_{2}, 1\right)}\right) ,
$$
and it is clear that  $A$  is a compact subset of $X\times X $. 

We can cover  $A$  with finite nonempty open sets of diameter less than $1 $, that is,  $A \subset \bigcup_{i=1}^{N_{1}} A_{1}^{i} $ and  $\operatorname{diam}\left(A_{1}^{i}\right)<1 $. Then for each  $m \in \mathbb{N} $, there is  $1 \leq N_{1}^{m} \leq N_{1}$  such that
$$\limsup _{n \rightarrow+\infty} \frac{1}{n} \inf _{\sigma \in S_{n}}\#\{1\leq i \leq n:\left(T^{i} x_m^1, T^{\sigma(i)} x_m^2\right) \in  \overline{A_{1}^{N_{1}^{m}}} \cap A \}> 0.$$
Without loss of generality, we assume $ N_{1}^{m}=1 $ for all  $m \in \mathbb{N} $. Namely,
$$\limsup _{n \rightarrow+\infty} \frac{1}{n} \inf _{\sigma \in S_{n}}\#\{1\leq i \leq n:\left(T^{i} x_m^1, T^{\sigma(i)} x_m^2\right) \in  \overline{A_{1}^{1}} \cap A \}> 0\quad \text { for all } m \in \mathbb{N} .$$

Repeating the above procedure, for  $l \geq 1 $, we can cover  $\overline{A_{l}^{1}} \cap A $ with finite nonempty open sets of diameter less than  $1 /(l+1) $, that is,  $\overline{A_{l}^{1}} \cap A \subset \bigcup_{i=1}^{N_{l+1}} A_{l+1}^{i} $ and  $\operatorname{diam}\left(A_{l+1}^{i}\right)<1 /(l+1) $. Then for each  $m \in \mathbb{N} $, there is  $1 \leq N_{l+1}^{m} \leq N_{l+1}$  such that
$$\limsup _{n \rightarrow+\infty} \frac{1}{n} \inf _{\sigma \in S_{n}}\#\{1\leq i \leq n:\left(T^{i} x_m^1, T^{\sigma(i)} x_m^2\right) \in  \overline{A_{l+1}^{N_{l+1}^{m}}} \cap A \}> 0 .$$
Without loss of generality, we assume  $N_{l+1}^{m}=1$  for all  $m \in \mathbb{N} $. Namely,
$$\limsup _{n \rightarrow+\infty} \frac{1}{n} \inf _{\sigma \in S_{n}}\#\{1\leq i \leq n:\left(T^{i} x_m^1, T^{\sigma(i)} x_m^2\right) \in  \overline{A_{l+1}^{1}} \cap A \}> 0 \quad \text { for all } m \in \mathbb{N} .$$

It is clear that there is a unique point  $\left(z_{1}^{1}, z_{2}^{1}\right) \in \bigcap_{l=1}^{\infty} \overline{A_{l}^{1}} \cap A $. We claim that  $\left(z_{1}^{1},z_{2}^{1}\right) \in M S_{\overline{F}}(X, T) $. In fact, for any  $\tau>0 $, there is  $l \in \mathbb{N} $ such that  $\overline{A_{l}^{1}} \cap A \subset  V_{1} \times V_{2} $, where  $V_{i}=B\left(z_{i}^{1}, \tau\right)$  for  $i=1,2 $. By the construction, for any $ m \in \mathbb{N}$ , there are  $x_{m}^{1}, x_{m}^{2} \in U_{m}$  such that
$$\limsup _{n \rightarrow+\infty} \frac{1}{n} \inf _{\sigma \in S_{n}}\#\{1\leq i \leq n:\left(T^{i} x_m^1, T^{\sigma(i)} x_m^2\right) \in  \overline{A_{l}^{1}} \cap A \}> 0$$
and so
$$\limsup _{n \rightarrow+\infty} \frac{1}{n} \inf _{\sigma \in S_{n}}\#\{1\leq i \leq n:\left(T^{i} x_m^1, T^{\sigma(i)} x_m^2\right) \in  V_{1} \times V_{2} \}> 0$$
for all  $m \in \mathbb{N} $. For any nonempty open set  $U \subset X $, since $x$  is a transitive point, there is  $s \in \mathbb{Z}$  such that  $T^{s} x \in U$. We can choose  $m \in \mathbb{Z}$  such that  $T^{s} U_{m} \subset U $. This implies that  $T^{s} x_{m}^{1}, T^{s} x_{m}^{2} \in U$  and
$$\limsup _{n \rightarrow+\infty} \frac{1}{n} \inf _{\sigma \in S_{n}}\#\{1\leq i \leq n:\left(T^{i}\left(T^{s} x_{m}^{1}\right), T^{\sigma(i)} \left(T^{s} x_{m}^{2}\right)\right) \in  V_{1} \times V_{2} \}>0.$$
Hence, we have  $\left(z_{1}^{1}, z_{2}^{1}\right) \in M S_{\overline{F}}(X, T) $.
Similarly, for each  $p \in \mathbb{N} $, there exists  $\left(z_{1}^{p},z_{2}^{p}\right) \in M S_{\overline{F}}(X, T) \cap \prod_{i=1}^{2}   \pi^{-1}\left(\overline{B\left(y_{i}, 1 / p\right)}\right) $. Set  $z_{i}^{p} \rightarrow z_{i} $ as  $p \rightarrow \infty $. Then  $\left(z_{1}, z_{2}\right) \in M S_{\overline{F}}(X, T) \cup \Delta(X) $ and  $\pi\left(z_{i}\right)=y_{i}$ for $i=1,2$.

\section{  weakly mean equicontinuity and weakly equicontinuity in the mean  with respect to a function}

Motivated by [11], in this section,  We introduce the notions of  weakly mean equicontinuity (resp. weakly equicontinuity in the mean) and strong mean sensitivity (resp. strong sensitivity in the mean) of a topological dynamical system respect to a given continuous function $f$. We show that in the case of a minimal topological dynamical system, weakly mean equicontinuity (resp. weakly equicontinuity in the mean)  respect to $f$ and strong mean sensitivity (resp. strong sensitivity in the mean) respect to $f$ are complementary notions. Finally, we show that a topological dynamical system is weakly mean equicontinuity iff the system is weakly mean equicontinuity with respect to every continuous function.

Let  $(X, T)$  be a  topological dynamical system and  $f \in C(X) $. For any  $x, y \in X,$ $n\in\mathbb{N}$, we define
 $$d_{f}(x, y)=\limsup _{n \rightarrow \infty} \inf _{\sigma \in S_{n}}\frac{1}{n} \sum_{i=1}^{n}\left|f\left(T^{i} x\right)-f\left(T^{\sigma(i)} y\right)\right|$$
 
  $$\underline{d}_{f}(x, y)=\liminf_{n \rightarrow \infty} \inf _{\sigma \in S_{n}}\frac{1}{n} \sum_{i=1}^{n}\left|f\left(T^{i} x\right)-f\left(T^{\sigma(i)} y\right)\right|$$
 $$d_{f}^{n}(x, y)=\inf _{\sigma \in S_{n}}\frac{1}{n} \sum_{i=1}^{n}\left|f\left(T^{i} x\right)-f\left(T^{\sigma(i)} y\right)\right|$$
 where  $S_{n}$  is the $n$-order permutation group.  
 
 {\bf Definition 8.1}  Let $(X, T)$  be  a topological dynamical system  and  $f \in C(X)$.
 
  We say  $x \in X $ is a  $f$-weakly mean equicontinuity point if for every  $\varepsilon>0$  there exists  $\delta>0$  such that if  $d(x, y) < \delta $ then  $d_{f}(x, y) < \varepsilon $. We say  $(X, T)$  is  $f$-weakly mean equicontinuous if all  $x \in X$  are  $f$-weakly mean equicontinuous points. In this case we also say that  $f$  is weakly mean equicontinuous function for  $(X, T)$; we denote the set of weakly mean equicontinuous functions by  $C_{\text {me}}$. We say  $(X, T)$  is  $f$-almost weakly mean equicontinuous if the set of  $f$-weakly mean equicontinuity points is residual.
 
 We say  $x \in X $ is a  $f$-weakly  equicontinuous in the mean point if for every  $\varepsilon>0$  there exists  $\delta>0$  such that if  $d(x, y) < \delta $ then  $d_{f}^{n}(x, y) < \varepsilon $ for every $n\in\mathbb{N}$. We say  $(X, T)$  is  $f$-weakly equicontinuous in the mean  if all  $x \in X$  are  $f$-weakly  equicontinuous in the mean  points. In this case we also say that  $f$  is weakly equicontinuous in the mean  function for  $(X, T)$; we denote the set of weakly equicontinuous in the mean  functions by  $C_{\text {em}}$. We say  $(X, T)$  is  $f$-almost weakly equicontinuous in the mean if the set of  $f$-weakly equicontinuous in the mean  points is residual.
 
{\bf Definition 8.2} 
We say  $(X, T)$  is  $f$-strong mean sensitive if there exists  $\delta>0$  such that for every nonempty open set  $U$  there exists $ x, y \in U$  such that $d_{f}(x, y)>\delta.$
In this case we also say that  $f$  is strong mean sensitive for  $(X, T)$, and we denote the set of all strong mean sensitive functions for  $(X, T) $ by  $C_{ms}$.

$(X, T)$  is  $f$-strong  sensitive in the mean if there exists  $\delta>0$  such that for every nonempty open set  $U$  there exists $n\in\mathbb{N}$ and  $ x^{n}, y^{n} \in U$  such that 
$d_{f}^{n}(x^{n}, y^{n})>\delta.$
In this case we also say that  $f$  is strong  sensitive in the mean for  $(X, T)$, and we denote the set of all strong mean sensitive functions for  $(X, T) $ by  $C_{sm}$.

With the same methods in [12], we can easily obtain the propositions as follows.

{\bf Proposition 8.3} Let  $(X, T)$  be a topological dynamical system and  $f \in C(X)$. Then

(1) For any sequences  $\left\{x_{k}\right\}_{k=0}^{n-1}$  and  $\left\{y_{k}\right\}_{k=0}^{n-1}$  of $ X$, we have
$$
\inf _{\sigma \in S_{n}} \sum_{k=0}^{n-1} |f\left(x_{k}\right)-f(y_{\sigma(k)})|=\inf _{\sigma \in S_{n}} \sum_{k=0}^{n-1} |f\left(y_{k}\right)-f(x_{\sigma(k)})|
$$
In particular, for any  $x, y \in X$, we have
$$
\inf _{\sigma \in S_{n}} \frac{1}{n} \sum_{k=0}^{n-1} \left|f\left(T^{k} x\right)-f\left(T^{\sigma(k)} y\right)\right|=\inf _{\sigma \in S_{n}} \frac{1}{n} \sum_{k=0}^{n-1} \left|f\left(T^{k} y\right)-f\left(T^{\sigma(k)} x\right)\right .
$$

(2) For any sequences  $\left\{x_{k}\right\}_{k=0}^{n-1},\left\{y_{k}\right\}_{k=0}^{n-1}$  and $\left\{z_{k}\right\}_{k=0}^{n-1}$  of  $X$, we have
$$
\inf _{\sigma \in S_{n}} \sum_{k=0}^{n-1} |f\left(x_{k}\right)-f(z_{\sigma(k)})| \leq \inf _{\sigma \in S_{n}} \sum_{k=0}^{n-1} |f\left(x_{k}\right)-f(y_{\sigma(k)})|+\inf _{\sigma \in S_{n}}\sum_{k=0}^{n-1} |f\left(y_{k}\right)-f(z_{\sigma(k)})|.$$
In particular, for any  $x, y, z \in X$, we have
$$
\begin{aligned}
\inf _{\sigma \in S_{n}} \frac{1}{n} \sum_{k=0}^{n-1} \left|f\left(T^{k} x\right)-f\left(T^{\sigma(k)} z\right)\right| \leq \inf _{\sigma \in S_{n}} \frac{1}{n} \sum_{k=0}^{n-1}  \left|f\left(T^{k} x\right)-f\left(T^{\sigma(k)} y\right)\right|+\inf _{\sigma \in S_{n}} \frac{1}{n} \sum_{k=0}^{n-1}\left|f\left(T^{k} y\right)-f\left(T^{\sigma(k)} z\right)\right|.
\end{aligned}
$$

(3) For any  $x, y \in X$, we have
$$
d_{f}(x, y)=d_{f}(y, x)
$$
and
$$
\underline{d}_{f}(x, y)=\underline{d}_{f}(y, x).
$$

(4) For any  $x, y, z \in X$, we have
$$
d_{f}(x, z) \leq d_{f}(x, y)+d_{f}(y, z)
$$
and
$$
\underline{d}_{f}(x, z) \leq \underline{d}_{f}(x, y)+d_{f}(y, z).
$$

{\bf Proposition 8.4} Let  $(X, T)$  be a topological dynamical system and  $f \in C(X)$. For any  $x, y \in X$  and  $r, s \in \mathbb{N}$, we have
$$
d_{f}(T^{r} x, T^{s} y)=d_{f}(x, y)
$$
and
$$
\underline{d}_{f}\left(T^{r} x, T^{s} y\right)=\underline{d}_{f}(x, y).
$$

Let  $(X, T)$  be a topological dynamical system and  $f \in C(X) $. Fix $ x, y \in X $. For every  $n \in \mathbb{N} $ and  $\sigma \in S_{n}$,  we define
$$
\Delta_{\sigma}^{f}(x, y, \delta)=\#\{1 \leq i \leq n: \left|f\left(T^{i} x\right)-f\left(T^{\sigma(i)} y\right)\right|>\delta\},
$$
where  $\#(\cdot) $ stands for the cardinality of a set, and
$$
d_{\widetilde{f}}(x, y)=\inf \left\{\varepsilon>0: \limsup _{n \rightarrow \infty} \frac{1}{n} \inf _{\sigma \in S_{n}} \Delta_{\sigma}^{f}(x, y, \varepsilon)<\varepsilon\right\} .
$$

{\bf Lemma 8.5} For a topological dynamical system  $(X, T)$ and  $f \in C(X)$, the pseudo-metrics  $d_{f}$  and  $d_{\widetilde{f}}$  are uniformly equivalent on  $X $.

{\bf Proof\  } Without loss of generality we assume that the diameter of  $X$  is 1. Fix  $x, y \in X$. For every  $\delta>0, n \in \mathbb{N}$  and  $\sigma \in S_{n} $, we have
$$
\delta \cdot \Delta_{\sigma}^{f}(x, y, \delta) \leq \sum_{i=1}^{n}\left|f\left(T^{i} x\right)-f\left(T^{\sigma(i)} y\right)\right| \leq \Delta_{\sigma}^{f}(x, y, \delta)+\delta \cdot\left(n-\Delta_{\sigma}^{f}(x, y, \delta)\right) .
$$
Put  $\Delta_{n}^{f}(x, y, \delta)=\inf_{\sigma \in S_{n}} \Delta_{\sigma}^{f}(x, y, \delta)$  and
$$
\bar{d}(\Delta^{f}(x, y, \delta))=\limsup _{n \rightarrow \infty} \frac{1}{n} \Delta_{n}^{f}(x, y, \delta) .
$$
Then for every  $n \in \mathbb{N}$  we have
$$
\delta \cdot \Delta_{n}^{f}(x, y, \delta) \leq \inf_{\sigma \in S_{n}} \sum_{i=1}^{n}\left|f\left(T^{i} x\right)-f\left(T^{\sigma(i)} y\right)\right| \leq \Delta_{n}^{f}(x, y, \delta)+\delta \cdot\left(n-\Delta_{n}^{f}(x, y, \delta)\right) .
$$
Furthermore, by  $n$  and passing with  $n$  to infinity we obtain
$$
\delta \cdot \bar{d}(\Delta^{f}(x, y, \delta)) \leq d_{f}(x, y) \leq(1-\delta) \bar{d}(\Delta^{f}(x, y, \delta))+\delta .
$$
Fix  $\varepsilon>0$  and assume that  $\delta>0$  satisfies  $\delta<\varepsilon^{2} $. If  $d_{f}(x, y)<\delta $, then we deduce that
$$\varepsilon \cdot \bar{d}(\Delta^{f}(x, y, \varepsilon)) \leq d_{f}(x, y)<\varepsilon^{2}.$$
Then  $\bar{d}(\Delta(x, y, \varepsilon))<\varepsilon $, and by the definition one has  $\widetilde{E}(x, y)<\varepsilon $. This implies that the identity map $id:  (X, d_{f}) \rightarrow(X, d_{\widetilde{f}})$  is uniformly continuous.

It remains to show that $id:  (X, d_{\widetilde{f}}) \rightarrow(X, d_{f})$  is uniformly continuous. Fix  $\varepsilon>0$  and take  $\delta>0$  such that  $\delta<\varepsilon / 2$. Assume  $d_{\widetilde{f}}(x, y)<\delta $, that is  $\bar{d}(\Delta^{f}(x, y, \delta))<\delta $, and then we get that $d_{f}(x, y)<\varepsilon $.

{\bf Remark 8.6 } By  Lemma 8.5  it is easy to find that the pseudo-metric  $d_{f}$  does not depend on the metric  $d$  on  $X$, that is if  $d^{\prime}$  is another compatible metric on  $X$, then the pseudo-metric  $d_{f}$  and  $d_{f}^{\prime}$  induced by  $d$  and  $d^{\prime}$  are uniformly equivalent. In additional, since any continuous function on  $X$  is bounded, we can use  Lemma 8.5 to characterize  $f$-weakly mean equicontinuity and  $f$-strong mean sensitivity by replacing  $d_{f}(x, y)$ with $d_{\widetilde{f}}(x, y)$  in the definition. 

Again using the idea that a continuous function on a compact metric set is uniformly continuous, it is not hard to see that if  $(X, T)$  is  $f$-weakly mean equicontinuous, then it is uniformly  $f$-weakly mean equicontinuous in the sense that for every  $\varepsilon>0$  there exists  $\delta>0$  such that if  $d(x, y) < \delta$  then $ d_{f}(x, y) < \varepsilon$.

{\bf Definition 8.7} Let  $(X, T)$  be a  topological dynamical system and $f \in C(X)$. We denote the set of  $f$-weakly mean equicontinuity points by  $E^{f} $ and we define
$$
E_{\varepsilon}^{f}:=\left\{x \in X: \exists \delta>0, \forall y, z \in B(x,\delta), \liminf_{n \rightarrow \infty} \frac{1}{n} \inf _{\sigma \in S_{n}} \#\{1 \leq i \leq n: \left|f\left(T^{i} y\right)-f\left(T^{\sigma(i)} z\right)\right|\leq\varepsilon\}\geq 1-\varepsilon\right\} .
$$

Note that  $E^{f}=\cap_{\varepsilon>0} E_{\varepsilon}^{f} $.

{\bf Lemma 8.8}  Let  $(X, T)$  be a  topological dynamical system and $f \in C(X)$. The sets  $E^{f}$, $E_{\varepsilon}^{f}$  are inversely invariant (i.e.  $T^{-1}(E^{f}) \subseteq E^{f}$, $T^{-1}(E_{\varepsilon}^{f}) \subseteq E_{\varepsilon}^{f} )$, and $E_{\varepsilon}^{f}$ is open.

{\bf Proof\  } Let  $\varepsilon>0 $, and  $x \in T^{-1} E_{\varepsilon}^{f} $. There exists  $\eta>0$  such that if  $d\left(T x, z\right) < \eta$  and  $d\left(T x, y\right) < \eta$  then
$$
\liminf_{n \rightarrow \infty} \frac{1}{n} \inf _{\sigma \in S_{n}} \#\{1 \leq i \leq n: \left|f\left(T^{i} y\right)-f\left(T^{\sigma(i)} z\right)\right|\leq\varepsilon\}\geq 1-\varepsilon.
$$
There also exists  $\delta>0$  such that if  $d(x, y) < \delta$  then  $d\left(Tx, Ty\right) < \eta $. So, if  $y, z \in B(x,\delta)$, then  $\liminf_{n \rightarrow \infty} \frac{1}{n} \inf _{\sigma \in S_{n}} \#\{1 \leq i \leq n: \left|f\left(T^{i} Ty\right)-f\left(T^{\sigma(i)} Tz\right)\right|\leq\varepsilon\}\geq 1-\varepsilon.$ 
Hence, We have that  $x \in E_{\varepsilon}^{f}$, and so  $E_{\varepsilon}^{f}$  is inversely invariant. It follows that $ E^{f}$  is also inversely invariant.

Let  $x \in E_{\varepsilon}^{f}$  and  $\delta>0$  satisfy the defining property of  $E_{\varepsilon}^{f} $. If  $d(x, w)<\delta / 2 $ then $ w \in E_{\varepsilon}^{f}$; indeed if  $y, z \in B(w,\delta / 2)$  then  $y, z \in B(x,\delta)$. So,  $E_{\varepsilon}^{f}$  is open. The  proof is completed.

It is not hard to see that  $f$-strong mean sensitive systems have no  $f $-weakly mean equicontinuity points. Next we will give the  dichotomies as follows.

{\bf Theorem 8.9} Let  $(X, T)$  be a  topological dynamical system and $f \in C(X)$. Then

(1) if  $(X, T)$ is transitive, then it is either $ f$-almost weakly mean equicontinuous or  $f$-strong mean sensitive;

(2) if  $(X, T)$  is minimal, then it is either  $f$-weakly mean equicontinuous or $ f$-strong mean sensitive.

{\bf Proof\  } First, we show that if  $(X, T) $ is a transitive system, then for every  $\varepsilon $,  $E_{\varepsilon}^{f}$  is either empty or dense. Assume  $E_{\varepsilon}^{f}$  is nonempty and not dense. Then  $U=X \backslash\overline{E_{\varepsilon}^{f}}$  is a nonempty open set. Since the system is transitive and $ E_{\varepsilon}^{f} $ is nonempty and open there exists  $k \in \mathbb{N}$  such that  $U \bigcap T^{-k}\left(E_{\varepsilon}^{f}\right)$  is nonempty. By Lemma 8.8 we have that  $U \bigcap T^{-k}\left(E_{\varepsilon}^{f}\right) \subset U \bigcap E_{\varepsilon}^{f}=\emptyset $, a contradiction.

If  $E_{\varepsilon}^{f}$  is nonempty for every  $\varepsilon>0$  then we have that  $E^{f}=\bigcap_{n \geq 1} E_{1 / n}^{f}$  is a residual set; hence the system is  $f$-almost weakly mean equicontinuous.

If there exists  $\varepsilon>0 $ such that  $E_{\varepsilon}^{f} $ is empty, then for any open set $ U$  there exist  $y, z \in U$  such that  $$\liminf_{n \rightarrow \infty} \frac{1}{n} \inf _{\sigma \in S_{n}} \#\{1 \leq i \leq n: \left|f\left(T^{i} y\right)-f\left(T^{\sigma(i)} z\right)\right|\leq\varepsilon\}< 1-\varepsilon,$$ 
this means that  $$\limsup_{n \rightarrow \infty} \frac{1}{n} \inf _{\sigma \in S_{n}} \#\{1 \leq i \leq n: \left|f\left(T^{i} y\right)-f\left(T^{\sigma(i)} z\right)\right|>\varepsilon\}\geq \varepsilon.$$ It follows that $ (X, T)$  is  $f$-strong mean sensitive.

Now suppose  $(X, T)$  is minimal and  $f$-almost weakly mean equicontinuous. For every  $x \in X $ and every  $\varepsilon>0 $ there exists  $k \in \mathbb{N}$  such that  $T^{k} x \in E_{\varepsilon}^{f}$. Since  $E_{\varepsilon}^{f} $ is inversely invariant, we have  $x \in E_{\varepsilon}^{f}.$ So,  $x \in E^{f}$, and then $(X, T)$ is $f$-weakly mean equicontinuous.

With the similar methods, we can obtain the dichotomies between $f$-weakly equicontinuous in the mean and $ f$-strong sensitive in the mean.

{\bf Theorem 8.10}  
Let  $(X, T)$  be a  topological dynamical system and $f \in C(X)$. Then

(1) if  $(X, T)$ is transitive, then it is either $ f$-almost weakly equicontinuous in the mean or  $f$-strong  sensitive in the mean;

(2) if  $(X, T)$  is minimal, then it is either  $f$-weakly equicontinuous in the mean or $ f$-strong sensitive in the mean.

{\bf Theorem 8.11} Let  $(X, T)$  be a minimal topological dynamical system and $f \in C(X)$. Then  $(X, T)$  is  $f$-weakly mean equicontinuous if and only if it is $f$-weakly equicontinuous in the mean.

{\bf Proof}  It is easy to show that $f$-weakly equicontinuity in the mean implies $f$-weakly mean equicontinuity.
Assume that  $(X, T)$  is $f$-weakly mean equicontinuous. For any  $\varepsilon>0$  there is  $\delta_{1}>0$  such that if  $d(x, y)<\delta_{1}$  then
$$
\limsup _{n \rightarrow \infty} \inf _{\sigma \in S_{n}}\frac{1}{n} \sum_{i=1}^{n} \left|f\left(T^{i} x\right)-f\left(T^{\sigma(i)} y\right)\right|<\frac{\varepsilon}{8}.
$$
Fix  $z \in X$. For each  $m\in \mathbb{N}$ , let
$$
A_{m}=\left\{x \in \overline{B\left(z, \delta_{1} / 2\right)}: \inf _{\sigma \in S_{n}}\frac{1}{n} \sum_{i=1}^{n} \left|f\left(T^{i} x\right)-f\left(T^{\sigma(i)} z\right)\right| \leq \frac{\varepsilon}{4}, n=m, m+1, \ldots\right\} .
$$
Then $A_{m}$  is closed and  $\overline{B\left(z, \delta_{1} / 2\right)}=\bigcup_{m=1}^{\infty} A_{m}$. By the Baire Category Theorem, there is  $m_{1} \in \mathbb{N}$  such that  $A_{m_{1}} $contains an open subset $ U $ of $ X$. By the minimality of $(X, T)$,  we have that there is  $m_{2} \in \mathbb{N}$  with  $\bigcup_{i=0}^{m_{2}-1} T^{-i} U=X $. Let  $\delta_{2}$  be the Lebesgue number of the open cover  $\left\{T^{-i} U: 0 \leq i \leq m_{2}-1\right\}$  of  $X $. Let  $m=\max \left\{m_{1}, 2m_{2}\right\} $. With the continuity of  $f$, there exists  $\delta_{3}>0$  such that if  $d(x, y)<\delta_{3}$  implies that $|f(x)-f(y)| < \varepsilon/4 $, and with the continuity of  $T$, there exists  $\delta_{4}>0$ such that if  $d(x, y)<\delta_{4}$ implies that  $d\left(T^{i} x, T^{i} y\right)<\delta_{3}$  for any  $0 \leq i \leq m $. Put  $\delta=\min \left\{\delta_{2}, \delta_{3},\delta_{4}\right\} $. Let  $x, y \in X$  with  $d(x, y)<\delta$  and  $n \in \mathbb{N} $. 

If  $n \leq m$, then
$$
\inf _{\sigma \in S_{n}}\frac{1}{n} \sum_{i=1}^{n} \left|f\left(T^{i} x\right)-f\left(T^{\sigma(i)} y\right)\right| \leq \frac{1}{n} \cdot n \cdot \frac{\varepsilon}{4}<\varepsilon.
$$

If  $n>m$, there exists  $0 \leq k \leq m_{2}-1 $ such that  $x, y \in T^{-k} U$, then  $T^{k} x, T^{k} y \in U$, one has
$$
\begin{aligned}
	&\inf _{\sigma \in S_{n}}\frac{1}{n} \sum_{i=1}^{n} \left|f\left(T^{i} x\right)-f\left(T^{\sigma(i)} y\right)\right|\\
	 &\leq \inf _{\sigma \in S_{k}}\frac{1}{n} \sum_{i=1}^{k} \left|f\left(T^{i} x\right)-f\left(T^{\sigma(i)} y\right)\right|+\inf _{\sigma \in S_{n}}\frac{1}{n} \sum_{i=1}^{n} \left|f\left(T^{i} T^{k}x\right)-f\left(T^{\sigma(i)} T^{k}y\right)\right| \\
	& \leq \frac{\varepsilon}{4}+\inf _{\sigma \in S_{n}}\frac{1}{n} \sum_{i=1}^{n} \left|f\left(T^{i} T^{k}x\right)-f\left(T^{\sigma(i)} z\right)\right|+\inf _{\sigma \in S_{n}}\frac{1}{n} \sum_{i=1}^{n} \left|f\left(T^{\sigma(i)} T^{k}y\right)-f\left(T^{i} z\right)\right| \\
	& = \frac{\varepsilon}{4}+\inf _{\sigma \in S_{n}}\frac{1}{n} \sum_{i=1}^{n} \left|f\left(T^{i} T^{k}x\right)-f\left(T^{\sigma(i)} z\right)\right|+\inf _{\sigma \in S_{n}}\frac{1}{n} \sum_{i=1}^{n} \left|f\left(T^{i} T^{k}y\right)-f\left(T^{\sigma(i)} z\right)\right|\\
	& \leq \frac{\varepsilon}{4}+\frac{\varepsilon}{4}+\frac{\varepsilon}{4}<\varepsilon.
\end{aligned}
$$
Therefore,  $\sup _{n \in \mathbb{N}} \inf _{\sigma \in S_{n}}\frac{1}{n} \sum_{i=1}^{n} \left|f\left(T^{i} x\right)-f\left(T^{\sigma(i)} y\right)\right|<\varepsilon$. This implies that  $(X, T)$  is $f$-weakly equicontinuous in the mean.

{\bf Theorem 8.12} Let  $(X, T)$  be a minimal topological dynamical system and $f \in C(X)$. Then  $(X, T)$  is  $f$-strong mean sensitive if and only if it is $f$-strong sensitive in the mean.

By Theorem 8.12, we have that for a minimal system, the equivalence between  $f$-strong mean sensitivity and  $f$-strong sensitivity in the mean holds,  while we have found that above conclusion remains true even for general dynamical system (see Appendix Theorem 1 and Theorem 3).

{\bf Theorem 8.13} Let  $(X, T)$  be a  topological dynamical system. If  $(X, T)$  is weakly mean equicontinuous, then it is  $f$-weakly mean  equicontinuous for every  $f \in C(X)$ (i.e. $C(X)=C_{m e}$).

{\bf Proof\  }    
Let $ (X, T)$  be weakly mean equicontinuous,  $f \in C(X) $ and  $\varepsilon>0 $. Since  $X$  is compact,  $f$  is uniformly continuous; thus there exists $ \delta^{\prime} \in(0, \varepsilon) $ such that if  $d(x, y) < \delta^{\prime}$  then  $|f(x)-f(y)| < \varepsilon $. Since  $(X, T)$  is weakly mean equicontinuous, there exists $ \delta>0 $ such that if  $d(x, y) <\delta$  then
$$
 \liminf_{n \rightarrow \infty} \frac{1}{n} \inf_{\sigma \in S_{n}} \#\{1 \leq i \leq n: d(T^{i}x,T^{\sigma(i)}y)\leq\delta^{\prime}\}\geq 1-\delta^{\prime}.
$$

This implies that
$$
\liminf_{n \rightarrow \infty} \frac{1}{n} \inf _{\sigma \in S_{n}} \#\{1 \leq i \leq n: \left|f\left(T^{i} x\right)-f\left(T^{\sigma(i)} y\right)\right|\leq\varepsilon\}\geq 1-\varepsilon.
$$
Hence $ (X, T)$  is  $f$-weakly mean equicontinuous.

\section{Appendix}
Recall the notions of mean sensitivity and  sensitivity in the mean.

A topological dynamical system  $(X, T)$  is mean sensitive if there exists  $\delta>0$  such that  for any nonempty open subset $U$ of $X$, there exist $x, y \in U$ satisfying
$$
\limsup _{n \rightarrow+\infty} \frac{1}{n} \sum_{k=1}^{n} d\left(T^{k} x, T^{k} y\right)>\delta.$$

A topological dynamical system  $(X, T)$  is sensitive in the mean if  there exists  $\delta>0$  such that  for any nonempty open subset $U$ of $X$,
there exist   $n\in\mathbb{N}$ and $x^{n}, y^{n} \in U$ satisfying
$$
\frac{1}{n} \sum_{k=1}^{n} d\left(T^{k} x^{n}, T^{k} y^{n}\right)>\delta.$$

{\bf Theorem 1} Let  $(X, T)$  be a topological dynamical system. Then  $(X, T)$  is  mean sensitive if and only if sensitive in the mean.

{\bf Proof\  }  Necessity is clear, we only need to show the sufficiency.  
 Assume that  $(X, T) $ is sensitive in the mean and  $\delta$  is the constant of sensitivity in the mean. Given any  $m \in \mathbb{N} $, denote  $$F_{m}=\left\{(x, y): \frac{1}{n} \sum_{k=1}^{n} d\left(T^{k} x, T^{k} y\right) \leq \frac{\delta }{3}, \forall n>m\right\} .$$ It is clear that  $F_{m}$  is a closed set.

{\bf Claim } for any  $m \in \mathbb{N}$, $ (F_{m})^{o}=\emptyset $, where $(F_{m})^{o}$ is the interior of the set $F_{m}$.

{\bf Proof of Claim\  } 
 If not, there exists some  $m \in \mathbb{N}$  such that $(F_{m})^{o} \neq \emptyset $, there exist nonempty open sets  $U, V \subset X$  such that  $U \times V \subset F_{m}$. This implies that for any   $(x, y) \in U \times V$, then  $\frac{1}{n} \sum_{k=1}^{n} d\left(T^{k} x, T^{k} y\right) \leq \frac{\delta }{3}$  holds for any  $n>m$. It follows that  for any points  $x_{1}, x_{2} \in U$  and any $ n>m,$
 
 $$\frac{1}{n} \sum_{k=1}^{n} d\left(T^{k} x_{1}, T^{k} x_{2}\right) \leq \frac{1}{n} \sum_{k=1}^{n} d\left(T^{k} x_{1}, T^{k}y\right)+\frac{1}{n} \sum_{k=1}^{n} d\left(T^{k} x_{2}, T^{k} y\right) \leq \frac{2\delta }{3}.$$
 
With the continuity of $T$, there exists $\tilde{\delta}>0$, a non-empty open set  $\tilde{U} \subset U$ with $diam(\tilde{U})<\tilde{\delta}$ such that for any   $x_{1}, x_{2} \in X$  with $d(x_{1}, x_{2})<\tilde{\delta}$  and any  $0 \leq n \leq m$, $\frac{1}{n} \sum_{k=1}^{n} d\left(T^{k} x_{1}, T^{k} x_{2}\right) \leq \frac{2\delta }{3}$. Then  for any  $x_{1}, x_{2} \in \tilde{U}$  and any  $n \in \mathbb{N},$  $\frac{1}{n} \sum_{k=1}^{n} d\left(T^{k} x_{1}, T^{k} x_{2}\right) \leq \frac{2\delta }{3}$,  which is a contradiction with  the sensitivity in the mean of  $(X, T)$.

Moreover, we can deduce that  the set  $F=\bigcup_{m \in \mathrm{N}} F_{m}$ is a first category set in  $X \times X $. This implies that the set $ F^{C}=\{(x, y): \forall m \in \mathbb{N}, \exists n>m$  such that  $\frac{1}{n} \sum_{k=1}^{n} d\left(T^{k} x, T^{k} y\right)>\frac{\delta }{3}\}$  is residual in  $X \times X $, and then $ F^{C}$ is dense on $ X \times X $.

Suppose that  $(X, T)$  is  not mean sensitive,  there exists nonempty open subset $U$ of $X$, for any $x, y \in U$ satisfying
$$
\limsup _{n \rightarrow+\infty} \frac{1}{n} \sum_{k=1}^{n} d\left(T^{k} x, T^{k} y\right)\leq \frac{\delta }{4}.$$

Since  $F^{C}$  is residual in  $X \times X $, then we have that  there exists a pair  $\left(y_{1}, y_{2}\right) \in (U \times U )\cap F^{C} $, and  one has 
$$\limsup _{n \rightarrow+\infty} \frac{1}{n} \sum_{k=1}^{n} d\left(T^{k} y_{1}, T^{k}y_{2}\right) \leq  \frac{\delta }{4}.$$
 which is a  contradiction, for   $\left(y_{1}, y_{2}\right) \in F^{C}$.
Hence  $(X, T)$  is  mean sensitive. The proof is completed.

{\bf Remark } In [8], the authors show that  for a minimal topological dynamical system, then  $(X, T)$  is  mean sensitive if and only if sensitive in the mean, Theorem 1 extend their conclusion, and in their paper, they leave a question as follows:  for a minimal system, does the equivalence between mean $n$-sensitivity and $n$-sensitivity in the mean still hold ? While with the above Theorem 1, we may conjecture that  a  topological dynamical system is  mean $n$-sensitive if and only if $n$-sensitive in the mean.

With the similar method, we can easily obtain the equivalence as follows:

{\bf Theorem 2} Let  $(X, T)$  be a topological dynamical system. Then  $(X, T)$  is  strong mean sensitive if and only if it is strong sensitive in the mean.

{\bf Theorem 3} Let  $(X, T)$  be a topological dynamical system and $f \in C(X)$. Then  $(X, T)$  is  $f$-strong mean sensitive if and only if it is  $f$-strong sensitive in the mean.

{\bf Theorem 4}\ \  Let  $(X, T)$  be a transitive topological dynamical system Then  $SM_{\overline{F}}(X,T)\neq\emptyset$ if and only if $MS_{\overline{F}}(X, T)\neq\emptyset$.

{\bf Proof \  } It is easy to obtain it, since a transitive dynamical system is strong mean sensitive  if and only if it admits a strong mean sensitive tuple;  a transitive dynamical system is strong sensitive in the meae  if and only if it admits a strong sensitive in the mean tuple.

{\bf Theorem 5}\ \  Let  $(X, T)$  be a transitive topological dynamical system. Then  $(X, T)$  is  almost weakly mean equicontinuous if and only if it is  almost weakly equicontinuous in the mean.

{\bf Theorem 6}\ \  Let  $(X, T)$  be a transitive topological dynamical system and $f \in C(X)$. Then  $(X, T)$  is  almost $f$-weakly mean equicontinuous if and only if it is  almost $f$-weakly equicontinuous in the mean.

{\bf Theorem 7}\ \  Let  $(X, T)$  be a transitive topological dynamical system. Then  $(X, T)$  is  almost mean equicontinuous if and only if it is  almost  equicontinuous in the mean.

\bigskip \noindent{\bf Acknowledgement}.
The research was supported by NSF of China (No. 11671057) and NSF of Chongqing (Grant No. cstc2020jcyj-msxmX0694).

%¸ñÊ½´ýµ÷Õû

 \bibliographystyle{amsplain}

\end{document}